\documentclass[10pt,a4paper]{article}
\usepackage{fontspec}
\usepackage{etoolbox}

\usepackage{setspace}
\usepackage{array}
\usepackage{multirow}
\usepackage{booktabs} 
\usepackage{makecell} 
\usepackage[head=28pt]{geometry}
\usepackage[dvipsnames,table]{xcolor}
\usepackage{csquotes}      

\usepackage[
  backend=biber,
  style=authoryear,
  sortcites = true,
  dashed=false,
  doi=true,
  giveninits=true,
  uniquename=init,
  url=false,
  eprint=false,
  maxcitenames=99,
  mincitenames=1,
  maxbibnames=99
]{biblatex}
\AtEveryBibitem{
  \clearfield{month}
}
\DeclareFieldFormat[article]{volume}{\bibstring{volume}\addnbspace#1}
\DeclareFieldFormat[article]{number}{\addcomma\space\bibstring{number}\addnbspace#1}

\usepackage{hyperref}

\usepackage{amsmath}
\usepackage{amssymb}
\usepackage{stmaryrd}
\SetSymbolFont{stmry}{bold}{U}{stmry}{m}{n}
\usepackage{mathrsfs}
\usepackage{mathtools}
\usepackage{dsfont}
\usepackage{bm}

\usepackage{cleveref}
\usepackage{scrhack}
\usepackage{adjustbox}

\usepackage{amsthm}
\usepackage{thmtools}
\usepackage{mdframed}
\usepackage{enumitem}
\usepackage{microtype}

\usepackage{graphicx}
\usepackage{subcaption}
\usepackage{mwe}
\usepackage{siunitx}
\usepackage[labelfont=bf]{caption} 
\usepackage{float}

\usepackage{tikz}
\usetikzlibrary{calc,positioning,shapes,arrows.meta}

\usepackage{pgfornament}

\usepackage[absolute,overlay]{textpos}

\usepackage{epigraph}

\usepackage{wrapfig}

\usepackage{thmtools}
\usepackage{mdframed}
\definecolor{cardinalred}{RGB}{140,21,21}
\definecolor{midnightblue}{RGB}{25,25,112}
\definecolor{forestgreen}{RGB}{34,94,34}
\definecolor{goldenrod}{RGB}{218,165,32}

\definecolor{goalblue}{RGB}{30,70,140}
\definecolor{goalgray}{RGB}{248,249,250}

\colorlet{goalbg}{goalblue!5}
\colorlet{questionbg}{cardinalred!5}

\declaretheoremstyle[
  headfont=\color{cardinalred}\scshape,
  headformat=\NAME\ \textcolor{gray}{\NUMBER}%
    \ifstrempty{\NOTE}{}{\textcolor{cardinalred}{\NOTE}},
  postheadspace=1em,
  qed=\raisebox{1pt}{\textcolor{cardinalred}{\openbox}}
]{theoremstyle}

\declaretheoremstyle[
  headfont=\color{midnightblue}\scshape,
  headformat=\NAME\ \textcolor{gray}{\NUMBER}%
  \ifstrempty{\NOTE}{}{\textcolor{midnightblue}{\NOTE}},
  postheadspace=1em,
  qed=\raisebox{1pt}{\textcolor{midnightblue}{\openbox}}
]{definitionstyle}

\declaretheoremstyle[
  headfont=\color{black}\scshape,
  headformat=\NAME\ \textcolor{gray}{\NUMBER}%
  \ifstrempty{\NOTE}{}{\textcolor{black}{\NOTE}},
  postheadspace=1em,
  qed=\raisebox{1pt}{\textcolor{black}{\openbox}}
]{assumptionstyle}

\declaretheoremstyle[
  headfont=\color{forestgreen}\scshape,
  headformat=\NAME\ \textcolor{gray}{\NUMBER}%
    \ifstrempty{\NOTE}{}{\textcolor{forestgreen}{\NOTE}},
  postheadspace=1em,
  qed=\raisebox{1pt}{\textcolor{forestgreen}{\openbox}}
]{examplestyle}

\declaretheoremstyle[
  headfont=\color{goldenrod}\scshape,
  headformat=\NAME\ \textcolor{gray}{\NUMBER}%
    \ifstrempty{\NOTE}{}{\textcolor{goldenrod}{\NOTE}},
  postheadspace=1em,
  qed=\raisebox{1pt}{\textcolor{goldenrod}{\openbox}}
]{notationstyle}

\declaretheoremstyle[
  headfont=\color{brown!70!black}\scshape,
  headformat=\NAME\ \textcolor{gray}{\NUMBER}%
    \ifstrempty{\NOTE}{}{\textcolor{brown!70!black}{\NOTE}},
  postheadspace=1em,
  qed=\raisebox{1pt}{\textcolor{brown!70!black}{\openbox}}
]{remarkstyle}

\declaretheorem[style=theoremstyle, numberwithin=section]{theorem}
\declaretheorem[style=theoremstyle, sibling=theorem]{lemma}
\declaretheorem[style=theoremstyle, sibling=theorem]{corollary}
\declaretheorem[style=theoremstyle, sibling=theorem]{proposition}
\declaretheorem[style=definitionstyle, sibling=theorem]{definition}
\declaretheorem[
  style=assumptionstyle,
  sibling=theorem,
  name=Assumption
]{assumption}

\declaretheorem[style=notationstyle, sibling=theorem]{notation}
\declaretheorem[style=remarkstyle, sibling=theorem]{remark}
\declaretheorem[style=examplestyle, sibling=theorem]{example}
\declaretheorem[style=examplestyle, sibling=theorem]{illustration}

\crefname{theorem}{\textcolor{cardinalred}{\textsc{Theorem}}}{\textsc{Theorems}}
\crefname{lemma}{\textcolor{cardinalred}{\textsc{Lemma}}}{\textcolor{cardinalred}{\textsc{Lemmas}}}
\crefname{proposition}{\textcolor{cardinalred}{\textsc{Proposition}}}{\textsc{Propositions}}
\crefname{corollary}{\textcolor{cardinalred}{\textsc{Corollary}}}{\textsc{Corollaries}}
\crefname{definition}{\textcolor{midnightblue}{\textsc{Definition}}}{\textcolor{midnightblue}{\textsc{Definitions}}}

\crefname{remark}{\textcolor{brown!70!black}{\textsc{Remark}}}{\textcolor{brown!70!black}{\textsc{Remarks}}}

\crefname{example}{\textcolor{forestgreen}{\textsc{Example}}}{\textcolor{forestgreen}{\textsc{Examples}}}
\crefname{illustration}{\textcolor{forestgreen}{\textsc{Illustration}}}{\textsc{Illustrations}}
\crefname{assumption}{\textsc{Assumption}}{\textsc{Assumptions}}
\crefname{condition}{\textsc{Assumption}}{\textsc{Assumptions}}
\crefname{notation}{\textcolor{goldenrod}{\textsc{Notation}}}{\textsc{Notations}}

\creflabelformat{theorem}{\textcolor{gray}{#2#1#3}}
\creflabelformat{lemma}{\textcolor{gray}{#2#1#3}}
\creflabelformat{proposition}{\textcolor{gray}{#2#1#3}}
\creflabelformat{corollary}{\textcolor{gray}{#2#1#3}}
\creflabelformat{definition}{\textcolor{gray}{#2#1#3}}
\creflabelformat{remark}{\textcolor{gray}{#2#1#3}}
\creflabelformat{example}{\textcolor{gray}{#2#1#3}}
\creflabelformat{illustration}{\textcolor{gray}{#2#1#3}}
\creflabelformat{assumption}{\textcolor{gray}{#2#1#3}}
\creflabelformat{condition}{\textcolor{gray}{#2#1#3}}
\creflabelformat{notation}{\textcolor{gray}{#2#1#3}}

\crefformat{section}{\textcolor{cardinalred}{\textsc{Section}}~#2\textcolor{gray}{#1}#3}
\Crefformat{section}{\textcolor{cardinalred}{\textsc{Section}}~#2\textcolor{gray}{#1}#3}

\crefformat{chapter}{\textcolor{cardinalred}{\textsc{Chapter}}~#2\textcolor{gray}{#1}#3}
\Crefformat{chapter}{\textcolor{cardinalred}{\textsc{Chapter}}~#2\textcolor{gray}{#1}#3}

\crefformat{part}{\textcolor{cardinalred}{\textsc{Part}}~#2\textcolor{gray}{#1}#3}
\Crefformat{part}{\textcolor{cardinalred}{\textsc{Part}}~#2\textcolor{gray}{#1}#3}

\crefname{figure}{\textcolor{cardinalred}{\textsc{Fig.}}}{\textcolor{cardinalred}{\textsc{Figs.}}}
\Crefname{figure}{\textcolor{cardinalred}{\textsc{Fig.}}}{\textcolor{cardinalred}{\textsc{Figs.}}}

\crefformat{figure}{%
  \textcolor{cardinalred}{\textsc{Fig.}}~#2\textcolor{gray}{#1}#3%
}
\Crefformat{figure}{%
  \textcolor{cardinalred}{\textsc{Fig.}}~#2\textcolor{gray}{#1}#3%
}

\crefname{table}{\textcolor{cardinalred}{\textsc{Tab.}}}{\textcolor{cardinalred}{\textsc{Tabs.}}}
\Crefname{table}{\textcolor{cardinalred}{\textsc{Tab.}}}{\textcolor{cardinalred}{\textsc{Tabs.}}}

\crefformat{table}{%
  \textcolor{cardinalred}{\textsc{Tab.}}~#2\textcolor{gray}{#1}#3%
}
\Crefformat{table}{%
  \textcolor{cardinalred}{\textsc{Tab.}}~#2\textcolor{gray}{#1}#3%
}

\crefname{equation}{\textsc{Eq.}}{\textsc{Eqs.}}
\Crefname{equation}{\textsc{Equation}}{\textsc{Equations}}

\crefname{enumi}{}{}
\Crefname{enumi}{}{}

\definecolor{mycitecolor}{RGB}{25,25,112}

\DeclareNameFormat{labelname}{%
  \nameparts{#1}%
  \ifnumgreater{\value{listtotal}}{2}
    {%
      \ifnumequal{\value{listcount}}{1}
        {%
          \usebibmacro{name:family}
            {\namepartfamily}
            {\namepartgiven}
            {\namepartprefix}
            {\namepartsuffix}%
          \addspace\bibstring{andothers}%
        }
        {}%
    }
    {%
      \usebibmacro{name:family}
        {\namepartfamily}
        {\namepartgiven}
        {\namepartprefix}
        {\namepartsuffix}%
    }%
}

\DeclareDelimFormat[cite]{finalnamedelim}{\addspace and\space}
\DeclareDelimFormat[cite]{multinamedelim}{\addcomma\space}

\DefineBibliographyStrings{english}{andothers = {et\addabbrvspace al\adddot}}
\DefineBibliographyStrings{ngerman}{andothers = {et\addabbrvspace al\adddot}}
\DefineBibliographyStrings{german}{andothers = {et\addabbrvspace al\adddot}}

\DeclareCiteCommand{\cite}[\mkbibparens]
  {\usebibmacro{prenote}}
  {%
  \usebibmacro{citeindex}%
  \printtext[bibhyperref]{%
    \textcolor{mycitecolor}{\printnames{labelname}}%
    \addcomma\space
    \textcolor{mycitecolor}{\printlabeldateextra}}}
{\multicitedelim}
  {\setunit{\addcomma\space}\printfield{postnote}}

\renewcommand*{\multicitedelim}{%
  \ifnumequal{\value{citecount}}{\value{citetotal}}
    {%
      \ifnumequal{\value{citetotal}}{2}
        {\space and\space}
        {\addcomma\space and\space}%
    }
    {\addsemicolon\space}%
}

\newcommand{\IC}{\mathbb{C}}

\newcommand{\IE}{\mathbb{E}}

\newcommand{\IL}{L}

\newcommand{\IN}{\mathbb{N}}

\newcommand{\IP}{\mathbb{P}}

\newcommand{\IR}{\mathbb{R}}

\newcommand{\IWb}{\rmZ^{\rmw,\operatorname{b}}}
\newcommand{\IWu}{\rmZ^{\rmw,\operatorname{u}}}

\newcommand{\rmD}{\mathrm{D}}

\newcommand{\rmL}{\Lambda}
\newcommand{\rmM}{\mathrm{M}}

\newcommand{\rmO}{\mathrm{O}}

\newcommand{\rmR}{\mathrm{R}}

\newcommand{\rmV}{\mathrm{V}}
\newcommand{\rmVhat}{\widehat{\mathrm{V}}}
\newcommand{\rmW}{\mathrm{W}}

\newcommand{\rmZ}{\mathrm{Z}}
\newcommand{\rmZcheck}{\mathrm{Z}}
\newcommand{\rmZbar}{\mathrm{Z}}

\renewcommand{\mho}{\widehat{\Omega}}

\newcommand{\rma}{\mathrm{a}}
\newcommand{\rmb}{\mathrm{b}}
\newcommand{\rmbos}{\mathrm{b}_{\textsc{os}}}
\newcommand{\rmbes}{\mathrm{b}_{\textsc{es}}}
\newcommand{\rmc}{\mathrm{c}}
\newcommand{\rmd}{\mathrm{d}}
\newcommand{\rme}{\mathrm{e}}

\newcommand{\rmp}{\mathrm{q}}
\newcommand{\rmq}{\mathrm{q}}
\newcommand{\rmr}{\mathrm{r}}
\newcommand{\rms}{\mathrm{s}}
\newcommand{\rmsos}{\mathrm{s}_{\textsc{os}}}
\newcommand{\rmses}{\mathrm{s}_{\textsc{es}}}

\newcommand{\rmv}{\mathrm{v}}
\newcommand{\rmw}{\mathrm{w}}
\newcommand{\rmx}{\mathrm{x}}

\newcommand{\calA}{\mathcal{A}}

\newcommand{\calK}{\mathcal{K}}
\newcommand{\calL}{\mathcal{L}}

\newcommand{\scrB}{\mathscr{B}}

\newcommand{\frake}{\mathfrak{e}}

\newcommand{\frakm}{\mathfrak{m}}

\newcommand{\frakI}{\mathfrak{I}}

\newcommand{\frakMl}{\mathcal{T}_{M}}
\newcommand{\indMl}{\1_{\mathcal{T}_{M}}}

\newcommand{\frakR}{\mathfrak{R}}

\newcommand{\br}[1]{\left(#1\right)}
\newcommand{\brK}[1]{\left\{#1\right\}}
\newcommand{\brE}[1]{\left[#1\right]}

\makeatletter
\newcommand{\norm}[1]{%
  \left\| #1 \right\|%
  \@ifnextchar\bgroup 
    {\norm@i}
    {}
}

\newcommand{\norm@i}[1]{_{#1}}

\makeatother

\newcommand{\pRz}{\IR_{>0}}
\newcommand{\pRzz}{\mathbb{R}_{\geq 0}}

\newcommand{\nset}[1]{\left\llbracket #1\right\rrbracket}

\newcommand{\1}{\mathds{1}}

\newcommand{\sigmabar}{\widetilde{\sigma}}

\newcommand{\ConstLlower}{\Const_{\rmL,\mathrm{Lower}}}
\newcommand{\ConstLupper}{\Const_{\rmL,\mathrm{Upper}}}

\newcommand{\arsinh}{\operatorname{arsinh}}

\DeclareMathOperator*{\argmin}{arg\,min}

\newcommand{\E}{\IE}

\newcommand{\abs}[1]{\left|#1\right|}

\newcommand{\ind}{\1}
\newcommand{\indk}[1]{\ind_{[-#1,#1]}}%
\newcommand{\indkco}[1]{\ind_{[-#1,#1]^C}}%
\makeatletter
\newcommand{\Const}{%
  \mathfrak{C}%
  \@ifnextchar\bgroup 
    {\Const@i}
    {}
}

\newcommand{\Const@i}[1]{_{\mathrm{#1}}}

\makeatletter

\newcommand{\const}{%
  \@ifnextchar\bgroup
    {\const@with} 
    {\const@without} 
}

\newcommand{\const@with}[1]{%
  \mathfrak{C}_{\mathrm{#1}}%
}

\newcommand{\const@without}{%
  \mathfrak{c}%
}

\makeatother

\newcommand{\Constfeps}{\Const_{f^\varepsilon}}
\newcommand{\ConstfY}{\Const_{f^Y}}
\makeatother

\newcommand{\Var}{\operatorname{Var}}

\newcommand{\sqBias}[1]{\operatorname{bias}_{#1}^2}

\newcommand{\maxrisklocal}{\mathrm{locR}_{\textsc{max}}^{N}}

\newcommand{\ddrisklocal}{\mathrm{locR}_{\textsc{dd}}^N}

\newcommand{\sizedMid}[2]{#1 \, \kern-\nulldelimiterspace\mathopen{}\left| \vphantom{#1}\,#2\right.\mathclose{}\kern-\nulldelimiterspace}

\newcommand{\assuRef}[1]{\texorpdfstring{\protect\hyperlink{assu#1}{\textbf{(\small{#1})}}}{}}
\newcommand{\newAssuRef}[1]{\hypertarget{assu#1}{\textbf{(\small{#1})}:}}
\newcommand{\lebesgue}{\calL} 
\newcommand{\poslebesgue}{\calL_{>0}} 
\renewcommand{\Im}[1]{\frakI\frakm\br{#1}}
\renewcommand{\Re}[1]{\frakR\frake\br{#1}}

\newcommand{\kdiamond}{k^\diamond}

\usepackage{xspace}

\newcommand{\M}{\rmM}

\newcommand{\MelU}{\M_{\varepsilon}}
\newcommand{\MelUrez}{\M_{\varepsilon}^\dagger}
\newcommand{\MelUhat}{\widehat{\M}_{\varepsilon}}
\newcommand{\MelUhatrez}{\widehat{\M}_{\varepsilon}^\dagger}

\newcommand{\MelX}{\M_{X}}
\newcommand{\MelXcheck}{\widecheck{\M}_{X}}
\newcommand{\MelXhat}{\widehat{\M}_{X}}
\newcommand{\MelY}{\M_{Y}}

\newcommand{\MelYb}{\M_{Y}^{\operatorname{b}}}
\newcommand{\MelYu}{\M_{Y}^{\operatorname{u}}}

\newcommand{\MelYhat}{\widehat{\M}_{Y}}
\newcommand{\MelYhatb}{\widehat{\M}_{Y}^{\operatorname{b}}}
\newcommand{\MelYhatu}{\widehat{\M}_{Y}^{\operatorname{u}}}

\newcommand{\SXhat}{ \widehat{S}^{X}}
\newcommand{\SX}{ S^{X}}
\newcommand{\SXbar}{ \widetilde{S}^{X}}

\newcommand{\rmwbar}{\bar{\mathrm{w}}}
\newcommand{\rmwrez}{\rmw^\dagger}

\newcommand{\rmwU}{\mathrm{w}_\varepsilon}

\newcommand{\rmwhatU}{\widehat{\mathrm{w}}_\varepsilon}

\newcommand{\khat}{\hat{k}}

\newcommand{\Rz}{\mathbb R}

\makeatletter
\newcommand{\ocon}{\mathbin{\mathpalette\make@circled *}}
\newcommand{\make@circled}[2]{%
	\ooalign{$\m@th#1\smallbigcirc{#1}$\cr\hidewidth$\m@th#1#2$\hidewidth\cr}%
}
\newcommand{\smallbigcirc}[1]{%
	\vcenter{\hbox{\scalebox{0.77778}{$\m@th#1\bigcirc$}}}%
}
\makeatother

\renewcommand{\and}{\text{and}}

\newcommand{\penhat}{\mathop{\mathrm{P}\mkern-1mu \widehat{\mathrm{en}}}\nolimits}
\newcommand{\pen}{\operatorname{pen}}
\newcommand{\penbar}{\mathop{\mathrm{Pen}}\nolimits}

\DeclareFontFamily{U}{mathx}{}
\DeclareFontShape{U}{mathx}{m}{n}{<-> mathx10}{}
\DeclareSymbolFont{mathx}{U}{mathx}{m}{n}
\DeclareMathAccent{\widehat}{0}{mathx}{"70}
\DeclareMathAccent{\widecheck}{0}{mathx}{"71}

\newcommand{\ko}[1][\circ]{k_{#1}}
\newcommand{\Mo}[1][\circ]{M_{#1}}

\newcommand{\MelXtilde}{\widetilde{\M}_{X}}
\newcommand{\skalar}[1]{\big\langle#1\big\rangle}

\newcommand*{\ob}[1]{%
	\makebox[0pt][l]{%
		\kern.3em
		\raisebox{1.8ex}{\relsize{-3}{$\circ$}}
	}%
	#1%
}

\makeatletter
\newsavebox\myboxA
\newsavebox\myboxB
\newlength\mylenA

\newcommand*\xoverline[2][0.75]{%
	\sbox{\myboxA}{$\m@th#2$}%
	\setbox\myboxB\null
	\ht\myboxB=\ht\myboxA%
	\dp\myboxB=\dp\myboxA%
	\wd\myboxB=#1\wd\myboxA
	\sbox\myboxB{$\m@th\overline{\copy\myboxB}$}
	\setlength\mylenA{\the\wd\myboxA}
	\addtolength\mylenA{-\the\wd\myboxB}%
	\ifdim\wd\myboxB<\wd\myboxA%
	\rlap{\hskip 0.5\mylenA\usebox\myboxB}{\usebox\myboxA}%
	\else
	\hskip -0.5\mylenA\rlap{\usebox\myboxA}{\hskip 0.5\mylenA\usebox\myboxB}%
	\fi}
\makeatother

\newcommand{\pRZ}{\mathbb{R}_{>0}}
\newcommand{\pRZz}{\mathbb{R}_{\geq0}}

\usepackage{etoolbox}
\makeatletter
\providecommand*{\cupdot}{%
	\mathbin{%
		\mathpalette\@cupdot{}%
	}%
}
\newcommand*{\@cupdot}[2]{%
	\ooalign{%
		$\m@th#1\cup$\cr
		\hidewidth$\m@th#1\cdot$\hidewidth
	}%
}
\makeatother
\makeatletter
\providerobustcmd*{\bigcupdot}{%
	\mathop{%
		\mathpalette\bigop@dot\bigcup
	}%
}
\newrobustcmd*{\bigop@dot}[2]{%
	\setbox0=\hbox{$\m@th#1#2$}%
	\vbox{%
		\lineskiplimit=\maxdimen
		\lineskip=-0.7\dimexpr\ht0+\dp0\relax
		\ialign{%
			\hfil##\hfil\cr
			$\m@th\cdot$\cr
			\box0\cr
		}%
	}%
}
\makeatother

\newcommand{\Hilbert}[1]{$(\mathbb{H},\langle\,\cdot\,|\,\cdot\,\rangle_{\mathbb{H}})$}

\def\DHLhksqrt#1#2{\setbox0=\hbox{$#1\sqrt{#2\,}$}\dimen0=\ht0
	\advance\dimen0-0.2\ht0
	\setbox2=\hbox{\vrule height\ht0 depth -\dimen0}%
	{\box0\lower0.4pt\box2}}

\newcommand{\orange}[1]{\textcolor{orange}{#1}}
\newcommand{\magenta}[1]{\textcolor{black}{#1}}

\newcommand{\contrastbar}{\operatorname{contr}^{\rmwU}}
\newcommand{\contrasthat}{\widehat{\operatorname{contr}}}

\newcommand{\kbar}{\tilde{k}}

\newcommand{\linhat}{\widehat{\lambda}^{\rmwhatU}}
\newcommand{\lin}{\lambda}
\newcommand{\linbar}{\widehat{\lambda}^{\rmwU}}
\newcommand{\lintilde}{\widehat{\lin}^{\rmwU}}
\newcommand{\lincheck}{\widecheck{\lin}}

\newcolumntype{C}[1]{>{\centering\arraybackslash}p{#1}}

\newcommand{\phib}{\phi^{\operatorname{b}}}
\newcommand{\Phib}{\Phi^{\operatorname{b}}}
\newcommand{\phiu}{\phi^{\operatorname{u}}}

\begin{document}


\begin{center}
    {\LARGE
        \textsc{Linear Functional Estimation\\ under Multiplicative Measurement Error\\ with Unknown Error Distribution}
    \par}

    \vspace{1cm}

 \begin{center}
{\Large\scshape
\begin{tabular}{c c c}
Jan Johannes\footnotemark[1] & \hspace{1cm}\&\hspace{1cm} & Maximilian Siebel\footnotemark[2]\\
[-0.3em]
{\small Heidelberg University} & & {\small  Humboldt University of Berlin}
\end{tabular}
}
\end{center}

    \vspace{0.5cm}

    {\small
        \textsc{\today}
    \par}
\end{center}

\footnotetext[1]{Institute for Mathematics, Heidelberg University; e-mail: \href{mailto:johannes@math.uni-heidelberg.de}{johannes@math.uni-heidelberg.de}}
\footnotetext[2]{Institute for Mathematics, Humboldt University of Berlin; e-mail: \href{mailto:maximilian.siebel@hu-berlin.de}{maximilian.siebel@hu-berlin.de}}

\begin{abstract}
    \noindent
    In this work we consider i.i.d. observations from a multiplicative deconvolution model, i.e. generated by $Y=X\varepsilon$, where $X$ and $\varepsilon$ are independent strictly positive random variables.
    The goal is to statistically infer the unknown signal $X$.
    Starting from the scenario where the probability density $f_\varepsilon$ of the error is known, we revisit local estimation strategies for several quantities of interest including the probability density function $f$ and the survival function $\SX$ of the signal $X$.
    We use these results to propose a thresholded plug-in estimation technique for linear functionals in situations where the distribution of the error $\varepsilon$ is unknown, while having some additional data at our disposal.
    This estimation procedure combines the reconstruction of the unknown Mellin transform of $f$ as well as its regularized inversion by a spectral cut-off.
    We establish oracle-type upper bounds and discuss rates of convergence in explicit scenarios including estimating $f$ and $\SX$ pointwise.
    In addition, we develop a fully data-driven choice of the cut-off parameter via model selection, which is able to obtain the oracle risk up to a constant without any previous
knowledge of the error distribution.
The results are illustrated in a simulation experiment, where we compare different classes of distributions in settings with known and unknown $f_\varepsilon$.

\end{abstract}

\tableofcontents

\section{Introduction}

Throughout this work, let $X$ and $\varepsilon$ be two independent, strictly positive random variables admitting probability density functions $f \coloneqq f_X$ and $f_\varepsilon$, respectively, with respect to the Lebesgue measure $\poslebesgue$ on the positive real line. We observe i.i.d. copies $Y_1,\dots,Y_N$ of the random variable generated by the multiplicative measurement model
\begin{equation*}
Y = X\varepsilon .
\end{equation*}
In this situation $Y$ also admits a probability density function $f_Y$, given by the multiplicative convolution
\begin{equation*}
\forall y\in\pRz: f_Y(y) = (f_X\ocon f_\varepsilon)(y) \coloneqq \int_{\pRz} f_X(x) f_\varepsilon(yx^{-1})x^{-1}\mathrm{d}\poslebesgue(x).
\end{equation*}
Our aim is to estimate the unknown distribution of the signal $X$ from the observations $Y_1,\dots,Y_N$, while the error density $f_\varepsilon$ is itself unknown,
which naturally leads to a linear inverse problem with an unknown operator, commonly referred to as multiplicative deconvolution.
Multiplicative measurement models historically arose from the \textit{multiplicative censoring problem}, which was first introduced and studied in \cite{Vardi_1989} and \cite{VardiZhang_1992}. We note that multiplicative censoring corresponds to the special case of multiplicative deconvolution in which the multiplicative error $\varepsilon$ is uniformly distributed on the interval $[0,1]$. Similarly, \cite{vanEsetal_2000} motivate and explain multiplicative censoring in the context of survival time analysis.
Building on these ideas from the late 1980s and 1990s, various reconstruction methods have been proposed for different special cases. While \cite{VardiZhang_1992} and \cite{AsgharianWolfson_2005} consider the estimation of the cumulative distribution function $F^X$ of $X$, \cite{AndersenHansen_2001} treat the problem as an inverse problem and study series expansion methods. Density estimation in the multiplicative censoring model is considered in \cite{Bruneletal_2016}, where a kernel estimator and a convolution power kernel estimator are proposed. Moreover, \cite{ComteDion_2016} propose a density estimator based on a projection method with respect to the Laguerre basis, assuming that $\varepsilon$ is uniformly distributed on an interval $[1-\alpha,1+\alpha]$ for $\alpha \in (0,1)$. A beta-distributed error $\varepsilon$ is studied in \cite{Belomestnyetal_2016}.

\medskip

\noindent\textbf{Related Work and Contribution}
The authors of \cite{BelomestnyGoldenshluger_2020} study pointwise density estimation, while \cite{Migueletal_2023} cast the problem as the estimation of the value of a known linear functional.
Both works assume the error distribution to be known. Global estimation techniques in multiplicative measurement models were considered in \cite{Migueletal_2021}.
The authors use the Mellin transform together with a spectral cut-off regularization of its inverse to define an estimator for the unknown density $f$.
Similarly, under multiplicative measurement errors, \cite{BrennerMiguel_2022} studies global multivariate density estimation, while global estimation of the survival function is considered in \cite{BrennerMiguelPhandoidaen_2022}.
Local and global estimation techniques with an unknown error distribution have been studied in the context of additive deconvolution problems, for instance in \cite{NeumannHossjer_1997} and \cite{Johannes_2009}.
For the multiplicative measurement model, the authors of \cite{Migueletal_2024} study global density and survival function estimation with an unknown error distribution.
The theory developed in the present paper is inspired by these works on multiplicative measurement models.
In the spirit of \cite{BrennerMiguel_2023} and \cite{Migueletal_2024}, we study local estimation techniques for linear functionals of the unknown density $f$ in the multiplicative measurement model,
while the error distribution is unknown.
We revisit the strategies and results proposed in \cite{Migueletal_2023} for the case of a known error distribution, and adapt them to the statistical regime in which the error distribution is unknown.
By comparing oracle-type risk bounds in the cases of known and unknown error densities, we characterize the influence of estimating the error density on the quality of the resulting estimator.
This comparison naturally involves the two sample sizes $N$ and $M$, corresponding to the observations from the signal-contaminated model and from the error distribution, respectively.
The proposed estimation strategy depends on a tuning parameter $k$, which has to be chosen by the practitioner.
In the case of a known error density, \cite{Migueletal_2023} propose a data-driven choice of the spectral cut-off parameter $k$ by means of model selection techniques, exploiting the theory of \cite{Barronetal_1999}.
We further refer to \cite{Massart_2007} for an extensive overview.\\

Our aim is to develop a fully data-driven procedure for selecting the tuning parameter $k$ in a local estimation framework
when the error distribution is unknown, following the well-established adaptive selection methodology of \cite{GoldenshlugerLepski_2011}.
The main challenge in our setting lies in constructing a suitable data-driven penalty: existing adaptive procedures, such as those in \cite{Migueletal_2021},
involve the error distribution in a form that is not directly accessible from the given data.
The extension to the setting when the error distribution is unknown therefore requires a new penalty.
To this end, we propose an alternative penalty and first establish the resulting selection procedure when the error distribution is known.
We then extend the construction and its analysis to the unknown-error setting.
Related procedures from known to unknown error distributions have been studied in additive deconvolution,
for instance in \cite{ComteLacour_2011} and \cite{JohannesSchwarz_2013}, and for global estimation techniques in multiplicative measurement models in \cite{Migueletal_2024}.

\newpage

\noindent\textbf{Outline.}
The paper is structured as follows.
We begin with a detailed introduction to the statistical setting in \cref{sec:background} and revisit the most essential definitions and properties of the Mellin transform.
In \cref{sec:lifumu:known} we revisit the techniques derived in \cite{BrennerMiguel_2023} when the error density is known.
In particular, we introduce the spectral cut-off estimator for linear functionals of the unknown density $f$,
and provide an oracle-type inequality for the mean squared error and discuss minimax-optimal rates of convergence in concrete scenarios; see \cref{example:lifumu:known:oracle:rates}.
In \cref{sec:lifumu:datadriven:known}, we revisit data-driven estimation procedures for the tuning parameter $k$ based on model selection techniques.
In \cref{sec:lifumu:unknown:oracle} and \cref{sec:lifumu:unknown:datadriven}, we adapt these findings to the statistical regime in which the error distribution is unknown and provide oracle-type upper bounds as well as a data-driven model selection procedure for the tuning parameter $k$.
Corresponding rates of convergence are discussed in \cref{example:lifumu:unknown:oracle:rates} and \cref{example:lifumu:unknown:datadriven:rates}, respectively.
Lastly, we provide illustrations of the proposed estimation techniques in \cref{sec:simulation-experiments}.
Technical proofs are postponed to the appendix.

\section{Model Assumptions and Theoretical Background}\label{sec:background}

In the following we introduce the multiplicative measurement
model in detail, and recapitulate the definition of the Mellin
transform, its empirical counterpart as well as their properties.
This introduction follows the lines of \cite{Migueletal_2024}.

\subsection{The Multiplicative Measurement Model}

In the following, $(\IR,\scrB_{\IR},\lebesgue)$ denotes the real line $\IR$ equipped with its Borel $\sigma$-algebra $\scrB_{\IR}$ and the usual Lebesgue measure $\lebesgue$.
For any $a\in\IR$ we set $\IR_{>a}\coloneqq (a,\infty)$.
Thus, $(\pRz,\scrB_{\pRz},\poslebesgue)$ denotes the corresponding measure space on the strictly positive real line $\pRz$.
We assume the signal and error $X,\varepsilon\colon(\Omega,\calA,\IP)\to(\pRz,\scrB_{\pRz})$ to be two strictly positive and independent random variables
both defined on a joint probability space $(\Omega,\calA,\IP)$. 
Their marginal probability distributions are denoted by $\IP^X$ and $\IP^\varepsilon$, respectively, 
and are assumed to admit probability density functions both with respect to the Lebesgue measure $\poslebesgue$, which are denoted by 
 $f \coloneqq f_X$ and $f_\varepsilon$.
The \textit{multiplicative measurement model} characterizes observations arising from 
\begin{equation}
	Y = X\varepsilon.
\end{equation}
In this situation $Y\sim\IP^Y$ also admits a probability density function which is given by the \textit{multiplicative convolution} of $f$ and $f_\varepsilon$, i.e.
\begin{equation*}
	\forall y\in\pRz: f_Y(y) = (f\ocon f_\varepsilon)(y) \coloneqq \int_{\pRz} f(x) f_\varepsilon(yx^{-1})x^{-1}\mathrm{d}\poslebesgue(x).
\end{equation*}
For an extensive overview of multiplicative convolution and its properties as an operator between function spaces we refer to \cite{BrennerMiguel_2023}.
In this work, we consider two statistical regimes:
	\begin{itemize}
		\item \newAssuRef{SRI} \textbf{Known Error Distribution.} Assuming the error probability density of the multiplicative error $\varepsilon$ to be known to us,
			we aim to statistically infer the unknown distribution $\IP^X$ of the signal $X$.
			To that end, given a fixed sample size $N\in\IN$, we observe $Y_1,\dots,Y_N$ generated as i.i.d. copies of the response variable $Y\sim\IP^Y$.
			If $f$ denotes the ground truth, i.e. the data generating probability density function of the signal $X$, we denote by $\IP_{f}^N$ the
		 	law of the observed data $(Y_i)_{i\in\nset{N}}$ on $(\pRz^N,\scrB_{\pRz}^{\otimes N})$. 
			The corresponding expectation under  $\IP_{f}^N$ is denoted by  $\IE_{f}^N$.

		\item \newAssuRef{SRII} \textbf{Unknown Error Distribution.} In the second statistical regime, we follow the same goal of statistically inferring the unknown distribution $\IP^X$ of the signal $X$, 
			but now without knowing the error distribution $\IP^\varepsilon$ of the multiplicative noise $\varepsilon$.
			Instead, we observe not only the sample $(Y_i)_{i\in\nset{N}}$, 
			but also \textit{additional data} $\varepsilon_1,\dots,\varepsilon_M$ with a second sample size $M\in\IN$, which is generated as i.i.d. copies of $\varepsilon\sim\IP^\varepsilon$
			and independently of the first sample $(Y_i)_{i\in\nset{N}}$.
			The law of the additional data on $(\pRz^M,\scrB_{\pRz}^{\otimes M})$ is denoted by $\IP_\varepsilon^M$ with corresponding expectation $\IE_{\varepsilon}^M$.
			As before, 
			if $f$ denotes the ground truth associated to the signal $X$, 
			we denote by 
			\begin{equation*}
				\IP_{f,\varepsilon}^{N,M}\coloneqq \IP_{f}^N \otimes \IP_{\varepsilon}^M
			\end{equation*}
			the joint law of the full data set 
			$$D_{N,M}\coloneqq\br{\br{Y_i}_{i\in\nset{N}},\br{\varepsilon_j}_{j\in\nset{M}}},$$
			which defines a probability measure on 
			$(\pRz^{N}\times\pRz^{M},\scrB_{\pRz}^{\otimes N}\otimes\scrB_{\pRz}^{\otimes M})$.
			The corresponding expectation is denoted by $\IE_{f,\varepsilon}^{N,M}$.
			Due to independence of the samples $(Y_i)_{i\in\nset{N}}$ and $(\varepsilon_j)_{j\in\nset{M}}$, we note that for any $\IP_{f,\varepsilon}^{N,M}$-integrable $h\colon\pRz^{N}\times\pRz^{M} \to \IC$,
			the function
				\begin{equation*}
					g\colon\pRz^{M}\ni e \mapsto g(e)=\IE_{f}^N\brE{h(\cdot,e)} \in \IC
				\end{equation*}
				is a version of the conditional expectation of $h$ given the additional data $\br{\varepsilon_j}_{j\in\nset{M}}$.
				In particular, we set
				\begin{equation*}
					\IE_{f|\varepsilon}^{N,M}\brE{h} \coloneqq g 
				\end{equation*}
	\end{itemize}

\subsection{The (empirical) Mellin Transform}
			Fix $p\in[1,\infty]$. We denote by $\IL^p \coloneqq \IL^p(\IR,\scrB_\IR,\lebesgue)$ the usual space of $p$-integrable functions $f\colon\IR\to \IC$ with respect to the Lebesgue measure $\lebesgue$.
			For any $a\in\IR$ we define the weight function 
			\begin{equation*}
				\rmx^{a}\colon\pRz\to\pRz,\; x\mapsto x^{a}.
			\end{equation*}
			On $(\pRz,\scrB_{\pRz})$ we consider the weighted Lebesgue measure $\rmx^a\poslebesgue$, and 
			use throughout the shorthand notation
			\begin{equation*}
				\IL_{\rmx^{a}}^p \coloneqq \IL^p(\pRz,\scrB_{\pRz},\rmx^a\poslebesgue)
			\end{equation*}
			to denote the usual space of $p$-integrable functions $f\colon\pRz\to \IC$ with respect to the measure $\rmx^{a}\poslebesgue$.
			Fix some $c\in\IR$.
		The Mellin
transform $\M_c$ is the linear and bounded operator between the
Hilbert spaces $\IL^2_{\rmx^{2c-1}}$ and $\IL^2$, which for each
$h\in\IL^2_{\rmx^{2c-1}}\cap\IL^1_{\rmx^{c-1}}$ and $t\in\IR$ satisfies
\begin{equation}\label{def:mellin}
	\M_c[h](t)\coloneqq\int_{\pRz}x^{c-1+\iota 2\pi t}h(x)\rmd\poslebesgue(x),
\end{equation}
where $\iota\coloneqq\sqrt{-1}\in\IC$ denotes the imaginary unit. Similar to the
Fourier transform (see
\cite{Meister_2009} for definitions and properties) the Mellin
transform $\M_c$  is unitary, i.e.
$\skalar{h,g}_{\IL^2_{\rmx^{2c-1}}} =
\skalar{\M_c[h],\M_c[g]}_{\IL^2}$ for any
$h,g\in\IL^2_{\rmx^{2c-1}}$. In particular, it satisfies a \textit{Plancherel type identity}
\begin{equation}\label{eq:Plancherel}
	\norm{h}_{\IL^2_{\rmx^{2c-1}}}^2 = \norm{\M_c[h]}_{\IL^2}^2,\quad\forall h\in\IL^2_{\rmx^{2c-1}}.
\end{equation}
Its adjoint (and inverse) $\M_c^{\star}:\IL^2\to\IL^2_{\rmx^{2c-1}}$ fulfils for each
$g\in\IL^2\cap\IL^1$ and $x\in\pRz$
\begin{equation}\label{def:mellin:inv}
	\M_c^{\star}[g](x)\coloneqq\int_{\Rz}x^{-c-\iota 2\pi t}g(t)\rmd\lambda(t).
\end{equation}
In analogy to the additive
convolution theorem of the Fourier transform, there is a
multiplicative convolution theorem for the Mellin
transform.  Namely, for any $h_1\in\IL^2_{\rmx^{2c-1}}$ and $h_2\in\IL^1_{\rmx^{c-1}}\cap
\IL^2_{\rmx^{2c-1}}$ the multiplicative convolution theorem states 
\begin{equation}\label{eq:multconvothm}
	\M_c[h_1\ocon h_2] = \M_c[h_1] \M_c[h_2].
\end{equation}
For a detailed discussion of the Mellin
transform and its properties we refer again to \cite{BrennerMiguel_2023}.
In order to simplify the presentation in the
following sections, we introduce some frequently used shorthand
notations.
\begin{notation}[The (empirical) Mellin Transform \& Constants]\label{not: MellinT}
	\mbox{}
	\begin{itemize}
		\item The Mellin transformations of the densities $f_Y$, $f$ and $f_\varepsilon$, respectively, are briefly denoted by 
	\begin{equation*}
		\MelY(t)\coloneqq \M_c[f_Y](t) = \MelX(t)\MelU(t),\quad	\MelX(t)\coloneqq \M_c[f](t)\quad\text{and}\quad \MelU(t)\coloneqq\M_c[f_\varepsilon](t)
	\end{equation*}
	for all $t\in\IR$.
	The reciprocal of a function $\rmw:\Rz\to\IC$ is well-defined on the set $\{\rmw\neq 0\}\coloneqq\{t\in\IR:\,\rmw(t)\neq 0\}$ and for each $t\in\IR$ we write for short
	\begin{equation*}
		\rmwrez(t)\coloneqq\frac{1}{\rmw(t)}\1_{\{\rmw\neq 0\}}(t)
	\end{equation*}
	with $\1_{\{\rmw\neq 0\}}$ denoting the indicator function. Precisely, for any set $A\subseteq\IR$ we write $\1_A$ for the indicator function, mapping $x\in\IR$ to $\1_A(x)\coloneqq1$, if $x\in A$, and  $\1_A(x)\coloneqq0$, otherwise.
	In particular, $\rmwrez(t)=0$ for $t\in\{\rmw = 0\}$. Since the Mellin transformation $\MelY$ of $f_Y$ is unknown we follow \cite{BrennerMiguel_2023} and introduce an empirical counterpart based on the observations $(Y_i)_{i\in\nset{N}}$.
	The \textit{empirical Mellin transformation} $\MelYhat$  is defined by 
	\begin{equation*}
		\MelYhat(t)\coloneqq\frac{1}{N}\sum_{i\in\nset N}Y_{i}^{c-1+\iota2\pi t}
	\end{equation*}
	for any $t\in\IR$. Observe that $\MelYhat$ is a pointwise unbiased estimator of $\MelY$, i.e. for all $t\in\IR$, we have $	\IE_f^N[\MelYhat(t)] = \MelY(t)$.
	In \assuRef{SRI}, when the error distribution is known, a corresponding estimator for the unknown Mellin transformation $\MelX$ is given by 
	\begin{equation*}
		\forall t\in\IR:\;\MelXtilde(t) \coloneqq\MelYhat(t)\MelUrez(t).
	\end{equation*}
	Further, for any $k\in\pRz$, we write $\MelXtilde^k(t) \coloneqq \MelXtilde(t)\1_{[-k,k]}(t)$. Similarly, in \assuRef{SRII}, we define the empirical Mellin transformation $\MelUhat$ by 
	\begin{equation*}
		\MelUhat(t)\coloneqq\frac{1}{M}\sum_{i\in\nset M}\varepsilon_{i}^{c-1+\iota2\pi t},
	\end{equation*}
	pointwise unbiased estimator of $\MelU$ in the sense that $\IE_\varepsilon^M[\MelUhat(t)] = \MelU(t)$ for all $t\in\IR$. 
	Given an additional (random) threshold event $\frakMl$ (to be specified later), we define a corresponding estimator $\MelXhat$ of $\MelX$ by 
	\begin{equation*}
		\forall t\in\IR:\; \MelXhat(t) \coloneqq \MelYhat(t)\MelUhatrez(t)\1_{\frakMl}(t)
	\end{equation*}
	and accordingly $\MelXhat^k(t)\coloneqq\MelXhat(t)\1_{[-k,k]}(t)$ for all $k\in\pRZ$.
	\item Throughout this work, $\const\in\pRz$ denotes a positive universal constant whose numerical value may vary from line to line.
	If such a constant depends on a family of objects $\rmO$, we indicate this dependence by writing $\const(\rmO)\in\pRz$.
	 Positive constants arising from (model) assumptions are generally denoted by $\Const\in\pRz$, with additional indices or labels as appropriate.
	\end{itemize}
\end{notation}

We have the following preliminary inequality.
\begin{lemma}[Preliminary Inequality]\label{lem:MellinTrick}
	\mbox{}\\
	Let $Y = X\varepsilon$ with $f_Y\in\IL^1_{\rmx^{2(c-1)}}$ and $f_\varepsilon\in\IL^1_{\rmx^{2(c-1)}}\cap\IL_{\rmx^{2c-1}}^\infty$. 
	 Let $h\colon\IR\to\IC$, such that $h\indk{k}\in\IL^2$ for all $k\in\pRz$. 
	 Defining
	\begin{equation*}
		\Constfeps \coloneqq  1\lor \tfrac{\norm{f^\varepsilon}_{\IL_{\rmx^{2c-1}}^\infty}}{\norm{f^\varepsilon}_{\IL_{\rmx^{2(c-1)}}^1}}\in\IR_{\geq 1}\quad\and\quad\ConstfY\coloneqq 1\lor \norm{f_Y}_{\IL_{\rmx^{2(c-1)}}^1}\in\IR_{\geq 1}, 
	\end{equation*}
	we have for all $k\in\pRz$
	\begin{equation*}
		\IE\brE{\abs{\int_{[-k,k]}Y^{c-1+2\pi\iota t}h(t)\rmd\lebesgue(t)}^2} \leq \Constfeps\ConstfY\times \norm{h\indk{k}}_{\IL^2}^2.
	\end{equation*}
\end{lemma}

\begin{proof}[Proof of \cref{lem:MellinTrick}]
	The claim is proven in \cite[Lemma 1.2.5]{Neubert_2025}.
\end{proof}

\subsection{Linear Functionals and Smoothness Classes}

We follow the goal formulated in \cite{Migueletal_2023} and aim to reconstruct linear functionals of the form
\begin{equation*}
	\lin(f) \coloneqq \int_{\IR}\rmL(-t)\MelX(t)\rmd\lebesgue(t) = \int_{\IR}\rmL(-t)\MelY(t)\MelUrez(t)\rmd\lebesgue(t)
\end{equation*}
in both statistical regimes \assuRef{SRI} and \assuRef{SRII}.
Here, $\Lambda\colon\IR\to\IC$ is a suitable function, subsequently also called \textit{representant}. 
	The following example summarizes different choices of the representant $\rmL$ that
		induce concrete quantities $\lin(f)$ with a statistical interpretation.
		

		\begin{example}[The Choice of $\rmL$]\label{example:Illustration1}
			\mbox{}
			Let $x_o\in\pRz$ be fixed.
			\begin{itemize}
				\item\label{lifumu:item:ex1} \newAssuRef{Probability Density Function} 
					Let $c\in\IR$ and assume that $f\in\IL_{\rmx^{c-1}}^1\cap\IL_{\rmx^{2c-1}}^2$ as well as $\MelX\in\IL^1$.
					We then have
					\begin{equation*}
						f(x_o) = \M_c^{-1}\brE{\MelX}(x_o) = \int_{\IR}x_o^{-c-2\pi\iota t}\MelX(t)\rmd\lebesgue(t).
					\end{equation*}
					Defining
					\begin{equation*}
						\forall t\in\IR:\;\rmL_f(t)\coloneqq x_o^{-c+2\pi\iota t},
					\end{equation*}
					we obtain $\lin(f) = f(x_o)$.
	
				\item\label{lifumu:item:ex2} \newAssuRef{Derivative of the Probability Density Function} 
						Let $c\in\IR$ and $\beta\in\IN$. Let $f\in C_{\rmc}^\infty(\pRz)$,
						such that $f\in\IL_{\rmx^{c-1}}^1\cap\IL_{\rmx^{2c-1}}^2$. We denote by $\rmD^\beta f$ the $\beta^{th}$ derivative of $f$.
						Exploiting \cite[Proposition 2.3.11]{BrennerMiguel_2023}, we have 
						\begin{equation*}
							\forall t\in\IR:\;\M_{c+\beta}[\rmD^\beta f](t) = (-1)^\beta\prod_{j\in\nset{\beta}}\br{c+\beta-j+2\pi\iota t}\MelX(t).
						\end{equation*}
						Thus, if $(-1)^\beta\prod_{j\in\nset{\beta}}\br{c+\beta-j+2\pi\iota \cdot}\MelX\in\IL^1\cap\IL^2$, we write 
						\begin{align*}
							\rmD^\beta f(x_o) & = \M_{c+\beta}^{-1}\brE{(-1)^\beta\prod_{j\in\nset{\beta}}\br{c+\beta-j+2\pi\iota \cdot}\MelX}(x_o)\\
								& = \int_{\IR}(-1)^\beta\prod_{j\in\nset{\beta}}\br{c+\beta-j+2\pi\iota t}x_o^{-c-\beta-2\pi\iota t}\MelX(t)\rmd\lebesgue(t),
						\end{align*}
						and hence define
						\begin{equation*}
							\forall t\in\IR:\; \rmL_{\rmD^\beta f}(t)\coloneqq (-1)^\beta\prod_{j\in\nset{\beta}}\br{c+\beta-j-2\pi\iota t}x_o^{-c-\beta+2\pi\iota t},
						\end{equation*}
						such that $\lin(f) = \rmD^\beta f(x_o)$.

			\item\label{lifumu:item:ex3} \newAssuRef{Survival Function and CDF} 
						We denote by $F^X(x)\coloneqq \IP(X\leq x)$, $x\in\pRz$, the usual cumulative distribution function (CDF) of $X$. 
						The survival function $\SX$ is defined by $\SX(x)\coloneqq 1 - F^X(x)$. Let $c\in\IR_{>1}$. Following \cite{BrennerMiguelPhandoidaen_2022,Migueletal_2024}, we have $f\in\IL_{\rmx^{c-1}}^1\cap\IL_{\rmx^{2c-1}}^2$ if and only if
						$\SX\in\IL_{\rmx^{c-2}}^1\cap\IL_{\rmx^{2c-3}}^2$ and in particular
						\begin{equation*}
							\forall t\in\pRz:\; \M_{c-1}\brE{\SX}(t) = \br{c-1+2\pi\iota t}^{-1}\MelX(t),
						\end{equation*}
						as well as 
						\begin{equation*}
							 \SX(x_o) = \int_{\IR}\frac{x_o^{-c+1-2\pi\iota t}}{c-1+2\pi\iota t}\MelX(t)\rmd\lebesgue(t),
						\end{equation*}
						whenever $\MelX\in\IL^1\cap\IL^2$.
						We hence define
						\begin{equation*}
							\forall t\in\IR:\; \rmL_{\SX}(t)\coloneqq  x_o^{-c+1+2\pi\iota t}\br{c-1-2\pi\iota t}^{-1}
						\end{equation*}
						to obtain $\lin(f) = \SX(x_o)$.
						With similar arguments, one also defines for $c\in\IR_{<1}$
						\begin{equation*}
							\forall t\in\IR:\; \rmL_{F^X}(t)\coloneqq \br{1-c+2\pi\iota t}^{-1}x_o^{1-c+2\pi\iota t}
						\end{equation*} 
						to obtain $\lin(f) = F^X(x_o)$.
			\end{itemize} 
		\end{example}

\subsubsection{Smoothness Classes}\label{sec:Muduner:Smoothness}

Below, we propose estimation strategies for the general quantity $\lin(f)$. 
To quantify the resulting estimation accuracy we consider the Mean Squared Error (MSE) and 
derive corresponding upper bounds.
To obtain concrete rates of convergence for the MSE
we need to specify the regularity of the unknown probability density function 
$f$ of the signal $X$, as well as the decay of the Mellin transform $\MelU$ of the error probability density function $f^\varepsilon$.

\begin{definition}[Mellin Sobolev Spaces]\label{def:ordinaryandexponentially}
	\mbox{}\\
	Let $c\in\IR$ be fixed.
	\begin{itemize}
		\item For all $\rmb,\rms\colon\IR\to\pRzz$ continuous and symmetric (that is $\rmb(t)=\rmb(-t)$ and $\rms(t)=\rms(-t)$ for all $t\in\IR$), we define 
			\begin{equation*}
				H_\rms\coloneqq\brK{h\in\IL^2_{\rmx^{2c-1}}:\;\norm{\rms\M_c[h]}_{\IL^2}^2<\infty}
			\end{equation*}
			and the corresponding ball of radius $\rmR\in\pRz$ by 
			\begin{equation*}
				H_\rms(\rmR)\coloneqq\brK{h\in H_\rms:\;\norm{\rms\M_c[h]}_{\IL^2}^2\leq\rmR^2}.
			\end{equation*}
			Further, we set 
			\begin{equation*}
				H_\rmb^\infty(\rmR)\coloneqq\brK{h\in\IL^2_{\rmx^{2c-1}}:\;\rmb^\dagger\leq\rmR\abs{\M_c[h]}\text{ $\lebesgue$-a.s.}}.
			\end{equation*}
		\item 
	Now, let $b,s\in\pRzz$, and define for all $t\in\IR$
	\begin{equation*}
		\rmb_{\textsc{os}}(t)\coloneqq \br{1+\abs{t}^2}^{\frac{b}{2}}\quad\text{and}\quad \rms_{\textsc{os}}(t)\coloneqq \br{1+\abs{t}^2}^{\frac{s}{2}}
	\end{equation*} 
	as well as
	\begin{equation*}
		\rmb_{\textsc{es}}(t)\coloneqq \exp\br{\abs{t}^{2b}}\quad\text{and}\quad \rms_{\textsc{es}}(t)\coloneqq \exp\br{\abs{t}^{2s}}.
	\end{equation*}
	We distinguish between the following regularities of the probability density function $f^\varepsilon$.
	\begin{enumerate}
		\item[i)] We say that $f_\varepsilon$ is \textit{ordinary smooth} (\textsc{os}), if $f_\varepsilon\in H_{\rmb_{\textsc{os}}}^\infty(\rmr_\infty)$ for some $\rmr_\infty\in\IR_{\geq 1}$.
		\item[ii)] We say that $f_\varepsilon$ is \textit{super smooth} (or \textit{exponentially smooth} (\textsc{es})), if $f_\varepsilon\in H_{\rmb_{\textsc{es}}}^\infty(\rmr_\infty)$ for some $\rmr_\infty\in\IR_{\geq 1}$.
	\end{enumerate} 
	Similarly, we distinguish between the following smoothness classes for the probability density function $f$.
	\begin{enumerate}
		\item[iii)] We say that $f$ is \textit{ordinary smooth} (\textsc{os}), if $f\in H_{\rms_{\textsc{os}}}(\rmr_2)$ for some $\rmr_2\in\pRz$.
		\item[iv)] We say that $f$ is \textit{super smooth} (or \textit{exponentially smooth} (\textsc{es})), if $f\in H_{\rms_{\textsc{es}}}(\rmr_2)$ for some $\rmr_2\in\pRz$.
	\end{enumerate}
	\end{itemize}
\end{definition}

\section{Linear Functional Estimation - Known Error Distribution}\label{sec:lifumu:known}

			We follow the approach of \cite{Migueletal_2023} (keeping \cref{not: MellinT} in mind) and propose a 
			\textit{spectral cut-off estimator} of $\lin(f)$ by
			\begin{equation}\label{def:linearfunctionalestimator:known}
				\forall k\in\pRz:\; \lintilde_k \coloneqq \int_{\IR}\rmL(-t)\MelXtilde^k(t)\rmd\lebesgue(t)
					 = \int_{[-k,k]}\rmL(-t)\MelYhat(t)\MelUrez(t)\rmd\lebesgue(t),
			\end{equation}
			whenever the integral is well-defined,
			where $\rmwU\coloneqq\MelUrez$ indicates the dependence of the Mellin transform $\MelU$ in what follows.
			For the further risk analysis, let us summarize the following assumptions.
			\begin{assumption}[Linear Functional Estimation]\label{ass:lifumu:A1}
				\mbox{}\\
				Suppose the statistical regime \assuRef{SRI}.
				\begin{itemize}
					\item Fix $c\in\IR$. Let $f^\varepsilon\in\IL_{\rmx^{c-1}}^1\cap\IL_{\rmx^{2(c-1)}}^1\cap \IL_{\rmx^{2c-1}}^\infty$, $f\in\IL_{\rmx^{2(c-1)}}^1\cap\IL_{\rmx^{2c-1}}^2$ and $f_Y\in\IL_{\rmx^{2(c-1)}}^1$.\\
						We set 
						\begin{equation*}
							\Constfeps \coloneqq  1\lor \frac{\norm{f^\varepsilon}_{\IL_{\rmx^{2c-1}}^\infty}}{\norm{f^\varepsilon}_{\IL_{\rmx^{2(c-1)}}^1}}\in\IR_{\geq 1}\quad\and\quad\ConstfY\coloneqq 1\lor \norm{f_Y}_{\IL_{\rmx^{2(c-1)}}^1}\in\IR_{\geq 1}. 
						\end{equation*}
					\item Let  $\rmL\colon\IR\to\IC$ be a measurable function, such that $\overline{\rmL(t)} = \rmL(-t)$ for all $t\in\IR$.
					\item Let $\MelU(t)\neq 0$ for $\lebesgue$-almost all $t\in\IR$, and assume that $\rmL\MelX\in\IL^1$
						as well as $\rmL\MelUrez\indk{k}\in\IL^\infty$ for all $k\in\pRz$.
					
				\end{itemize}
			\end{assumption}

			We note that the conditions on $f^\varepsilon$ imposed in \cref{ass:lifumu:A1} already imply
$
\norm{\MelU}_{\IL^\infty}<\infty.
$
Moreover, the condition
$
f\in\IL_{\rmx^{2(c-1)}}^1
$
can be interpreted in terms of moment requirements.
Square integrability is a standard assumption in additive deconvolution and is recovered here as the particular case $c=1/2$. The choice $c=0$, on the other hand, corresponds to the naive approach of log-transforming the data and subsequently applying an additive deconvolution method.
Allowing $c$ to vary over $\IR$ makes explicit how the assumptions depend on this parameter. For instance, when $c=0$, the assumptions above require $X^{-1}$ and $\varepsilon^{-1}$ to have finite (second) moments. Consequently, whenever our assumptions are satisfied for $c=0$, the following results also apply to the naive log-transformation approach.

\subsection{Oracle-type Inequality}\label{sec:lifumu:oracle:known}

We start analyzing the MSE of $\lintilde_k$ in the statistical regime \assuRef{SRI}, i.e. with a known error distribution (KED).
\begin{lemma}[MSE - Known Error Distribution]\label{lem:Lifumu:known:MSE}
	\mbox{}\\
	Grant \cref{ass:lifumu:A1}. We then have for all $N\in\IN$ and $k\in\pRz$
	\begin{equation*}
		\IE_f^N\brE{\abs{\linbar_k-\lin(f)}^2}  \leq \norm{\rmL\indkco{k}\MelX}_{\IL^1}^2 + \frac{1}{N}\Constfeps\ConstfY\norm{\rmL\indk{k}\MelUrez}_{\IL^2}^2.
	\end{equation*}
\end{lemma}

\begin{proof}[Proof of \cref{lem:Lifumu:known:MSE}]
	Fix $N\in\IN$ and $k\in\pRz$. We recap the proof of \cite[Proposition 2.3]{Migueletal_2023} and start with a squared bias-variance decomposition
	from which we obtain
	\begin{align*}
		\IE_f^N\brE{\abs{\linbar_k-\lin(f)}^2} & 
		 = \IE_f^N\brE{\abs{\linbar_k-\IE_f^N\brE{\linbar_k}}^2}+\abs{\IE_f^N\brE{\linbar_k}-\lin(f)}^2\\
		&\eqqcolon \Var_f^N\br{\linbar_k}+ \sqBias{k}.
	\end{align*}
	Keeping in mind that $\MelYhat$ is a pointwise unbiased estimator of $\MelY$, an application of Fubini's theorem yields
	\begin{equation*}
		\IE_f^N\brE{\linbar_k} = \int_{[-k,k]}\rmL(-t)\MelY(t)\MelUrez(t)\rmd\lebesgue(t)\eqqcolon \lin_k^{\rmwU}(f).
	\end{equation*}
	As a consequence, for the squared bias expression we conclude
	\begin{align*}
		 \sqBias{k} &= \abs{ \lin_k^{\rmwU}(f)-\lin(f)}^2 =  
		  \abs{\int_\IR \rmL(-t)\MelX(t)\indkco{k}(t)\rmd\lebesgue(t)}^2
		 \leq \norm{\rmL\MelX\indkco{k}}_{\IL^1}^2,
	\end{align*}
	where the right-hand side is finite by \cref{ass:lifumu:A1}.
		For the variance term, observe first 
		\begin{equation*}
			\linbar_k = \frac{1}{N}\sum_{i\in\nset{N}} \int_{[-k,k]}\rmL(-t)Y_i^{c-1+2\pi\iota t}\MelUrez(t)\rmd\lebesgue(t).
		\end{equation*}
		If we denote by $\Var$ the variance associated to $\IP$, we further conclude
		\begin{align*}
			\Var_f^N\br{\linbar_k} & = \frac{1}{N}\Var\br{\int_{[-k,k]}\rmL(-t)Y_1^{c-1+2\pi\iota t}\MelUrez(t)\rmd\lebesgue(t)}\\
				& \leq \frac{1}{N}\IE\brE{\abs{\int_{[-k,k]}\rmL(-t)Y_1^{c-1+2\pi\iota t}\MelUrez(t)\rmd\lebesgue(t)}^2}.
		\end{align*}
		By assumption $\rmL\MelUrez\indk{k}\in\IL^2$ and thus we can apply \cref{lem:MellinTrick} to obtain
		\begin{equation*}
			\frac{1}{N}\IE\brE{\abs{\int_{[-k,k]}\rmL(-t)Y_1^{c-1+2\pi\iota t}\MelUrez(t)\rmd\lebesgue(t)}^2} \leq \frac{1}{N}\Constfeps\ConstfY\times\norm{\rmL\MelUrez\indk{k}}_{\IL^2}^2.
		\end{equation*}
		Combining the bounds for the squared bias and the variance expressions shows the claim.
\end{proof}

\begin{corollary}[Local Estimation - Oracle-type Inequality - KED]\label{co:lifumu:known:oracle}
	\mbox{}\\
	Let $N\in\IN$. Under \cref{ass:lifumu:A1} there exists an optimal tuning parameter $\kdiamond = \kdiamond(N)\in\pRz$, such that 
		\begin{align*}
			\IE_f^N\brE{\abs{\linbar_{\kdiamond}-\lin(f)}^2}  
			&\leq \Constfeps\ConstfY\times\inf_{k\in\pRz}\brK{ \norm{\rmL\MelX\indkco{k}}_{\IL^1}^2+\frac{1}{N}\norm{\rmL\indk{k}\MelUrez}_{\IL^2}^2}.
		\end{align*}
		In particular, we have
		\begin{equation}\label{eq:lifumu:rateofconvergence}
			\lim_{N\to\infty}\IE_f^N\brE{\abs{\linbar_{\kdiamond}-\lin(f)}^2} =0.
		\end{equation}
\end{corollary}

\begin{proof}[Proof of \cref{co:lifumu:known:oracle}]
	Fix $N\in\IN$. By \cref{lem:Lifumu:known:MSE} we have for all $k\in\pRz$
	\begin{equation*}
		\IE_f^N\brE{\abs{\linbar_k-\lin(f)}^2}  \leq \Constfeps\ConstfY\times\br{\norm{\rmL\indkco{k}\MelX}_{\IL^1}^2 + \frac{1}{N}\norm{\rmL\indk{k}\MelUrez}_{\IL^2}^2}.
	\end{equation*}
	We now show that the right-hand side attains its global minimum at some $\kdiamond\in\pRz$. 
	Let us define $A_1,A_2\colon\pRz\to\pRzz$ by 
	\begin{equation*}
		A_1(k)\coloneqq \norm{\rmL\indkco{k}\MelX}_{\IL^1}^2\quad\and\quad A_2(k)\coloneqq \norm{\rmL\indk{k}\MelUrez}_{\IL^2}^2.
	\end{equation*}
	Evidently, $A_1$ defines a positive and non-increasing function, while $A_2$
	defines a non-decreasing function. By hypothesis $\rmL\MelX\in\IL^1$, and thus $\lim_{k\to\infty} A_1(k) = 0$
	due to the continuity of the Lebesgue-integral.
	Thus, to show existence of $\kdiamond = \kdiamond(N)\in\pRz$,
	such that 
	\begin{equation*}
		R_N(\kdiamond)\coloneqq A_1(\kdiamond) + \frac{1}{N}A_2(\kdiamond) = \inf_{k\in\pRz}\brK{A_1(k)+\frac{1}{N}A_2(k)},
	\end{equation*}
	it suffices to verify that $A_1$ and $A_2$ are continuous.
	Let $\delta\in\pRz$ be arbitrary. The continuity of the Lebesgue-integral implies
	\begin{equation*}
		\forall k\in\pRz:\;\abs{A_1(k+\delta)^\frac{1}{2}- A_1(k)^\frac{1}{2}} \leq \int_{\IR}\abs{\rmL(t)}\abs{\MelX(t)}\1_{(k,k+\delta]}(\abs{t})\rmd\lebesgue(t) \overset{\delta\to 0}{\longrightarrow}0.
	\end{equation*}
	Similarly,
	\begin{align*}
		\forall k\in\pRz:\;\abs{A_2(k+\delta) - A_2(k)} \leq \int_{\IR}\abs{\rmL(t)}^2\abs{\MelUrez(t)}^2\1_{(k,k+\delta]}(\abs{t})\rmd\lebesgue(t)  \overset{\delta\to 0}{\longrightarrow}0.
	\end{align*}
	It remains to show \cref{eq:lifumu:rateofconvergence}.
	Indeed, by definition of $\kdiamond$, we have 
	\begin{equation*}
		\forall k\in\pRz:\;R_N(\kdiamond) \leq R_N(k) = A_1(k) + \frac{1}{N}A_2(k).
	\end{equation*}
	As $\lim_{k\to\infty}A_1(k)=0$, we can find for any $\delta\in\pRz$ a cut-off parameter $\ko=\ko(\delta)\in\pRz$, such that $A_1(\ko)\leq \delta$.
	Thus,
	\begin{equation*}
		\limsup_{N\to\infty} R_N(\kdiamond) \leq \delta.
	\end{equation*}
	Since $\delta\in\pRz$ was chosen arbitrarily, we have $\lim_{N\to\infty} R_N(\kdiamond) = 0$.
\end{proof}

\begin{remark}[Oracle Choice]\label{rem:lifumu:oracle:known}
	\mbox{}\\
	The optimal choice $\kdiamond$ as derived in \cref{co:lifumu:known:oracle} depends on the unknown Mellin transform $\MelX$ and hence is computationally infeasible. 
	Thus, $\kdiamond$ is also called \textit{oracle} (or \textit{oracle choice}) and the whole inequality in \cref{co:lifumu:known:oracle} is called \textit{oracle-type inequality}.
	As an alternative that does not rely on the knowledge of $\MelX$, we consider in \cref{sec:lifumu:datadriven:known} a data driven selection method for the tuning parameter $k$.
\end{remark}

We now investigate rates of convergence that result from the oracle-type inequality derived in \cref{co:lifumu:known:oracle}.
To that end, we specify the regularity of the unknown probability density functions $f$ and the behavior of the Mellin transform $\MelU(t)$ in the asymptotic regime $\abs{t}\to\infty$.
Recall the Mellin-Sobolev spaces as defined in \cref{def:ordinaryandexponentially}.

\begin{assumption}[Regularity of $f$ and decay of $\MelU$]\label{ass:regularity:lifumu}
	\mbox{}\\
	Grant \cref{ass:lifumu:A1}. Let $\rms,\rmb\colon\IR\to\pRz$ be symmetric, positive and continuous weight functions, such that $\rmL\rms^\dagger\in\IL^2$ and $\rmL\rmb\indk{k}\in\IL^2$ for all $k\in\pRz$.
	Let $\rmr_2\in\pRz$ and $\rmr_\infty\in\IR_{\geq 1}$.
	Assume $f\in H_{\rms}(\rmr_2)$ and $f^\varepsilon\in H^\infty_{\rmb}(\rmr_\infty)$.
\end{assumption}
	
\begin{corollary}[Local Estimation - Maximal Inequality - Known Error Distribution]\label{lem:lifumu:known:maximialinequality}
	\mbox{}\\
	Let $N\in\IN$. Grant \cref{ass:regularity:lifumu}. We define the \textit{local maximal risk}
	\begin{align*}
		\forall k\in\pRz:\; \maxrisklocal\br{k|\rms,\rmb,\rmL}&\coloneqq \norm{\rmL\indkco{k}\rms^\dagger}_{\IL^2}^2+\frac{1}{N}\norm{\rmL\indk{k}\rmb}_{\IL^2}^2. 
	\end{align*}
	Then there exists an optimal tuning parameter $\kdiamond = \kdiamond(N)\in\pRz$, such that
	\begin{align*}
		\IE_f^N\brE{\abs{\linbar_{\kdiamond}-\lin(f)}^2} 
		&\leq\Constfeps\ConstfY(\rmr_2^2\lor\rmr_\infty^2)\times\inf_{k\in\pRz}\maxrisklocal\br{k|\rms,\rmb,\rmL}.
	\end{align*}
\end{corollary}

\begin{proof}[Proof of \cref{lem:lifumu:known:maximialinequality}]
	We first note, by definition, \cref{ass:regularity:lifumu} implies 
	\begin{equation*}
		\norm{\rms\MelX}_{\IL^2}^2 \leq \rmr_2^2 \quad\and\quad \abs{\MelUrez(t)} \leq \rmr_\infty\rmb(t)\quad\text{for $\lebesgue$-almost all $t\in\IR$.}
	\end{equation*}
	Thus, an application of the Cauchy-Schwarz inequality yields
	\begin{align*}
		\norm{\rmL\MelX\indkco{k}}_{\IL^1}^2 & 
			 = \br{\int_{\IR}\abs{\rmL(t)}\indkco{k}\abs{\MelX(t)}\rms(t)\rms^\dagger(t)\rmd\lebesgue(t)}^2\\
			& \leq \norm{\rms\MelX}_{\IL^2}^2\times\norm{\rmL\rms^\dagger\indkco{k}}_{\IL^2}^2
			\leq \rmr_2^2\times\norm{\rmL\indkco{k}\rms^\dagger}_{\IL^2}^2.
	\end{align*}
	Further,
	\begin{align*}
		\frac{1}{N}\norm{\rmL\indk{k}\MelUrez}_{\IL^2}^2 &\leq \frac{1}{N}\int_{\IR}\abs{\rmL(t)}^2\indk{k}(t)\abs{\MelUrez(t)}^2\rmd\lebesgue(t)\leq \rmr_\infty^2\times \frac{1}{N}\norm{\rmL\indk{k}\rmb}_{\IL^2}^2.
	\end{align*}
	As we can now follow the proof of \cref{co:lifumu:known:oracle}, the details are omitted.
\end{proof}

To obtain rates of convergence, we further specify the behavior of the 
representant $\rmL(t)$ in the asymptotic regime $\abs{t}\to\infty$.

\begin{assumption}[Asymptotic Behaviour of $\rmL$]\label{lifumu:representant:decay:ass}
	\mbox{}\\
	Let  $p\in\IR$. 
	There exist finite constants $\ConstLlower,\ConstLupper\in\pRz$ such that 
	\begin{equation}\label{lifumu:eq:decay:L}
		\forall t\in\IR:\;\ConstLlower \times \br{1+\abs{t}^2}^{\frac{p}{2}} \leq\abs{\rmL(t)}\leq\ConstLupper\times	\br{1+\abs{t}^2}^{\frac{p}{2}} \tag{$\rmL$}.
	\end{equation}
\end{assumption}

In the following illustration we summarize rates of convergence for the local maximal risk as derived in \cref{lem:lifumu:known:maximialinequality}
by considering scenarios in which $f$ and $f^\varepsilon$ are \textit{ordinary smooth} or \textit{super smooth}, for parameters $b,s\in\pRzz$, where we refer to \cref{def:ordinaryandexponentially}.

\begin{illustration}[Rates of Convergence - Local Maximal Risk - KED]\label{example:lifumu:known:oracle:rates}
\mbox{}\\
Grant \cref{lifumu:representant:decay:ass}.
\Cref{tab:lifumu:oracle:known} provides the order of the \textit{local maximal risk}
\begin{equation*}
	\maxrisklocal\coloneqq \inf_{k\in\pRz}\maxrisklocal \br{k|\rms,\rmb,\rmL}
\end{equation*}
as derived in \cref{lem:lifumu:known:maximialinequality} in the asymptotic regime $N\to\infty$. 
Throughout assume $s>p + \frac{1}{2}$.
The derivation of the presented rates of convergence can be found in \cref{sec:app:lifumu:known:oracle}.
\begin{table}[H]
	\centering
	\footnotesize
	\setlength{\tabcolsep}{6pt}        
	\renewcommand{\arraystretch}{1} 
	\begin{tabular}{@{}
		C{0.1cm}
		C{2.0cm}!{\vrule width 0.6pt}
		C{0.2cm}
		C{0.2cm}!{\vrule width 0.6pt}
		C{2.4cm}
		C{2.4cm}
		C{2cm}
		C{2.4cm}
	  @{}}

	\toprule
	&\ref{lifumu:eq:decay:L}& $f$ & $f_\varepsilon$ &
	$\norm{\rmL\indkco{k}\rms^\dagger}_{\IL^2}^2$ &
	$\norm{\rmL\indk{k}\rmb}_{\IL^2}^2$ &
	$\kdiamond$ & $\maxrisklocal$\\
	\midrule
	\end{tabular}
	\renewcommand{\arraystretch}{3} 
	\begin{tabular}{@{}
		C{0.1cm}
		C{2.0cm}!{\vrule width 0.6pt}
		C{0.2cm}
		C{0.2cm}!{\vrule width 0.6pt}
		C{2.4cm}
		C{2.4cm}
		C{2cm}
		C{2.4cm}
	  @{}}
	\midrule
	(1)
	&$p+b > -\frac{1}{2}$
	& \textsc{os} & \textsc{os}
	& $\displaystyle k^{-2(s-p)+1}$
	& $\displaystyle k^{2(p+b)+1}$
	& $\displaystyle N^{\frac{1}{2(s+b)}}$
	& $\displaystyle N^{-\frac{2(s-p)-1}{2(s+b)}}$
	\\
	
	(2)
	& $p+b=-1/2$
	&  & 
	& 
	& $\displaystyle \ln(k)$
	& $\displaystyle  N^{\frac{1}{2(s-p)-1}}$
	& $\displaystyle \frac{\ln(N)}{N}$
	\\
	
	(3)
	& $p+b<-\frac{1}{2}$
	&  & 
	& 
	& $\displaystyle 1$
	& $\displaystyle N^{\frac{1}{2(s-p)-1}}$
	& $\displaystyle \frac{1}{N}$
	\\
	
	\midrule
	
	(4)
	&$\displaystyle s>p+\frac{1}{2}$
	& \textsc{os}
	& \textsc{es}
	& $\displaystyle k^{-2(s-p)+1}$
	& $\displaystyle k^{(2(p+b)+1)_+}\rme^{k^{2b}}$
	& $\displaystyle \ln(N)^{\frac{1}{2b}}$
	& $\displaystyle \ln(N)^{-\frac{2(s-p)-1}{2b}}$
	\\
	
	\midrule
	
	(5)
	&$p+b>-\frac{1}{2}$
	& \textsc{es} & \textsc{os}
	& $\displaystyle \rme^{-k^{2s}}$
	& $\displaystyle k^{2(p+b)+1}$
	& $\displaystyle \ln(N)^{\frac{1}{2s}}$
	& $\displaystyle \frac{\ln(N)^{\frac{2(p+b)+1}{2s}}}{N}$
	\\
	
	(6)
	&$p+b = - \frac{1}{2}$
	&  & 
	& 
	& $\displaystyle \ln(k)$
	& $\displaystyle \ln(N)^{\frac{1}{2s}}$
	& $\displaystyle \frac{\ln\ln(N)}{N}$
	\\
	
	(7)
	&$p+b<-\frac{1}{2}$
	&  & 
	& 
	& $\displaystyle 1$
	& $\displaystyle \ln(N)^{\frac{1}{2s}}$
	& $\displaystyle \frac{1}{N}$
	\\
	
	\bottomrule
	\end{tabular}
	
	\captionof{table}{Order of the local maximal risk $\maxrisklocal$ in the asymptotic regime $N\to\infty$ when the error distribution is known.}
	\label{tab:lifumu:oracle:known}
	\end{table}
\end{illustration}

We briefly discuss that the representants $\rmL$ as introduced in \cref{example:Illustration1} satisfy \cref{lifumu:representant:decay:ass}.

\begin{example}[\cref{example:Illustration1} Continued]\label{example:lifumu:representants:assumptions:statisfied}
	\mbox{}\\
	Let $x_o\in\pRz$ be fixed.
	\begin{itemize}
		\item \assuRef{Probability Density Function}:
			\cref{lifumu:representant:decay:ass} is satisfied with $p=0$ as $\abs{\rmL_f(t)} = x_o^{-c}$ for all $t\in\IR$.

		\item\assuRef{Derivative of the Probability Density Function}:
		 \cref{lifumu:representant:decay:ass} is satisfied with $p=\beta$ as $\abs{\rmL_{\rmD^\beta f}(t)} = x_o^{-c-\beta}\prod_{j\in\nset{\beta}}\br{(c+\beta-j)^2+4\pi^2t^2}^{\frac{1}{2}}$ for all $t\in\IR$.

		\item \assuRef{Survival Function and CDF}:
		\cref{lifumu:representant:decay:ass} is satisfied with $p=-1$ as $\abs{\rmL_{\SX}(t)} = x_o^{-c+1}\br{(c-1)^2+4\pi^2 t^2}^{-\frac{1}{2}}$ for all $t\in\IR$.
		Similarly, for $\rmL_{F^X}$.
	\end{itemize}
\end{example}

The following corollary provides rates of convergence for estimating $f(x_o)$, $\rmD^\beta f(x_o)$, and $\SX(x_o)$, respectively.
To this end, we consider the spectral cut-off estimators $\widetilde{f}_{k}(x_o)$, $\widetilde{(\rmD^\beta f)}_{k}(x_o)$, and $\SXbar_{k}(x_o)$, obtained by inserting $\rmL_f$, $\rmL_{\rmD^\beta f}$, and $\rmL_{\SX}$, respectively, into the generic formula in \cref{def:linearfunctionalestimator:known}; see also \cref{example:Illustration1}.

\begin{corollary}[Local Estimation - Rates of Convergence - KED]\label{lifumu:maximalinequalities:known:oracle:examples}
	\mbox{}\\
	Grant \cref{ass:regularity:lifumu} with $\rms \coloneqq \rmsos$ for $s\in\pRzz$.
	Further, let $x_o\in\pRz$.
	\begin{enumerate}
		\item\label{lifumu:density:known:rates}
				 Let $\rmb = \rmbos$ with $b>-\frac{1}{2}$ and $s>\frac{1}{2}$. If we choose $\kdiamond \coloneqq\kdiamond(N)\coloneqq N^{\frac{1}{2(s+b)}}$, we obtain
						\begin{equation*}
							 \IE_{f}^{N}\brE{\abs{\widetilde{f}_{\kdiamond}(x_o)-f(x_o)}^2} \lesssim N^{-\frac{2s-1}{2(s+b)}}.
						\end{equation*}

				 Let $\rmb = \rmbes$ with $b\in\pRz$ and $s>\frac{1}{2}$. If we choose $\kdiamond \coloneqq\kdiamond(N)\coloneqq \ln(N)^{\frac{1}{2b}},$ we obtain
						\begin{equation*}
							 \IE_{f}^{N}\brE{\abs{\widetilde{f}_{\kdiamond}(x_o)-f(x_o)}^2} \lesssim \ln(N)^{-\frac{2s-1}{2b}}.
						\end{equation*}
		\item 
				 Let $\beta\in\IN$. Let $\rmb = \rmbos$ with $b>-\beta-\frac{1}{2}$ and $s>\beta + \frac{1}{2}$. If we choose $\kdiamond \coloneqq\kdiamond(N)\coloneqq N^{\frac{1}{2(s+b)}}$, we obtain
					\begin{equation*}
						 \IE_{f}^{N}\brE{\abs{\widetilde{(\rmD^\beta f)}_{\kdiamond}(x_o)-\rmD^\beta f(x_o)}^2} \lesssim N^{-\frac{2(s-\beta)-1}{2(s+b)}}.
					\end{equation*}
				
				 Let $\beta\in\IN$. Let $\rmb = \rmbes$ and $s>\beta + \frac{1}{2}$. If we choose $\kdiamond \coloneqq\kdiamond(N)\coloneqq \ln(N)^{\frac{1}{2b}}$, we obtain
					\begin{equation*}
						\IE_{f}^{N}\brE{\abs{\widetilde{(\rmD^\beta f)}_{\kdiamond}(x_o)-\rmD^\beta f(x_o)}^2} \lesssim \ln(N)^{-\frac{2(s-\beta)-1}{2b}}.
					\end{equation*}

		\item 
				 Let $\rmb = \rmbos$ with $b>\frac{1}{2}$ and $s>\frac{1}{2}$. If we choose $\kdiamond \coloneqq\kdiamond(N)\coloneqq N^{\frac{1}{2(s+b)}}$, we obtain
						\begin{equation*}
							 \IE_{f}^{N}\brE{\abs{\SXbar_{\kdiamond}(x_o)-\SX(x_o)}^2} \lesssim N^{-\frac{2s+1}{2(s+b)}}.
						\end{equation*}
				 Let $\rmb = \rmbes$ and $s>\frac{1}{2}$. If we choose $\kdiamond \coloneqq\kdiamond(N)\coloneqq \ln(N)^{\frac{1}{2b}}$, we obtain
					\begin{equation*}
						 \IE_{f}^{N}\brE{\abs{\SXbar_{\kdiamond}(x_o)-\SX(x_o)}^2} \lesssim \ln(N)^{-\frac{2s+1}{2b}}.
					\end{equation*}
	\end{enumerate}
	Here, all inequalities are understood in the asymptotic regime $N\to\infty$.
\end{corollary}

\begin{proof}[Proof of \cref{lifumu:maximalinequalities:known:oracle:examples}]
	The proof is a direct consequence of \cref{lem:lifumu:known:maximialinequality}, \cref{example:lifumu:known:oracle:rates} and \cref{example:lifumu:representants:assumptions:statisfied}.
	The details are thus omitted. We note that the appearing multiplicative constant $\const\in\pRz$ depends on
	the parameters and constants $\rmr_2,\rmr_\infty,x_o,b,s, f^\varepsilon$, and $c$. In case of ii) also on $\beta$.
\end{proof}

\begin{remark}[Minimax Optimality]
	\mbox{}\\
	The rates of convergence obtained for the probability density function estimation in \cref{lifumu:maximalinequalities:known:oracle:examples} i)
			are proven to be optimal in the minimax sense in \cite{Migueletal_2023}.
			Moreover, the rates are obtained for probability density functions $f$,
			which are $s$-regular in the Mellin-Sobolev sense. 
			The rates of convergence obtained for the survival function in \cref{lifumu:maximalinequalities:known:oracle:examples} ii) are also proven to be minimax-optimal in \cite{BrennerMiguelPhandoidaen_2022}.
			Keep in mind that analogous rates of convergence can be derived for estimating the CDF $F^X$.
\end{remark}

\subsection{Data-Driven Model Selection}\label{sec:lifumu:datadriven:known}

Following \cref{rem:lifumu:oracle:known}, the optimal choice $\kdiamond$ leading to optimal rates of convergence of the local maximal risk in the minimax sense 
is typically infeasible from a computational point of view.
In this section we provide a fully data-driven selection method for the tuning parameter $k$,
which does not necessitate any additional knowledge of the unknown Mellin transform $\MelX$.
In contrast to the selection method in \cite{Migueletal_2021}, the proposed penality below makes the whole analysis also applicable in 
the statistical regime \assuRef{SRII}.\\

Heuristically speaking, an optimal choice of the tuning parameter $k$
is supposed to minimize the MSE
\begin{align*}
	\pRz\ni k\mapsto \IE_f^N\brE{\abs{\linbar_k-\lin(f)}^2}  &= \abs{\IE_f^N\brE{\linbar_k}-\lin(f)}^2+ \IE_f^N\brE{\abs{\linbar_k-\IE_f^N\brE{\linbar_k}}^2}\\
		& \approx \abs{\lin_k^{\rmwU}-\lin(f)}^2 +  \frac{1}{N}\norm{\rmL\indk{k}\MelUrez}_{\IL^2}^2
\end{align*}
derived from a squared bias-variance trade-off. As the 
squared bias term depends on the unknown $\MelX$, the main idea is, given a finite family $\magenta{\calK\subseteq \pRz}$ of models, to
define a function $B\colon\calK\to\IR$ imitating the bias
and to
select
	\begin{equation*}
		\kbar:\in\argmin_{\magenta{k\in\calK}}\brK{B(k)+P(k)},
	\end{equation*}
	where $P\colon\calK\to\IR$ is an additional penalty to avoid overfitting.
	\magenta{In the subsequent we further define for any $k\in\calK$ the subset $\calK_{>k}\coloneqq \brK{j\in\calK:\; j>k}$.
	Analogously $\calK_{\geq k}$ is defined.}\\

Concretely, given a family of penalties $(\pen_k)_{k\in\calK}\subseteq\pRzz$,
we define
\begin{align}
	\forall k\in\magenta{\calK}:\; B(k)\coloneqq \contrastbar_k &\coloneqq \max_{j\in\magenta{\calK_{>k}}}\left\{\left(|\linbar_j-\linbar_k|^2-6\pen_j-6\pen_k\right)_+\right\} \notag    \\
	&= \max_{j\in\magenta{\calK_{>k}}}\brK{\br{|\linbar_{k\lor j}-\linbar_{k\land j}|^2-6\pen_{k\lor j}-6\pen_{k\land j}}_+}\label{eq:lifumu:known:contrast}.
\end{align}
 Choosing $P(k)\coloneqq 6\pen_k$ for each $k\in\magenta{\calK}$, we define a data-driven selection method of the tuning parameter $k$ by 
\begin{equation}\label{eq:lifumu:known:datadriven}
	\kbar:\in \argmin_{k\in\magenta{\calK}}\brK{\contrastbar_k +6\pen_k}.
\end{equation}
As $\kbar$ is implicitly defined by utilizing the estimator $\linbar_k$,
$\kbar$ itself defines a random variable depending on the data $(Y_i)_{i\in\nset{N}}\sim\IP_{f}^N$.
In what follows, we derive an upper bound for the local \textit{data-driven} risk
\begin{equation*}
	\IE_{f}^N\brE{\abs{\linbar_{\kbar}-\lin(f)}^2},
\end{equation*}
and show that selecting the tuning parameter $k$ through the data-driven selection method above
does not change the rates of convergence derived above (up to constants).
We start with a \textit{Key Argument} to find a first upper bound of $$\abs{\linbar_{\kbar}-\lin(f)}^2.$$

\begin{lemma}[Local Estimation - Key Argument I - Known Error Distribution]\label{lem:lifumu:known:KeyArgument}
	\mbox{}\\
	For any family $(\pen_k)_{k\in\magenta{\calK}}\subseteq\IR_{\geq 0}$ and \magenta{$\calK$ as above}, we select the tuning parameter $\kbar$ as in \cref{eq:lifumu:known:datadriven}.
	Defining
	\begin{align}
		\forall k\in\magenta{\calK}:\;\rmZbar_k \coloneqq \linbar_k - \lin_k^{\rmwU} &= \int_{[-k,k]}\rmL(-t)\br{\MelXtilde(t)-\MelX(t)}\mathrm{d}\lebesgue(t)\nonumber\\
			&= \int_{[-k,k]}\rmL(-t)\MelUrez(t)\br{\MelYhat(t)-\MelY(t)}\mathrm{d}\lebesgue(t), \label{eq:lifumu:known:empiricalprocess}
	\end{align}
	we obtain
	\begin{align*}
		|\linbar_{\kbar}-\lin(f)|^2 &\leq \min_{k\in\magenta{\calK}}\brK{12\norm{\rmL\MelX\mathds{1}_{[-k,k]^C}}_{\IL^1}^2+24\pen_k + 2|\linbar_k-\lin(f)|^2} \\
		&\hspace{1cm}+ 24\max_{k\in\magenta{\calK}}\brK{\br{|\rmZbar_k^{\rmwU}|^2-2\pen_k}_+},
	\end{align*}
	where the inequality is understood to hold $\IP_f^N$-almost surely.
\end{lemma}

\begin{proof}[Proof of \cref{lem:lifumu:known:KeyArgument}]
	Fix any realization of $\linbar_{\kbar}$.
	From the definition of the contrast function $\contrastbar_k$ as well as of the definition of $\kbar$ in \cref{eq:lifumu:known:datadriven}, we deduce for all $k\in\magenta{\calK}$
	\begin{align*}
		\abs{\linbar_{\kbar}-\lin(f)}^2 
		& \leq 2\br{\abs{\linbar_{\kbar\lor k}-\linbar_{\kbar\land k}}^2  \mp 6\pen_{\kbar\land k} \mp 6\pen_{\kbar\lor k} +\abs{\linbar_{k}-\lin(f)}^2}\\
		& \leq 2 \br{\contrastbar_{\kbar\land k} + 6\pen_{\kbar\land k} + \contrastbar_{\kbar\lor k} + 6\pen_{\kbar\lor k}  +\abs{\linbar_{k}-\lin(f)}^2}\\
		& \leq 4 \br{\contrastbar_{k} + 6\pen_{k}} + 2\abs{\linbar_{k}-\lin(f)}^2.
	\end{align*}
	Further, for all $k\in\magenta{\calK}$ and $j\in\magenta{\calK_{>k}}$ we have 
	\begin{align*}
		\abs{\lin_j(f)-\lin_k(f)} 
			 = \abs{\int_{\IR}\rmL(-t)\MelX(t)\1_{[-j,-k)\cup(k,j]}(t)\rmd\lebesgue(t)}
			& \leq \norm{\rmL\MelX\indkco{k}}_{\IL^1}.
	\end{align*}
	Thus,
	\begin{align*}
		\contrastbar_k &= \max_{j\in\magenta{\calK_{>k}}}\brK{\br{|\linbar_{ j}-\linbar_{k}|^2-6\pen_{j}-6\pen_{k}}_+}\\
		& \leq \max_{j\in\magenta{\calK_{>k}}}\brK{\br{3|\rmZbar_j^{\rmwU}|^2+3|\rmZbar_k^{\rmwU}|^2+ 3 \norm{\rmL\MelX\indkco{k}}_{\IL^1}^2-6\pen_{j}-6\pen_{k}}_+}\\
		& \leq \max_{j\in\magenta{\calK_{>k}}}\brK{\br{3|\rmZbar_j^{\rmwU}|^2-6\pen_{j}}_+} + \br{3|\rmZbar_k^{\rmwU}|^2-6\pen_{k}}_+  + 3 \norm{\rmL\MelX\indkco{k}}_{\IL^1}^2\\
		& \leq 6\max_{j\in\magenta{\calK}}\brK{\br{|\rmZ_j^{\rmwU}|^2-2\pen_{j}}_+} + 3 \norm{\rmL\MelX\indkco{k}}_{\IL^1}^2.
	\end{align*}
	Together with above computations, we obtain
	\begin{align*}
		\abs{\linbar_{\kbar}-\lin(f)}^2 & \leq 12\norm{\rmL\MelX\mathds{1}_{[-k,k]^C}}_{\IL^1}^2+24\pen_k + 2|\linbar_k-\lin(f)|^2\\
			&\hspace{1cm} + 24\max_{j\in\magenta{\calK}}\brK{\br{|\rmZ_j^{\rmwU}|^2-2\pen_{j}}_+}.
	\end{align*}
	As the last bound holds for arbitrary $k\in\magenta{\calK}$, we apply the minimum on both sides, which shows the claim.
\end{proof}

Exploiting the key argument derived in \cref{lem:lifumu:known:KeyArgument} to find an upper bound of the data driven risk, 
we need to choose the
family of penalties $(\pen_k)_{k\in\pRz}$ in such a way that 
\begin{equation*}
	\IE_f^N\brE{\max_{k\in\magenta{\calK}}\brK{\br{|\rmZbar_k^{\rmwU}|^2-2\pen_k}_+}}
\end{equation*}
is bounded appropriately. This is part of the next assertion.
To that end, we impose some additional moment assumptions on the observables $(Y_i)_{i\in\nset{N}}$.
\begin{assumption}[Data-driven Selection - Known Error Distribution]\label{ass:lifumu:A2}
	\mbox{}
	\begin{itemize}
		\item Grant \cref{ass:lifumu:A1}. Assume in addition, that there exists $\rmp\in \IR_{\geq 2}$, such that $\sigmabar_{Y,2\rmp}^2\coloneqq 1+\IE[Y^{2\rmp(c-1)}]  \in\pRz$,
	i.e. $f_Y\in\IL_{\rmx^{2\rmp(c-1)}}^1$.

	\item \magenta{Let $\kappa\in\pRz$, and choose $\calK\subseteq \kappa\IN\coloneqq\brK{\kappa l:\; l\in\IN}$.}
	\end{itemize}
\end{assumption}
\begin{notation}[Local Estimation - Model Selection]\label{not:lifumu:datadriven:stuff}
	\mbox{}\\
	Fix $N\in\IN$ and $\rmp\in\IR_{\geq 2}$. Let $\rmv,\rmw\colon\IR\to\IR$ be weight functions chosen in such a way that $\rmL\rmv\indk{k},\rmL\rmw\indk{k}\in\IL^2$ for all $k\in\IR_{\geq \magenta{\kappa}}$. 
	Let
	\begin{equation}\label{eq:notation:lifumu:Delta}
		\forall k\in\IR_{\geq \magenta{\kappa}}:\; \Delta_k^\rmv\coloneqq  2\magenta{\br{\frac{k}{\kappa}}^2}\br{1\lor\magenta{k}\norm{\rmL\rmv\mathds{1}_{[-k,k]}}_{\IL^2}^2}\in\IR_{\geq 1}.
	\end{equation}
	Further, we define
	\begin{equation}\label{eq:notation:lifumu:rho}
		\forall k\in\IR_{\geq \magenta{\kappa}}:\;\rho_N\br{k|\rmv}\coloneqq\br{1\lor\br{ \frac{2k}{N}\cdot \br{\Delta_k^\rmv}^{\frac{1}{\rmp-1}}\br{1\lor\ln\br{\Delta_k^\rmv}}}} \br{1\lor\ln\br{\Delta_k^\rmv}}\in\IR_{\geq 1}.
	\end{equation}
	Evidently, if there is some $A\in\IR_{\geq 1}$, such that $\abs{\rmv(t)} \leq A\times \abs{\rmw(t)}$ for all $t\in[-k,k]$,
	elementary computations show 
	\begin{equation}\label{eq:lifumu:rho:calculus}
		\forall k\in\IR_{\geq \magenta{\kappa}}:\; \Delta_k^\rmv \leq A^2\times \Delta_k^\rmw\quad\text{and}\quad \br{1\lor\ln\br{\Delta_k^\rmv}} \leq \br{1+\ln\br{A^2}}\times\br{1\lor\ln\br{\Delta_k^{\rmw}}},
	\end{equation}
	as well as 
	\begin{equation}\label{eq:lifumu:rho:calculus2}
		\forall k\in\IR_{\geq \magenta{\kappa}}:\; \rho_N(k|\rmv) \leq \br{1+\ln(A^2)}^2A^{\frac{2}{\rmp-1}}\times\rho_N\br{k|\rmw} \leq \frac{4}{\rme}\times A^4\times\rho_N\br{k|\rmw}.
	\end{equation}
	Let us next define
	\begin{equation*}
		\forall k\in\IR_{\geq \magenta{\kappa}}:\; \nu_k^\rmv\coloneqq \br{\frac{1}{128}\br{\Delta_k^\rmv}^\frac{1}{\rmp-1}}^\frac{1}{2}\in\pRz
	\end{equation*}
	as well as 
	\begin{equation*}
		\forall k\in\IR_{\geq \magenta{\kappa}}:\; \rmV_k^\rmv \coloneqq \IE\brE{\mathds{1}_{(-\infty,\nu_k^\rmv]}(Y_1^{c-1})\abs{\int_{[-k,k]}\rmL(-t)\rmv(t)Y_1^{c-1+2\pi\iota t}\mathrm{d}\lebesgue(t)}^2}\in\pRzz.
	\end{equation*}
	By an application of \cref{lem:MellinTrick}, we have for all $k\in\IR_{\geq \magenta{\kappa}}$
	\begin{equation}\label{eq:lifumu:Vkv:estimate}
		 \rmV_k^\rmv \leq \IE\brE{\abs{\int_{[-k,k]}\rmL(-t)\rmv(t)Y_1^{c-1+2\pi\iota t}\mathrm{d}\lebesgue(t)}^2} \leq  \Constfeps\ConstfY\times\norm{\rmL\rmv\indk{k}}_{\IL^2}^2.
	\end{equation}
\end{notation}

\begin{corollary}[Local Estimation - Concentration of $\rmZ_k^{\rmwU}$ - Known Error Distribution]\label{co:lifumu:known:Concentration}
	\mbox{}\\
	Grant \cref{ass:lifumu:A2}. For $k\in\IR_{\geq \magenta{\kappa}}$ define $\Delta_k^{\rmwU}$, $\nu_k^{\rmwU}$ and $\rmV_k^{\rmwU}$ as in \cref{not:lifumu:datadriven:stuff} with  $\rmwU\coloneqq \MelUrez$.
	Further, set
	\begin{equation*}
		\forall k\in\IR_{\geq \magenta{\kappa}}:\;\penbar_k^{\rmwU}\coloneqq 8\br{\rmV_k^{\rmwU}+ 2k\norm{\rmL\MelUrez\indk{k}}_{\IL^2}^2\frac{\br{\Delta_k^{\rmwU}}^{\frac{1}{\rmp-1}}\ln(\Delta_k^{\rmwU})}{N}} \frac{\ln(\Delta_k^{\rmwU})}{N}.
	\end{equation*}
	We then have for all $N\in\IN$
	\begin{equation*}
		\IE_f^N\brE{\max_{k\in\magenta{\calK}}\brK{\br{|\rmZbar_k^{\rmwU}|^2-2\penbar_k^{\rmwU}}_+}} \leq \frac{1}{\magenta{\kappa^2}N}\times \br{40\ConstfY\Constfeps+\frac{128^\rmp}{32}\sigmabar_{Y,2\rmp}^2}.
	\end{equation*}
\end{corollary}
\begin{proof}[Proof of \cref{co:lifumu:known:Concentration}]
	Defining $\rmw\coloneqq \rmL\MelUrez$, the proof follows directly from \cref{propo:lifumu:known:concentration:process:app}.
\end{proof}

The last Corollary provides a first proposal for a suitable family of penalties $\br{\penbar_k^{\rmwU}}_{k\in\magenta{\calK}}$ in order to control
the local data-driven risk
\begin{equation*}
	\IE_f^N\brE{\abs{\linbar_{\kbar}-\lin(f)}^2}.
\end{equation*}
However, $\br{\penbar_k^{\rmwU}}_{k\in\magenta{\calK}}$ does depend
on the unknown distribution of $Y$ via the term $\rmV_k^{\rmwU}$.
In other words, $\kbar$ is not fully known and thus computationally infeasible.
Having the data $(Y_i)_{i\in\nset{N}}$ at our disposal,
we replace the term $\rmV_k^{\rmwU}$ by its empirical counterpart $\rmVhat_k^{\rmwU}$ defined by
\begin{equation*}
	\forall k\in\magenta{\IR_{\geq\kappa}}:\;\rmVhat_k^{\rmwU}\coloneqq \frac{1}{N}\sum_{i\in\nset{N}}\mathds{1}_{(-\infty,\nu_k^{\rmwU}]}(Y_i^{c-1})\left|\int_{[-k,k]}\rmL(-t)\MelUrez(t)Y_i^{c-1+2\pi\iota t}\mathrm{d}\lebesgue(t)\right|^2.
\end{equation*}
Clearly, $\rmVhat_k^{\rmwU}$ is an unbiased estimator of $\rmV_k^{\rmwU}$, i.e.
$$\forall k\in\IR_{\geq \magenta{\kappa}}:\;\IE_f^N\brE{\rmVhat_k^{\rmwU}} = \rmV_k^{\rmwU}.$$
We now define an \textit{estimated penalty} for each $k\in\IR_{\geq \magenta{\kappa}}$ by 
\begin{equation}\label{eq:ligfumu:known:penhat}
	\penhat_k^{\rmwU}\coloneqq 8\br{\orange{2\rmVhat_k^{\rmwU}}+\orange{3}\times 2k\norm{\rmL\MelUrez\indk{k}}_{\IL^2}^2\frac{\br{\Delta_k^{\rmwU}}^{\frac{1}{\rmp-1}}\ln(\Delta_k^{\rmwU})}{N}} \frac{\ln(\Delta_k^{\rmwU})}{N} \in\pRzz,
\end{equation}
which is fully known in the statistical regime \assuRef{SRI}.
Consequently, a \textit{fully data-driven} model selection procedure is given by 
\begin{equation}\label{eq:lifumu:khat:known:fullydatadriven}
	\khat:\in \argmin_{k\in\magenta{\calK}}\brK{\contrasthat_k^{\rmwU} +6\penhat_k^{\rmwU}},
\end{equation}
where
\begin{equation*}
	\contrasthat_k^{\rmwU} \coloneqq\max_{j\in\magenta{\calK_{>k}}}\brK{\br{|\linbar_j-\linbar_k|^2-6\penhat_j^{\rmwU}-6\penhat_k^{\rmwU}}_+}.
\end{equation*}

\begin{lemma}[Local Estimation - Key Argument II - Known Error Distribution]\label{lem:lifumu:riskbound:dd:known}
	\mbox{}\\
	Grant \cref{ass:lifumu:A2}.
	Select $\khat$ as in \cref{eq:lifumu:khat:known:fullydatadriven}. We then have
	\begin{align*}
		\abs{\linbar_{\khat}-\lin(f)}^2 &\leq \min_{k\in\magenta{\calK}}\left\{12\norm{\rmL\MelX\mathds{1}_{[-k,k]^C}}_{\IL^1}^2+24\penhat_k^{\rmwU} + 2|\linbar_k-\lin(f)|^2\right\}\\
		& + 24\max_{k\in\magenta{\calK}}\left\{\left(|\rmZbar_k^{\rmwU}|^2-2\penbar_k^{\rmwU}\right)_+\right\} + 48\max_{k\in\magenta{\calK}}\left\{\left(\penbar_k^{\rmwU}-\penhat_k^{\rmwU}\right)_+\right\},
	\end{align*}
	where the inequality is understood to hold $\IP_f^N$-almost surely.
\end{lemma}

\begin{proof}[Proof of \cref{lem:lifumu:riskbound:dd:known}]
	Observing that
	\begin{align*}
		\forall k\in\magenta{\calK}:\;\left(|\rmZbar_k^{\rmwU}|^2-2\penhat_k^{\rmwU}\right)_+ \leq \left(|\rmZbar_k^{\rmwU}|^2-2\penbar_k^{\rmwU}\right)_+ + 2\left(\penbar_k^{\rmwU}-\penhat_k^{\rmwU}\right)_+,
	\end{align*}
	and utilizing \cref{lem:lifumu:known:KeyArgument}, we immediately obtain
	\begin{align*}
	|\linbar_{\kbar}-\lin(f)|^2 &\leq \min_{k\in\magenta{\calK}}\left\{12\norm{\rmL\MelX\mathds{1}_{[-k,k]^C}}_{\IL^1}^2+24\penhat_k^{\rmwU} + 2|\linbar_k-\lin(f)|^2\right\}\\
	& + 24\max_{k\in\magenta{\calK}}\left\{\left(|\rmZbar_k^{\rmwU}|^2-2\penbar_k^{\rmwU}\right)_+\right\} + 48\max_{k\in\magenta{\calK}}\left\{\left(\penbar_k^{\rmwU}-\penhat_k^{\rmwU}\right)_+\right\},
	\end{align*}
	which shows the claim.
\end{proof}

\begin{remark}[Local Estimation - Key Argument II - Known Error Distribution]
	\mbox{}\\
	In contrast to $\penbar_k^{\rmwU}$, the definition of $\penhat_k^{\rmwU}$ is adjusted by a factor $\orange{2}$ when replacing $\rmV_k^{\rmwU}$ by $\rmVhat_k^{\rmwU}$, and by a factor $\orange{3}$
	concerning the remaining term. This is due to the use of Bernstein's inequality (see \cref{lem:BernsteinInequality_general}) in the proof of the following assertions.
\end{remark}

To find an upper bound of the data driven risk $$\IE_f^N\brE{\abs{\linbar_{\kbar}-\lin(f)}^2},$$ it hence suffices to 
find upper bounds for 
\begin{equation*}
	\IE_f^N\brE{\penhat_k^{\rmwU}}\quad\text{and}\quad\IE_f^N\brE{\max_{k\in\magenta{\calK}}\left\{\left(\penbar_k^{\rmwU}-\penhat_k^{\rmwU}\right)_+\right\}},
\end{equation*}
respectively.

\begin{lemma}[Local Estimation - Estimated Penalty - Known Error Distribution]\label{lem:lifumu:penality:known:dd}
	\mbox{}\\
	Let $N\in\IN$. Grant \cref{ass:lifumu:A2}. For $k\in\IR_{\geq \magenta{\kappa}}$ let $\penbar_k^{\rmwU}$ as defined in \cref{co:lifumu:known:Concentration} and $\penhat_k^{\rmwU}$ as defined in \cref{eq:ligfumu:known:penhat}.
	We then have
	\begin{equation*}
		\forall k\in\IR_{\geq \magenta{\kappa}}:\;\IE_f^N\brE{\penhat_k^{\rmwU}} \leq \frac{40}{N}\Constfeps\ConstfY\times \norm{\rmL\MelUrez\indk{k}}_{\IL^2}^2\rho_N(k|\rmwU)
	\end{equation*}
	with $\rho_N(k|\rmwU)$ as defined in \cref{not:lifumu:datadriven:stuff}.
	Moreover,
	\begin{equation*}
		\IE_f^N\brE{\max_{k\in\magenta{\calK}}\left\{\left(\penbar_k^{\rmwU}-\penhat_k^{\rmwU}\right)_+\right\}} \leq \Constfeps\ConstfY\frac{1}{\magenta{\kappa^2}N}.
	\end{equation*}
\end{lemma}

\begin{proof}[Proof of \cref{lem:lifumu:penality:known:dd}]
	This is a direct consequence of \cref{lifumu:unknown:app:estimatedpenality} with $\rmw = \rmL\MelUrez$.
\end{proof}


\begin{proposition}[Local Estimation - Data-Driven Risk - KED]\label{prop:lifumu:neworacle:known}
	\mbox{}\\
	Let $N\in\IN$. Grant \cref{ass:lifumu:A2}. 
	Set
	\begin{equation}\label{eq:lifumu:known:KN:def}
		\magenta{\calK^N}\coloneqq\brK{k\in\kappa\IN:\;\br{1\lor\ln(k)}\norm{\rmL\indk{k}}_{\IL^2}^2 \leq N\norm{\rmL\indk{1}}_{\IL^2}^2}.
	\end{equation}
	Let
	\begin{equation}\label{eq:khat:known:final}
		\khat:\in \argmin_{k\in\magenta{\calK^N}}\brK{\contrasthat_k^{\rmwU} +6\penhat_k^{\rmwU}}
	\end{equation}
	with
	\begin{equation*}
		\contrasthat_k^{\rmwU} \coloneqq\max_{j\in\magenta{\calK^N_{>k}}}\brK{\br{|\linbar_j-\linbar_k|^2-6\penhat_j^{\rmwU}-6\penhat_k^{\rmwU}}_+}
	\end{equation*}
	and $\penhat_k^{\rmwU}$ as defined in \cref{eq:ligfumu:known:penhat}.
	Then, for all $N\in\IN$ large enough such that
	\begin{equation*}
		\norm{\MelU}_{\IL^\infty}^2\min_{k\in\magenta{\kappa\IN}}\brK{\norm{\rmL\MelX\mathds{1}_{[-k,k]^C}}_{\IL^1}^2 + \frac{1}{N}\norm{\rmL\MelUrez\indk{k}}_{\IL^2}^2\rho_N(k|\rmwU)}\leq\norm{\rmL\indk{1}}_{\IL^2}^2,
	\end{equation*}
	we have
	\begin{align*}
		\IE_f^N\left[|\linbar_{\khat}-\lin(f)|^2\right] &\leq \const_1\times\ConstfY\Constfeps \min_{k\in\magenta{\kappa\IN}}\brK{\norm{\rmL\MelX\mathds{1}_{[-k,k]^C}}_{\IL^1}^2 + \frac{1}{N}\norm{\rmL\MelUrez\indk{k}}_{\IL^2}^2\rho_N(k|\rmwU)}\\
			&\hspace{1cm} + \const_2\times \frac{1}{\magenta{\kappa^2}N}\br{\ConstfY\Constfeps+128^\rmp\sigmabar_{Y,2\rmp}^2}
	\end{align*}
	for some numerical constants $\const_1,\const_2\in\pRz$.
\end{proposition}

\begin{proof}[Proof of \cref{prop:lifumu:neworacle:known}]
	The proof is divided in four steps.
	\begin{itemize}
		\item \textbf{Step I:} Firstly, the key argument \cref{lem:lifumu:riskbound:dd:known}
			provides an upper bound for the data-driven risk, which for any $\ko\in\magenta{\calK^N}$ reads as 
			\begin{align*}
				\IE_f^N\brE{\abs{\linbar_{\khat}-\lin(f)}^2} &\leq 12\norm{\rmL\MelX\mathds{1}_{[-\ko,\ko]^C}}_{\IL^1}^2+24\IE_f^N\brE{\penhat_{\ko}^{\rmwU}} + 2\IE_f^N\brE{|\linbar_{\ko}-\lin(f)|^2}\\
					&\hspace{1cm} + 24\IE_f^N\brE{\max_{k\in\magenta{\calK^N}}\left\{\left(|\rmZbar_k^{\rmwU}|^2-2\penbar_k^{\rmwU}\right)_+\right\}}\\
					&\hspace{2cm} + 48\IE_f^N\brE{\max_{k\in\magenta{\calK^N}}\brK{\left(\penbar_k^{\rmwU}-\penhat_k^{\rmwU}\right)_+}}\\
					& \eqqcolon T_1 + T_2 + T_3 + T_4 + T_5.
			\end{align*}
		
		\item \textbf{Step II:} Following \cref{co:lifumu:known:Concentration} and \cref{lem:lifumu:penality:known:dd}, the last two terms are upper bounded by 
		\begin{equation*}
			T_4 + T_5\leq  \const_2\times \frac{1}{\magenta{\kappa^2}N}\br{\ConstfY\Constfeps+128^\rmp\sigmabar_{Y,2\rmp}^2}
		\end{equation*}
		with a numerical constant $\const_2\in\pRz$.
		On the other hand, the first three terms are upper bounded by 
		\begin{align*}
			T_1+T_2+T_3 & \leq \const_1\Constfeps\ConstfY\br{\norm{\rmL\MelX\mathds{1}_{[-\ko,\ko]^C}}_{\IL^1}^2 + \frac{1}{N}\times \norm{\rmL\MelUrez\indk{\ko}}_{\IL^2}^2\rho_N(\ko|\rmwU)}
		\end{align*}
		with a numerical constant $\const_1\in\pRz$ as a consequence of \cref{lem:lifumu:penality:known:dd} and \cref{lem:Lifumu:known:MSE}.

		\item \textbf{Step III:}
			Using similar arguments as in the proof of \cref{co:lifumu:known:oracle} one can show that
			\begin{equation*}
				\magenta{\kappa\IN}\ni k\mapsto \norm{\rmL\MelX\mathds{1}_{[-k,k]^C}}_{\IL^1}^2 + \frac{1}{N}\norm{\rmL\MelUrez\indk{k}}_{\IL^2}^2\rho_N(k|\rmwU)\in\pRz
			\end{equation*}
			attains a global minimum.	
			Thus, we define
			\begin{equation}\label{eq:lifumu:datadriven:known:kochoice}
				\ko:\in\argmin_{k\in\magenta{\kappa\IN}}\brK{\norm{\rmL\MelX\mathds{1}_{[-k,k]^C}}_{\IL^1}^2 + \frac{1}{N}\norm{\rmL\MelUrez\indk{k}}_{\IL^2}^2\rho_N(k|\rmwU)}.
			\end{equation}
			The claim is proven, if $\ko\in\magenta{\calK_N}$.

		\item \textbf{Step IV:}
		We now verify that $\ko\in\magenta{\calK^N}$.
		Keeping \cref{eq:notation:lifumu:Delta} in \cref{not:lifumu:datadriven:stuff} in mind, we 
		have 
		\begin{equation*}
			\ko \leq 2\ko^3 \leq \Delta_{\ko}^{\rmwU}\quad\and\quad \br{1\lor \ln(\ko)} \leq \rho_N(\ko|\rmwU).
		\end{equation*}
		Thus,
		\begin{align*}
			&\norm{\rmL\MelX\mathds{1}_{[-\ko,\ko]^C}}_{\IL^1}^2 + \frac{1}{N}\norm{\rmL\MelUrez\indk{\ko}}_{\IL^2}^2\rho_N(\ko|\rmwU)\\
				& \hspace{4cm}\geq   \frac{1}{N}\norm{\rmL\MelUrez\indk{\ko}}_{\IL^2}^2\times \br{1\lor \ln(\ko)}\\
				& \hspace{4cm}\geq   \frac{1}{N}\norm{\rmL\indk{\ko}}_{\IL^2}^2\times \norm{\MelU}_{\IL^\infty}^{-2} \times \br{1\lor \ln(\ko)}.
		\end{align*}
		By hypothesis, $N$ is chosen sufficiently large, such that
		\begin{equation*}
			\norm{\MelU}_{\IL^\infty}^2\min_{k\in\magenta{\kappa}\IN}\brK{\norm{\rmL\MelX\mathds{1}_{[-k,k]^C}}_{\IL^1}^2 + \frac{1}{N}\norm{\rmL\MelUrez\indk{k}}_{\IL^2}^2\rho_N(k|\rmwU)}\leq\norm{\rmL\indk{1}}_{\IL^2}^2,
		\end{equation*}
		where we note that the minimum decreases for an increasing $N\in\IN$.
		As a consequence
		\begin{equation*}
			\norm{\rmL\indk{1}}_{\IL^2}^2\geq \frac{1}{N}\norm{\rmL\indk{\ko}}_{\IL^2}^2\times \br{1\lor \ln(\ko)}.
		\end{equation*}
		The claim follows.
		\end{itemize}
\end{proof}

As in \cref{sec:lifumu:oracle:known}, we are now interested in deriving rates of convergence that result from the 
data-driven risk derived in \cref{prop:lifumu:neworacle:known}.
As before, we thus impose some additional regularity assumptions and derive a corresponding maximal inequality before illustrating 
resulting rates of convergence.

\begin{corollary}[Local Estimation - Data-Driven Maximal Inequality - KED]\label{lem:lifumu:known:maximialinequality:datadriven}
	\mbox{}\\
	Let $N\in\IN$. Grant \cref{ass:regularity:lifumu} as well as \cref{ass:lifumu:A2}.
	Let $\magenta{\calK^N}$ and $\khat$ be chosen as in \cref{prop:lifumu:neworacle:known}.
	Define the \textit{local data-driven maximal risk}
	\begin{align*}
		\forall k\in\magenta{\kappa}\IN:\; \ddrisklocal\br{k|\rms,\rmb,\rmL}&\coloneqq \norm{\rmL\indkco{k}\rms^\dagger}_{\IL^2}^2+\frac{1}{N}\rho_N(k|\rmb)\norm{\rmL\indk{k}\rmb}_{\IL^2}^2. 
	\end{align*}
	Then, for all $N$ sufficiently large, we have
	\begin{align*}
		\IE_f^N\brE{\abs{\linbar_{\khat}-\lin(f)}^2} 
		&\leq \const_1\Constfeps\ConstfY(\rmr_2^2\lor\rmr_\infty^2)\times\min_{k\in\magenta{\kappa}\IN}\ddrisklocal \br{k|\rms,\rmb,\rmL} \\
			&\hspace{1cm}+ \const_2\times \frac{1}{\magenta{\kappa^2}N}\br{\ConstfY\Constfeps+128^\rmp\sigmabar_{Y,2\rmp}^2}
	\end{align*}
	for some numerical constants $\const_1,\const_2\in\pRz$.
\end{corollary}

\begin{proof}[Proof of \cref{lem:lifumu:known:maximialinequality:datadriven}]
	The proof follows directly from an application of \cref{prop:lifumu:neworacle:known}
	while adapting the argumentation in the proof of \cref{co:lifumu:known:oracle} (similarly as in \cref{lem:lifumu:known:maximialinequality}).
	Thus, the details are omitted.
\end{proof}

In the following illustration we summarize the rates of convergence for the local data-driven maximal risk
as derived in \cref{lem:lifumu:known:maximialinequality:datadriven} similarly as we did in \cref{example:lifumu:known:oracle:rates}.
To that end we assume again \cref{lifumu:representant:decay:ass} and consider scenarios in which $f$ and $f^\varepsilon$ are \textit{ordinary smooth}
or \textit{super smooth} for parameters $b,s\in\pRzz$, respectively.

\begin{illustration}[Rates of Convergence - Local Data-Driven Maximal Risk - KED]\label{example:lifumu:known:datadriven:rates}
	\mbox{}\\
	Under the assumptions of \cref{lem:lifumu:known:maximialinequality:datadriven} grant \cref{lifumu:representant:decay:ass}.
	\Cref{tab:lifumu:dd:known} provides the order of the 
	\textit{local data-driven maximal risk}
	\begin{equation*}
		\ddrisklocal\coloneqq \min_{k\in\IN}\ddrisklocal \br{k|\rms,\rmb,\rmL}
	\end{equation*}
	as derived in \cref{lem:lifumu:known:maximialinequality:datadriven} in the asymptotic regime $N\to\infty$.
	We compare it to the order of $\maxrisklocal$ as discussed in \cref{example:lifumu:known:oracle:rates}.
	Throughout, we assume that $$ 2s-2p+(2(p+b)+1)_+-1\geq\frac{\rmp+2+(2(p+b)+1)_+}{\rmp-1}.$$
	The derivation of the presented rates of convergence can be found in \cref{sec:lifumu:appendix:datadriven:known}.
	\begin{table}[H]
		\centering
		\footnotesize
		\setlength{\tabcolsep}{6pt}        
		\renewcommand{\arraystretch}{1} 
		\begin{tabular}{@{}
			C{0.125cm}
			C{2.5cm}!{\vrule width 0.6pt}
			C{0.3cm}
			C{0.3cm}!{\vrule width 0.6pt}
			C{3.125cm}
			C{4.125cm}
		  @{}}
	
		\toprule
		&\ref{lifumu:eq:decay:L} & $f$ & $f_\varepsilon$ &
		$\maxrisklocal$ & $\ddrisklocal$\\
		\midrule
		\end{tabular}
		\renewcommand{\arraystretch}{3} 
		\begin{tabular}{@{}
			C{0.125cm}
			C{2.5cm}!{\vrule width 0.6pt}
			C{0.3cm}
			C{0.3cm}!{\vrule width 0.6pt}
			C{3.125cm}
			C{4.125cm}
		  @{}}
		\midrule
		(1)
		&$p+b > -\frac{1}{2}$
		& \textsc{os} & \textsc{os}
		& $\displaystyle N^{-\frac{2(s-p)-1}{2(s+b)}}$
		& $\displaystyle \br{\frac{\ln(N)}{N}}^{\frac{2(s-p)-1}{2(s+b)}}$
		\\
		
		(2)
		& $p+b=-1/2$
		&  & 
		& $\displaystyle \frac{\ln(N)}{N}$
		& $\displaystyle \frac{\ln(N)^2}{N}$
		\\
		
		(3)
		& $p+b<-\frac{1}{2}$
		&  & 
		& $\displaystyle \frac{1}{N}$
		& $\displaystyle \frac{\ln(N)}{N}$
		\\
		
		\midrule
		
		(4)
		&$\displaystyle s>p+\frac{1}{2}$
		& \textsc{os}
		& \textsc{es}
		& $\displaystyle \ln(N)^{-\frac{2(s-p)-1}{2b}}$
		& $\displaystyle \ln(N)^{-\frac{2(s-p)-1}{2b}}$
		\\
		
		\midrule
		
		(5)
		&$p+b>-\frac{1}{2}$
		& \textsc{es} & \textsc{os}
		& $\displaystyle \frac{\ln(N)^{\frac{2(p+b)+1}{2s}}}{N}$
		& $\displaystyle \frac{\br{\ln\ln(N)}\ln(N)^{\frac{2(p+b)+1}{2s}}}{N}$
		\\
		
		(6)
		&$p+b = - \frac{1}{2}$
		&  & 
		& $\displaystyle \frac{\ln\ln(N)}{N}$
		& $\displaystyle \frac{(\ln\ln(N))^2}{N}$
		\\
		
		(7)
		&$p+b<-\frac{1}{2}$
		&  & 
		& $\displaystyle \frac{1}{N}$
		& $\displaystyle \frac{\ln\ln(N)}{N}$
		\\
		
		\bottomrule
		\end{tabular}
		
		\captionof{table}{Order of the local data-driven maximal risk $\ddrisklocal$ in the asymptotic regime $N\to\infty$ when the error distribution is known.}
		\label{tab:lifumu:dd:known}
		\end{table}
	\end{illustration}

	\begin{remark}[Rates of Convergence - Local Data-Driven Maximal Risk - KED]
		\mbox{}
		\begin{itemize}
			\item Evidently, the asymptotic order of the local data-driven maximal risk presented in \cref{example:lifumu:known:datadriven:rates}
				corresponds to the (minimax-optimal) asymptotic order of the local maximal risk as presented in \cref{example:lifumu:known:oracle:rates}
				(up to a deterioration by an additional log-factor).

			\item Under the assumptions in \cref{lem:lifumu:known:maximialinequality:datadriven}, corresponding rates of convergence for estimating the probability density function $f$ (its derivatives, the survival function, the CDF)
				can be deduced from \cref{tab:lifumu:dd:known} similarly as in \cref{lifumu:maximalinequalities:known:oracle:examples}.
		\end{itemize}
	\end{remark}

\section{Linear Functional Estimation - Unknown Error Distribution}\label{sec:lifumu:unknown}
\subsection{Oracle-type Inequality}\label{sec:lifumu:unknown:oracle}

	We now turn to the linear functional estimation with an unknown error distribution (UED), i.e. in the statistical regime \assuRef{SRII}.
	Keeping \cref{not: MellinT} in mind, we replace $\MelXtilde^k$ in \cref{def:linearfunctionalestimator:known}
	with $\MelXhat^k$ and propose
	\begin{equation*}
		\forall k\in\pRz:\;\linhat_k \coloneqq \int_{[-k,k]}\rmL(-t)\MelXhat^k(t)\rmd\lebesgue(t) = \int_{[-k,k]}\rmL(t)\MelYhat(t)\MelUhatrez(t)\1_{\frakMl}(t)\rmd\lebesgue(t)
	\end{equation*}
	with random threshold set
	\begin{equation*}
		\frakMl\coloneqq\brK{t\in\IR:\;\abs{\MelUhat(t)}^2\geq M^{-1}}.
	\end{equation*}
	To shorten notation in the forthcoming analysis,
	we define
    \begin{equation*}
       \forall t\in\IR:\; \MelXcheck(t)\coloneqq \MelY(t)\MelUhatrez(t)\indMl(t),\quad \forall k\in\pRz:\;\MelXcheck^k(t)\coloneqq \MelXcheck(t)\indk{k}(t),
    \end{equation*}
    as well as 
	\begin{equation*}
		\lincheck \coloneqq \int_{\IR}\rmL(t)\MelXcheck(t)\rmd\lebesgue(t)\quad\text{and}\quad\lincheck_k \coloneqq \int_{\IR}\rmL(t)\MelXcheck^k(t)\rmd\lebesgue(t).
	\end{equation*}
	We impose some additional moment assumptions on the error variable $\varepsilon$.
	\begin{assumption}[Higher Moments - Unknown Error Distribution]\label{lifumu:ass:A3}
		\mbox{}\\
		In the statistical regime \assuRef{SRII} grant \cref{ass:lifumu:A1}. Assume further, that $\sigmabar_{\varepsilon,2\rmp}^2\coloneqq 1 + \IE\brE{\varepsilon^{2\rmp(c-1)}}\in\pRz$ for $\rmp\in\nset{2}$, i.e. $f^\varepsilon\in \IL_{\rmx^{2(c-1)}}^1\cap  \IL_{\rmx^{4(c-1)}}^1$.
	\end{assumption}

	As in \cref{sec:lifumu:oracle:known} we start with a first upper bound of the MSE.
	\begin{lemma}[MSE - Unknown Error Distribution]\label{lem:Lifumu:unknown:MSE}
		\mbox{}\\
		Grant \cref{lifumu:ass:A3}. We then have for all $N,M\in\IN$ and $k\in\pRz$
		\begin{align*}
			\IE_{f,\varepsilon}^{N,M}\brE{\abs{\linhat_k-\lin(f)}^2}& = \IE_{f,\varepsilon}^{N,M}\brE{\abs{\linhat_k-\lincheck_k}^2}+\IE_{\varepsilon}^M\brE{\abs{\lincheck_k -\lin(f)}^2}\\
				&\leq \frac{2}{N}\Constfeps\ConstfY\sigmabar_{\varepsilon,2}^2\norm{\rmL\indk{k}\MelUrez}_{\IL^2}^2 +  3\norm{\rmL\indkco{k}\MelX}_{\IL^1}^2\\
					&\hspace{1cm}+\const\times \br{\sigmabar_{\varepsilon,2}^2+\sigmabar_{\varepsilon,4}^2}\norm{\rmL\indk{k}\MelX\br{1\lor (M\abs{\MelU}^2)}^{-\frac{1}{2}}}_{\IL^1}^2
		\end{align*}
		with a numerical constant $\const\in\pRz$.\\
	\end{lemma}

	\begin{proof}[Proof of \cref{lem:Lifumu:unknown:MSE}]
		The proof is a direct consequence of the computations derived in \cref{lem:lifumu:basicIdentity} and \cref{lem:lifumu:unknown:errorestimates}.
		We start by considering the conditional expectation $\IE_{f|\varepsilon}^{N,M}$ and apply
		\cref{lifumu:unknown:biasvarianceI} and \cref{lifumu:unknown:biasvarianceIII}.
		Afterward, we apply the expectation w.r.t $\IE_{\varepsilon}^M$ and utilize \cref{lifumu:unknown:errorestimatesI}, \cref{lifumu:unknown:errorestimatesII}, and \cref{lifumu:unknown:errorestimatesIII}, respectively.
		The claim follows.
	\end{proof}

\begin{corollary}[Local Estimation - Oracle-type Inequality - UED]\label{co:lifumu:unknown:oraclerisk}
	\mbox{}\\
	 Let $N,M\in\IN$. Under \cref{lifumu:ass:A3}, there exists an optimal tuning parameter $\kdiamond = \kdiamond(N)\in\pRz$,
	 such that 
	\begin{align*}
			\IE_{f,\varepsilon}^{N,M}\brE{\abs{\linhat_{\kdiamond}-\lin(f)}^2}
				&\leq 3\Constfeps\ConstfY\times\sigmabar_{\varepsilon,2}^2\inf_{k\in\pRz}\brK{\norm{\rmL\indkco{k}\MelX}_{\IL^1}^2+\frac{1}{N}\norm{\rmL\indk{k}\MelUrez}_{\IL^2}^2 }\\
					&\hspace{1cm}+\const\times \br{\sigmabar_{\varepsilon,2}^2+\sigmabar_{\varepsilon,4}^2}\norm{\rmL\MelX\br{1\lor (M\abs{\MelU}^2)}^{-\frac{1}{2}}}_{\IL^1}^2
		\end{align*}
		with a numerical constant $\const\in\pRz$.
\end{corollary}

\begin{proof}[Proof of \cref{co:lifumu:unknown:oraclerisk}]
	The claim is a direct consequence of \cref{lem:Lifumu:unknown:MSE}
	by adapting the argumentation of the proof of \cref{co:lifumu:known:oracle} to the statistical regime \assuRef{SRII}.
	The details are thus omitted.
\end{proof}

\begin{remark}[Oracle Choice]\label{rem:lifumu:oracle:unknown}
	\mbox{}\\
	The optimal choice $\kdiamond$ as derived in \cref{co:lifumu:unknown:oraclerisk} depends on the unknown Mellin transforms $\MelX$ and $\MelUrez$, respectively,
	 and is hence computationally infeasible.
	Thus, $\kdiamond$ is also called \textit{oracle} (or \textit{oracle choice}) and the whole inequality in \cref{co:lifumu:unknown:oraclerisk} is called \textit{oracle-type inequality}.
	In \cref{sec:lifumu:unknown:datadriven} we propose a fully data driven choice of the tuning parameter $k$ that relies neither on the knowledge of $\MelX$ nor $\MelU$.
\end{remark}

We are following the steps in \cref{sec:lifumu:oracle:known} and are now investigating rates of convergence that result from the 
oracle-type inequality \cref{co:lifumu:unknown:oraclerisk} by specifying the regularity of $f$ and $f^\varepsilon$ respectively, 
which leads to the Corollary below.

\begin{corollary}[Local Estimation - Maximal Inequality - UED]\label{co:lifumu:unknwon:maximalrisk}
	\mbox{}\\
	Let $N,M\in\IN$. Recall the notation of \cref{lem:lifumu:known:maximialinequality}.
	Under \cref{ass:regularity:lifumu} and \cref{lifumu:ass:A3}, there exists an optimal tuning parameter $\kdiamond = \kdiamond(N)\in\pRz$,
	 such that 
	\begin{align*}
			\IE_{f,\varepsilon}^{N,M}\brE{\abs{\linhat_{\kdiamond}-\lin(f)}^2}
				&\leq 3\Constfeps\ConstfY\sigmabar_{\varepsilon,2}^2(\rmr_2^2\lor\rmr_\infty^2)\times\inf_{k\in\pRz}\maxrisklocal(k|\rms,\rmb,\rmL)\\
					&\hspace{1cm}+\const\rmr_2^2\rmr_\infty^2\br{\sigmabar_{\varepsilon,2}^2+\sigmabar_{\varepsilon,4}^2}\times\norm{\rmL\rms^\dagger\br{1\lor \br{M(\rmb^\dagger)^2}}^{-\frac{1}{2}}}_{\IL^2}^2
		\end{align*}
		with a numerical constant $\const\in\pRz$.
\end{corollary}

\begin{proof}[Proof of \cref{co:lifumu:unknwon:maximalrisk}]
	Under \cref{ass:regularity:lifumu}, we have 
	\begin{equation*}
		\norm{\rmL\MelX\br{1\lor (M\abs{\MelU}^2)}^{-\frac{1}{2}}}_{\IL^1}^2 \leq \rmr_2^2\rmr_\infty^2\norm{\rmL\rms^\dagger\br{1\lor \br{M(\rmb^\dagger)^2}}^{-\frac{1}{2}}}_{\IL^2}^2
	\end{equation*}
	by an application of the Cauchy-Schwarz inequality.
	The claim is a direct consequence of \cref{co:lifumu:unknown:oraclerisk} adapting the argumentation of \cref{co:lifumu:known:oracle} to the present setting.
	The details are hence omitted.
\end{proof}

In the following illustration we summarize rates of convergence for the local maximal risk as derived in \cref{co:lifumu:unknwon:maximalrisk} by considering 
scenarios in which $f$ and $f^\varepsilon$ are \textit{ordinary smooth} or \textit{super smooth} for parameters $b,s\in\pRzz$, where we refer again to \cref{def:ordinaryandexponentially}.

\begin{illustration}[Rates of Convergence - Local Maximal Risk - UED]\label{example:lifumu:unknown:oracle:rates}
	\mbox{}\\
	Grant \cref{lifumu:representant:decay:ass}.
	\Cref{tab:lifumu:oracle:unknwon} provides the order of the local maximal risk as derived in \cref{co:lifumu:unknwon:maximalrisk} 
	in the asymptotic regime $N,M\to\infty$.
	While $\maxrisklocal$ remains as derived in \cref{example:lifumu:known:oracle:rates},
	we discuss the additional term
	\begin{equation*}
		\norm{\rmL\rms^\dagger\br{1\lor \br{M(\rmb^\dagger)^2}}^{-\frac{1}{2}}}_{\IL^2}^2,
	\end{equation*}
	which appears due to estimating the unknown Mellin transform $\MelU$. 
	As before, throughout we assume $s>p+\frac{1}{2}$.
	\begin{table}[H]
		\centering
		\footnotesize
		\setlength{\tabcolsep}{6pt}        
		\renewcommand{\arraystretch}{1} 
		\begin{tabular}{@{}
			C{0.125cm}
			C{2.5cm}!{\vrule width 0.6pt}
			C{0.3cm}
			C{0.3cm}!{\vrule width 0.6pt}
			C{3.125cm}
			C{4cm}
		  @{}}
	
		\toprule
		&\ref{lifumu:eq:decay:L} & $f$ & $f_\varepsilon$ &
		$\maxrisklocal$ & $\norm{\rmL\rms^\dagger\br{1\lor \br{M(\rmb^\dagger)^2}}^{-\frac{1}{2}}}_{\IL^2}^2$\\
		\midrule
		\end{tabular}
		\renewcommand{\arraystretch}{3} 
		\begin{tabular}{@{}
			C{0.125cm}
			C{2.5cm}!{\vrule width 0.6pt}
			C{0.3cm}
			C{0.3cm}!{\vrule width 0.6pt}
			C{3.125cm}
			C{4cm}
		  @{}}
		\midrule
		(1)
		&$p+b > -\frac{1}{2}$
		& \textsc{os} & \textsc{os}
		& $\displaystyle N^{-\frac{2(s-p)-1}{2(s+b)}}$
		& $\displaystyle M^{-\br{1\land \frac{2(s-p)-1}{2b}}}$
		\\
		
		(2)
		& $p+b=-1/2$
		&  & 
		& $\displaystyle \frac{\ln(N)}{N}$
		&  $\displaystyle\frac{1}{M}$
		\\
		
		(3)
		& $p+b<-\frac{1}{2}$
		&  & 
		& $\displaystyle \frac{1}{N}$
		& $\displaystyle\frac{1}{M}$
		\\
		
		\midrule
		
		(4)
		&$\displaystyle s>p+\frac{1}{2}$
		& \textsc{os}
		& \textsc{es}
		& $\displaystyle \ln(N)^{-\frac{2(s-p)-1}{2b}}$
		& $\displaystyle \ln(M)^{-\frac{2(s-p)-1}{2b}}$
		\\
		
		\midrule
		
		(5)
		&$p+b>-\frac{1}{2}$
		& \textsc{es} & \textsc{os}
		& $\displaystyle \frac{\ln(N)^{\frac{2(p+b)+1}{2s}}}{N}$
		&  $\displaystyle\frac{1}{M}$
		\\
		
		(6)
		&$p+b = - \frac{1}{2}$
		&  & 
		& $\displaystyle \frac{\ln\ln(N)}{N}$
		&  $\displaystyle\frac{1}{M}$
		\\
		
		(7)
		&$p+b<-\frac{1}{2}$
		&  & 
		& $\displaystyle \frac{1}{N}$
		& $\displaystyle\frac{1}{M}$
		\\
		
		\bottomrule
		\end{tabular}
		
		\captionof{table}{Order of the local maximal risk in the asymptotic regime $N,M\to\infty$ when the error distribution is unknown.}
		\label{tab:lifumu:oracle:unknwon}
		\end{table}

		We briefly discuss the derivation of (1) and (4). The other cases follow simultaneously.
			\begin{itemize}
				\item (1):
					Let $p+b>-\frac{1}{2}$ and fix some $M\in\IN$.
					We first see  
					\begin{align*}
						\Vert\rmL\rmsos^\dagger(1\lor (M(\rmbos^\dagger)^2))^{-1/2}\Vert_{\IL^2}^2 & = \Vert\rmL\rmsos^\dagger(1\lor (M(\rmbos^\dagger)^2))^{-1/2}\indk{M^{\frac{1}{2b}}}\Vert_{\IL^2}^2\\
							&\hspace{1cm} + \Vert\rmL\rmsos^\dagger(1\lor (M(\rmbos^\dagger)^2))^{-1/2}\indkco{M^{\frac{1}{2b}}}\Vert_{\IL^2}^2.
					\end{align*}
					The first term is of order $M^{-1} \lor M^{-\frac{2(s-p)-1}{2b}}$ while the second term is of order $ M^{-\frac{2(s-p)-1}{2b}}$.
					Distinguishing between $\frac{2(s-p)-1}{2b} \leq 1$ and $\frac{2(s-p)-1}{2b}>1$ yields the claim.

				\item (4): Follows from (1) by replacing $\rmbos$ by $\rmbes$ and 
					$M^{\frac{1}{2b}}$ by $\ln(M)^{\frac{1}{2b}}$.
			\end{itemize}
	\end{illustration}
	\begin{remark}[Rates of Convergence - Local Maximal Risk - UED]
		\mbox{}
		\begin{itemize}
			\item 
			We note that in all scenarios derived in \cref{example:lifumu:unknown:oracle:rates},
			we need to pay a price in terms of an additional rate of convergence in $M$ to estimate the unknown Mellin transform $\MelU$.
			However, the additional term $\Vert\rmL\rmsos^\dagger(1\lor (M(\rmbos^\dagger)^2))^{-1/2}\Vert_{\IL^2}^2 $
			is of lower order than $\maxrisklocal$, and is hence negligible if $M\geq N$.  

			\item In the spirit of \cref{lifumu:maximalinequalities:known:oracle:examples}, we can deduce from \cref{tab:lifumu:oracle:unknwon}
			corresponding rates of convergence for estimating the probability density function $f$ (its derivatives, the survival function or its CDF).
		\end{itemize}
	\end{remark}

\subsection{Data Driven Model Selection}\label{sec:lifumu:unknown:datadriven}

In the spirit of  \cref{sec:lifumu:datadriven:known}, where we have investigated model selection methods 
for the tuning parameter $k$ in the statistical regime \assuRef{SRI}, 
we propose in this section a fully data-driven model selection method in the statistical regime \assuRef{SRII}.\\

More precisely, under \cref{ass:lifumu:A2} we recall \cref{not:lifumu:datadriven:stuff}.
The idea is to replace the unknown Mellin transform $\rmwU = \MelUrez$ 
by its empirical counterpart $\rmwhatU\coloneqq \MelUhatrez\1_{\frakMl}$.
Indeed, we define
\begin{equation*}
	\forall k\in\IR_{\geq \magenta{\kappa}}:\;\rmVhat_k^{\rmwhatU}\coloneqq \frac{1}{N}\sum_{i\in\nset{N}}\mathds{1}_{\big(-\infty,\nu_k^{\rmwhatU}\big]}(Y_i^{c-1})\left|\int_{[-k,k]}\rmL(-t)\MelUhatrez(t)\1_{\frakMl}(t)Y_i^{c-1+2\pi\iota t}\mathrm{d}\lebesgue(t)\right|^2
\end{equation*}
with 
\begin{equation*}
		\forall k\in\IR_{\geq \magenta{\kappa}}:\; \nu_k^{\rmwhatU}\coloneqq \br{\frac{1}{128}\br{\Delta_k^{\rmwhatU}}^\frac{1}{\rmp-1}}^\frac{1}{2}\in\pRz,
\end{equation*}
Conditioned on the additional data $\br{\varepsilon_i}_{i\in\nset{M}}$, $\rmVhat_k^{\rmwhatU}$ defines an unbiased estimator of $\rmV_k^{\rmwhatU}$,
i.e.  $$\forall k\in\IR_{\geq \magenta{\kappa}}:\;\IE_{f|\varepsilon}^{N,M}\brE{\rmVhat_k^{\rmwhatU}} = \rmV_k^{\rmwhatU}\quad \IP_{\varepsilon}^M - \text{a.s.}$$
We define an \textit{estimated penalty} for each $k\in\IR_{\geq \magenta{\kappa}}$ by 
\begin{equation}\label{eq:ligfumu:unknown:penhat}
	\penhat_k^{\rmwhatU}\coloneqq 8\br{\orange{2\rmVhat_k^{\rmwhatU}}+\orange{3}\times 2k\norm{\rmL\MelUhatrez\1_{\frakMl}\indk{k}}_{\IL^2}^2\frac{\br{\Delta_k^{\rmwhatU}}^{\frac{1}{\rmp-1}}\ln(\Delta_k^{\rmwhatU})}{N}} \frac{\ln(\Delta_k^{\rmwhatU})}{N} \in\pRzz,
\end{equation}
which is clearly fully known and thus computationally feasible.
We propose a \textit{fully data-driven} model selection method by 
\begin{equation}\label{eq:lifumu:khat:unknown:fullydatadriven}
	\khat:\in \argmin_{k\in\magenta{\calK}}\brK{\contrasthat_k^{\rmwhatU} +6\penhat_k^{\rmwhatU}}
\end{equation}
with contrast function
\begin{equation}\label{eq:lifumu:contrast:unknown:datadriven}
	\forall k\in\magenta{\calK}:\;\contrasthat_k^{\rmwhatU} \coloneqq\max_{j\in\magenta{\calK_{>k}}}\brK{\br{|\linhat_j-\linhat_k|^2-6\penhat_j^{\rmwhatU}-6\penhat_k^{\rmwhatU}}_+}.
\end{equation}
In the rest of this section we derive an upper bound for the local data-driven estimation risk
\begin{equation*}
	\IE_{f,\varepsilon}^{N,M}\brE{\abs{\linhat_{\khat}-\lin(f)}^2}.
\end{equation*}
As in \cref{sec:lifumu:datadriven:known} we start by providing a key argument.
\begin{lemma}[Local Estimation - Key Argument III - Unknown Error Distribution]\label{lifumu:unknown:dd:maindecomposition}
	\mbox{}\\
	Grant \cref{ass:lifumu:A2}.
	Select $\khat$ as in \cref{eq:lifumu:khat:unknown:fullydatadriven}.
	Defining 
	\begin{align*}
		\forall k\in\pRz:\;\rmZcheck_k^{\rmwhatU} \coloneqq  \linhat_k - \lincheck_k 
			& = \int_{[-k,k]}\rmL(t)\MelUhatrez(t)\1_{\frakMl}(t)\br{\MelYhat(t)-\MelY(t)}\rmd\lebesgue(t),
	\end{align*}
	and $\penbar_k^{\rmwhatU}$ as $\penbar^{\rmwU}_k$ in \cref{co:lifumu:known:Concentration}
	by replacing $\rmwU = \MelUrez$ with $\rmwhatU = \MelUhatrez\1_{\frakMl}$,
	we obtain
	\begin{align*}
		\abs{\linhat_{\khat}-\lin(f)}^2 & 
			\leq  4\times \min_{k\in\magenta{\calK}}\brK{6\norm{\rmL\indkco{k}\MelXcheck}_{\IL^1}^2 + 12\penhat_{k}^{\rmwhatU} + \abs{\linhat_{k}-\lincheck}^2}\\
				&\hspace{1cm}  + 48\max_{k\in\magenta{\calK}}\brK{\br{\abs{\rmZ_k^{\rmwhatU}}^2-2\penbar_k^{\rmwhatU}}_+} + 96\max_{k\in\magenta{\calK}}\brK{\br{\penbar_k^{\rmwhatU}-\penhat_k^{\rmwhatU}}_+}\\
				&\hspace{1cm} + 2\abs{\lincheck-\lin(f)}^2,
	\end{align*}
	where the inequality is understood to hold $\IP_{f,\varepsilon}^{N,M}$-almost surely.
\end{lemma}

\begin{proof}[Proof of \cref{lifumu:unknown:dd:maindecomposition}]
	In the following considerations, all inequalities and identities are understood to hold $\IP_{f,\varepsilon}^{N,M}$-almost surely without further notice.
	The general idea is to combine \cref{lem:lifumu:known:KeyArgument} and \cref{lem:lifumu:riskbound:dd:known}.
	Indeed, noting first
	\begin{align*}
		\abs{\linhat_{\khat}-\lin(f)}^2 \leq 2\abs{\linhat_{\khat}-\lincheck(f)}^2 + 2\abs{\lincheck-\lin(f)}^2,
	\end{align*}
	replacing $\lin(f)$ by $\lincheck(f)$ and setting $\pen_k \coloneqq \penhat_{k}^{\rmwhatU}$ in the proof of \cref{lem:lifumu:known:KeyArgument},  
	we obtain for any fixed $k\in\magenta{\calK}$
	\begin{align*}
		|\linhat_{\khat}-\lincheck(f)|^2 &\leq 12\norm{\rmL\MelXcheck\mathds{1}_{[-k,k]^C}}_{\IL^1}^2+24\penhat_{k}^{\rmwhatU} + 2|\linhat_k-\lincheck(f)|^2 \\
		&\hspace{1cm}+ 24\max_{k\in\magenta{\calK}}\brK{\br{\abs{\rmZ^{\rmwhatU}_k}^2-2\penhat_k^{\rmwhatU}}_+}.
	\end{align*}
	Evidently,
	\begin{equation*}
		\br{\abs{\rmZ^{\rmwhatU}_k}^2-2\penhat_k^{\rmwhatU}}_+ \leq \br{\abs{\rmZ^{\rmwhatU}_k}^2-2\penbar_k^{\rmwhatU}}_+ +2\br{\penbar_k^{\rmwhatU}-\penhat_k^{\rmwhatU}}_+,
	\end{equation*}
	such that 
	\begin{align*}
		|\linhat_{\khat}-\lincheck(f)|^2 &\leq 12\norm{\rmL\MelXcheck\mathds{1}_{[-k,k]^C}}_{\IL^1}^2+24\penhat_{k}^{\rmwhatU} + 2|\linhat_{k}-\lincheck(f)|^2 \\
		&\hspace{1cm}+ 24\max_{k\in\magenta{\calK}}\brK{\br{\abs{\rmZ^{\rmwhatU}_k}^2-2\penbar_k^{\rmwhatU}}_+} + 48\max_{k\in\magenta{\calK}}\brK{\br{\penbar_k^{\rmwhatU}-\penhat_k^{\rmwhatU}}_+}.
	\end{align*}
	As the last display holds for all $k\in\magenta{\calK}$, we can apply the minimum on both sides, which yields the claim.
\end{proof}

To upper bound the local data-driven risk 
\begin{equation*}
	\IE_{f,\varepsilon}^{N,M}\brE{|\linhat_{\khat}-\lin(f)|^2}
\end{equation*}
using \cref{lifumu:unknown:dd:maindecomposition},
it suffices to find appropriate upper bounds of
\begin{equation*}
	\IE_{f,\varepsilon}^{N,M}\brE{\penhat_{k}^{\rmwhatU}},\quad \IE_{f,\varepsilon}^{N,M}\brE{\max_{k\in\magenta{\calK}}\brK{\br{\abs{\rmZ^{\rmwhatU}_k}^2-2\penbar_k^{\rmwhatU}}_+}}
\end{equation*}
and 
\begin{equation*}
	\IE_{f,\varepsilon}^{N,M}\brE{\max_{k\in\magenta{\calK}}\brK{\br{\penbar_k^{\rmwhatU}-\penhat_k^{\rmwhatU}}_+}}.
\end{equation*}
Generally speaking, to derive such bounds
we condition on the additional data 
$(\varepsilon_j)_{j\in\nset{M}}$ first, apply the conditional expectation, and 
obtain similar results as derived in 
\cref{sec:lifumu:datadriven:known}.
The final result is then derived by applying the expectation associated to $(\varepsilon_j)_{j\in\nset{M}}$.

\begin{lemma}[Conditional Bound of the Estimated Penalty]\label{lem:lifumu:penhat:unknown:estimate}
	\mbox{}\\
	Grant \cref{ass:lifumu:A2}. Recall \cref{not:lifumu:datadriven:stuff}.
	Select $\khat$ as in \cref{eq:lifumu:khat:unknown:fullydatadriven}.
	Under \cref{lifumu:ass:A3}, we have for all $N,M\in\IN$ and $k\in\magenta{\calK}$.
		\begin{align*}
			\IE_{f|\varepsilon}^{N,M}\brE{48\penhat_{k}^{\rmwhatU}+4\abs{\linhat_{k}-\lincheck(f)}^2} &\leq \frac{1924}{N}\Constfeps\ConstfY\norm{\rmL\indk{k}\MelUhatrez\1_{\frakMl}}_{\IL^2}^2 \rho_N(k|\MelUhatrez\1_{\frakMl})\\
				&\hspace{1cm} + 4\times  \norm{\rmL\MelXcheck\indkco{k}}_{\IL^1}^2\quad\IP_\varepsilon^M-a.s.
		\end{align*}
\end{lemma}

\begin{proof}[Proof of \cref{lem:lifumu:penhat:unknown:estimate}]
	We start by noting that all (in-)equalities in this proof are to be understood $\IP_\varepsilon^M$-almost surely, 
	whenever we consider the conditional expectation.
	\begin{itemize}
		\item\textbf{Step I:} First note, an application of \cref{lifumu:unknown:app:estimatedpenality} yields
			\begin{equation*}
				\IE_{f|\varepsilon}^{N,M}\brE{\penhat_{k}^{\rmwhatU}} \leq \frac{40}{N}\times\Constfeps\ConstfY\norm{\rmL\indk{k}\MelUhatrez\1_{\frakMl}}_{\IL^2}^2 \rho_N(k|\MelUhatrez\1_{\frakMl}),
			\end{equation*}
			which holds for all $k\in\magenta{\calK}$.
		\item\textbf{Step II:} Secondly, for all fixed $k\in\magenta{\calK}$, we note that an application of \cref{lifumu:unknown:biasvarianceII} and \cref{lifumu:unknown:biasvarianceIII} in \cref{lem:lifumu:basicIdentity} leads to 
			\begin{equation*}
				\IE_{f|\varepsilon}^{N,M}\brE{\abs{\linhat_{k}-\lincheck(f)}^2} \leq \frac{1}{N}\Constfeps\ConstfY\norm{\rmL\MelUhatrez\1_{\frakMl}\indk{k}}_{\IL^2}^2 + \norm{\rmL\MelXcheck\indkco{k}}_{\IL^1}^2.
			\end{equation*}
			Recalling that $\rho_N(k|\MelUhatrez\1_{\frakMl})\geq 1$ and combining both bounds yields the claim.
		\end{itemize}
\end{proof}

The next assertion provides an upper bound for $\rho_N(k|\MelUhatrez\1_{\frakMl})$, which depends on the estimated Mellin transform $\MelUhatrez$.
We distinguish between the event when $\MelUhat$ is \textit{sufficiently close} to $\MelU$, and its complement.

\begin{lemma}[Event of Exception I]\label{lem:lifumu:unknwon:rhoestimate}
	\mbox{}\\
	Grant \cref{ass:lifumu:A2}. Define for $k\in\pRz$ the \textit{Event of Exception}
	\begin{equation}\label{eq:mho}
		\mho_{k}\coloneqq\left\{\sup_{t\in[-k,k]}\left|\MelUhat(t)\MelUrez(t)-1\right|\leq\frac{1}{3}\right\}\subseteq\left\{\sup_{t\in[-k,k]}\left|\MelUhat^\dagger(t)\MelU(t)\right|^2\leq\frac{9}{4}\right\},
		\end{equation}
		where its complement evidently satisfies
		\begin{equation*}
			\mho_{k}^C=\left\{\exists t\in[-k,k]:\,\left|\MelUhat(t)-\MelU(t)\right|>\frac{1}{3}|\MelU(t)|\right\}.
		\end{equation*}
		Keeping \cref{not:lifumu:datadriven:stuff} in mind, we have for all $N,M\in\IN$ and $k\in\pRz$
	\begin{align*}
		&\frac{1}{N}\norm{\rmL\indk{k}\MelUhatrez\1_{\frakMl}}_{\IL^2}^2 \rho_N(k|\MelUhatrez\1_{\frakMl})\\
			 & \hspace{2cm}\leq \const\times \frac{1}{N} \rho_N(k|\MelUrez)\norm{\MelUrez \rmL\indk{k}}_{\IL^2}^2\\
			&\hspace{3cm}\times\br{\1_{\mho_k}+ \br{1\lor\norm{\MelU}_{\IL^\infty}^6}\br{1\lor\ln(M)}^2M^{1+\frac{1}{\rmp-1}}\1_{\mho_k^C}}\quad \IP_\varepsilon^M-a.s.,
	\end{align*}
	where $\const\in\pRz$ is a numerical constant.
\end{lemma}

\begin{proof}[Proof of \cref{lem:lifumu:unknwon:rhoestimate}]
	We first note that $\mho_{k}$ and its complement $\mho_{k}^C$ are measurable events.
	Fix $k\in\pRz$. Clearly, on $\mho_k$ we have
		\begin{equation*}
			\abs{\MelUhatrez(t)(\omega)\1_{\frakMl}(t)} \leq \frac{3}{2}\abs{\MelUrez(t)}
		\end{equation*}
		for $\lebesgue$-almost all $t\in[-k,k]$. With \cref{eq:lifumu:rho:calculus2} it follows
		\begin{align*}
			\norm{\rmL\indk{k}\MelUhatrez\1_{\frakMl}}_{\IL^2}^2 \rho_N(k|\MelUhatrez\1_{\frakMl})
				&\leq \frac{9}{4}\norm{\rmL\indk{k}\MelUrez}_{\IL^2}^2 \times 2\times \br{\frac{3}{2}}^4\rho_N(k|\MelUrez).
		\end{align*}
		On the other hand, for $\lebesgue$-almost all $t\in[-k,k]$
		\begin{align*}
		\abs{\MelUhatrez(t)\1_{\frakMl}(t)} \leq M^\frac{1}{2}
			&\leq M^\frac{1}{2} \br{1\lor\norm{\MelU}_{\IL^\infty}}\abs{\MelUrez(t)}.
		\end{align*}
		Further, defining $A\coloneqq M^{\frac{1}{2}}\br{1\lor\norm{\MelU}_{\IL^\infty}}\in\IR_{\geq 1}$, elementary computations show
		\begin{align*}
			(1+2\ln(A))^2A^{\frac{2}{\rmp-1}} 
				& \leq 9(1\lor\ln(M))^2\times \br{1\lor\norm{\MelU}_{\IL^\infty}^{4}}\times M^{\frac{1}{\rmp-1}},
		\end{align*}
		which together with \cref{eq:lifumu:rho:calculus2} implies
		\begin{align*}
			\rho_N(k|\MelUhatrez\1_{\frakMl}) \leq  9(1\lor\ln(M))^2\times \br{1\lor\norm{\MelU}_{\IL^\infty}^{4}}\times M^{\frac{1}{\rmp-1}}\rho_N(k|\MelUrez).
		\end{align*}
		As a consequence,
		\begin{align*}
			&\norm{\rmL\indk{k}\MelUhatrez\1_{\frakMl}}_{\IL^2}^2 \rho_N(k|\MelUhatrez\1_{\frakMl})\\ 
				&\hspace{1cm} \leq 9\norm{\rmL\indk{k}\MelUrez}_{\IL^2}^2 \times  (1\lor\ln(M))^2\times \br{1\lor\norm{\MelU}_{\IL^\infty}^{6}}\times M^{1+\frac{1}{\rmp-1}}\rho_N(k|\MelUrez).
		\end{align*}
\end{proof}

The next assertion provides conditional upper bounds for the terms
\begin{equation*}
	\IE_{f|\varepsilon}^{N,M}\brE{\max_{k\in\magenta{\calK}}\brK{\br{\abs{\rmZ^{\rmwhatU}_k}^2-2\penbar_k^{\rmwhatU}}_+}}\quad\and\quad\IE_{f|\varepsilon}^{N,M}\brE{\max_{k\in\magenta{\calK}}\brK{\br{\penbar_k^{\rmwhatU}-\penhat_k^{\rmwhatU}}_+}},
\end{equation*}
respectively.

\begin{lemma}[Concentration of $\rmZ^{\rmwhatU}_k$ and $\penbar_k^{\rmwhatU}$]\label{lem:lifumu:unknown:concentration}
	\mbox{}\\
	Under \cref{ass:lifumu:A2} and \cref{lifumu:ass:A3}, let $\penbar_k^{\rmwhatU}$ be defined as $\penbar^{\rmwU}_k$ in \cref{co:lifumu:known:Concentration}
	with $\rmwU = \MelUrez$ replaced by $\rmwhatU = \MelUhatrez\1_{\frakMl}$.
	Then for all $M,N\in\IN$, we have
	\begin{equation*}
		\IE_{f|\varepsilon}^{N,M}\brE{\max_{k\in\magenta{\calK}}\brK{\br{\abs{\rmZ^{\rmwhatU}_k}^2-2\penbar_k^{\rmwhatU}}_+}} \leq \frac{1}{\magenta{\kappa^2}N}\br{40\Constfeps\ConstfY+\frac{128^\rmp}{32}\sigmabar_{Y,2\rmq}^2}\quad \IP_\varepsilon^M-a.s.,
	\end{equation*}
	as well as 
	\begin{equation*}
		\IE_{f|\varepsilon}^{N,M}\brE{\max_{k\in\magenta{\calK}}\brK{\br{\penbar_k^{\rmwhatU}-\penhat_k^{\rmwhatU}}_+}} \leq \Constfeps\ConstfY\frac{1}{\magenta{\kappa^2}N}\quad \IP_\varepsilon^M-a.s.
	\end{equation*}
\end{lemma}
\begin{proof}[Proof of \cref{lem:lifumu:unknown:concentration}]
	The proof follows directly from \cref{propo:lifumu:known:concentration:process:app} and \cref{lifumu:unknown:app:estimatedpenality} both by choosing $\rmw = \rmL\rmwhatU$.
\end{proof}
To obtain a final upper bound of
	\begin{equation}\label{eq:dd:explanation}
		\IE_{f,\varepsilon}^{N,M}\brE{\abs{\linhat_{\khat}-\lin(f)}^2},
	\end{equation}
	the presented proof strategy thus requires
		to find an appropriate upper bound for $\IP_{\varepsilon}^M(\mho_{k}^C)$.
		To that end, we make use of empirical process theory to upper control the 
		supremum of a normalized Mellin function process.

\begin{proposition}[Normalized Mellin Function Process]\label{lem:8maxi}
	\mbox{}\\
	Let $(Z_j)_{j\in\IN}$ be a family of i.i.d. $\pRz$-valued random variables
	and assume there exists a constant $\eta\in\Rz_{\geq1}$, such that 
	$$\eta\geq\max\brK{\br{\E\brE{Z_1^{2\beta}}}^{\frac{1}{2}},\E\brE{Z_1^{2\beta}\abs{\ln(Z_1)}^\gamma}}$$ for some
	$\beta\in\Rz$ and $\gamma\in\pRz$. For $m\in\IN$, define the normalized
	Mellin function process by
	\begin{equation}\label{lem:8maxi:e1}
		\forall t\in\IR:\;\rmW_m(t)\coloneqq\frac{1}{\sqrt{m}}\sum_{j\in\nset{m}}\br{Z_{j}^{\beta+2\pi\iota
			t}-\IE\left[Z_j^{\beta+2\pi\iota
			t}\right]}.
	\end{equation}
	Further, for $\rho\in\pRz$ define the weight function $\rmwbar\colon\IR\to(0,1]$ by 
	\begin{equation}\label{lem:8maxi:e2}
		\rmwbar(t)\coloneqq\left(\ln(\rme+\abs{t})\right)^{-\frac{1}{2}-\rho},\quad\forall t\in\Rz.
	\end{equation}
	Then there exists a constant $\Const(\eta,\rho)\in\Rz_{\geq1}$ only
	depending on $\eta$ and $\rho$, such that
	\begin{equation}\label{eq:EPT:Lifumu}
		\sup_{m\in\IN}\IE\left[\norm{\rmW_m}_{\IL_{\rmwbar}^\infty}\right]\leq \Const(\eta,\rho).
	\end{equation}
\end{proposition}

\begin{proof}[Proof of \cref{lem:8maxi}]
	The proof is shown in \cite[Proposition 4.4.]{Migueletal_2024}.
\end{proof}

\begin{proposition}[Event of Exception II]\label{lifumu:bound:event}
	\mbox{}\\
    Grant \cref{ass:lifumu:A2}. 
	In addition, assume
	$f^\varepsilon\in\IL_{\rmx^{7(c-1)}}^1\cap\IL_{\rmx^{2(c-1)}|\ln\rmx|^\gamma}^1$ for some $\gamma\in\pRz$
	and $\inf_{t\in[-\magenta{\kappa},\magenta{\kappa}]}\{|\MelU(t)|\}\in\pRz$.
	Let $\eta_\varepsilon\in\pRz$, such that
	 $$\eta_\varepsilon\geq\max\left\{\norm{f^{\varepsilon}}_{\IL_{\rmx^{7(c-1)}}^1}^{1/7},\norm{f^{\varepsilon}}_{\IL_{\rmx^{2(c-1)}|\ln\rmx|^\gamma}^1}\right\}.$$
	Let $\rho\in\pRz$ be
	arbitrary but fixed.  Consider the weight-function
	$\rmwbar:\IR\to(0,1]$ defined in \cref{lem:8maxi:e2},
	and let
	$\Const(\eta_\varepsilon,\rho)\in\Rz_{\geq1}$ be given by \cref{lem:8maxi}
	(with $\beta=c-1$). Given the universal numerical constant
	$\Const_{\mathrm{tal}}\in\IR_{>0}$ determined by Talagrand's
	inequality in \cite[Lemma C.1]{Migueletal_2024} we set 
	\begin{align}\nonumber
		\forall M\in\IN:\;\tau_M&\coloneqq\tau_M(\eta_\varepsilon,\rho)\coloneqq 2\eta_\varepsilon\Const_{\mathrm{tal}}^{-1/2}\ln(
		M)^{1/2}+2\Const(\eta_\varepsilon,\rho)
	\end{align}
	and 
	\begin{align}
		M_{\circ}&\coloneqq M_{\circ}(\MelU,\eta_\varepsilon,\rho)\nonumber\\
		&\coloneqq\min\left\{m\in\IN_{\geq3}:\,\inf_{t\in\magenta{[-\kappa,\kappa]}}\rmwbar(t)|\MelU(t)|\geq
		\min\brK{6\tau_m m^{-1/2},\Const_{\mathrm{tal}} \eta_\varepsilon^2m}\geq\ln(m)\right\}.\label{bound:event:mo}
	\end{align}
	For $M\in\nset{M_{\circ}}$ let  $K_M\in\pRz$ be arbitrary.
     For $M\in\IN_{\geq M_{\circ}}$ set
	\begin{equation}\label{bound:event:km}
		K_M\coloneqq\sup\left\{k\in\magenta{\kappa\IN}:\,\inf_{t\in[-k,k]}\brK{\rmwbar(t)|\MelU(t)|}\geq
		6\tau_M M^{-1/2}\right\}, 
	\end{equation}
	where the defining set is not empty for all $M\in\IN_{\geq
		M_{\circ}}$. For any $M\in\IN$ consider the event
	\begin{equation*}
		\mho_{K_M}^C\coloneqq\left\{\exists t\in[-K_M,K_M]:\,\left|\MelUhat(t)-\MelU(t)\right|>\frac{1}{3}|\MelU(t)|\right\}.
	\end{equation*}
	Then, there exists a universal numerical constant
	$\Const{event}\coloneqq3+11\Const_{\mathrm{tal}}^{-3}\in\IR_{\geq 1}$, such that
	\begin{equation*}
		\forall M\in\IN:\;\IP_{\varepsilon}^M(\mho_{K_M}^C)\leq (M_{\circ}^2\vee
                \Const{event}\eta_\varepsilon)M^{-2}.
	\end{equation*}
\end{proposition}

\begin{proof}[Proof of \cref{lifumu:bound:event}]
	\magenta{The proof is shown for $\kappa = 1$ in \cite[Proposition 4.5.]{Migueletal_2024}.
	For general $\kappa\in\pRz$, the proof follows similar lines and is hence omitted.}
\end{proof}

\begin{theorem}[Local Estimation - Data-Driven Risk - UED]\label{lifumu:unknown:datadriven}
	\mbox{}\\
	Fix $N,M\in\IN$. Grant \cref{ass:lifumu:A2} with $\rmq \in\IR_{\geq 3}$, \cref{lifumu:ass:A3}, and additionally the conditions imposed in \cref{lifumu:bound:event}.
	Define
	\begin{equation*}
		\magenta{\calK^N}\coloneqq\brK{k\in\magenta{\kappa}\IN:\;\br{1\lor\ln(k)}\norm{\rmL\indk{k}}_{\IL^2}^2 \leq N\norm{\rmL\indk{1}}_{\IL^2}^2}.
	\end{equation*}
	Let
	\begin{equation*}
		\khat:\in \argmin_{k\in\magenta{\calK^N}}\brK{\contrasthat_k^{\rmwhatU} +6\penhat_k^{\rmwhatU}}
	\end{equation*}
	with
	\begin{equation*}
		\contrasthat_k^{\rmwhatU} \coloneqq\max_{j\in\magenta{\calK^N_{>k}}}\brK{\br{|\linhat_j-\linhat_k|^2-6\penhat_j^{\rmwhatU}-6\penhat_k^{\rmwhatU}}_+}
	\end{equation*}
	and $\penhat_k^{\rmwhatU}$ as in \cref{eq:ligfumu:unknown:penhat}.
	Let $K_M\in\IN$ as in \cref{lifumu:bound:event}.
	Then for all $N\in\IN$ large enough such that
	\begin{equation*}
		\norm{\MelU}_{\IL^\infty}^2\min_{k\in\magenta{\kappa}\IN}\brK{\norm{\rmL\MelX\mathds{1}_{[-k,k]^C}}_{\IL^1}^2 + \frac{1}{N}\norm{\rmL\MelUrez\indk{k}}_{\IL^2}^2\rho_N(k|\rmwU)}\leq\norm{\rmL\indk{1}}_{\IL^2}^2,
	\end{equation*}
	we have
	\begin{align*}
		&\IE_{f,\varepsilon}^{N,M}\brE{\abs{\linhat_{\khat}-\lin(f)}^2} \leq \const_1\times\Constfeps\ConstfY\br{\Mo^2\lor \Const{event}\eta_\varepsilon}\br{1\lor\norm{\MelU}_{\IL^\infty}^6}\\
			&\hspace{4cm}\times \min_{k\in\magenta{\kappa}\IN}\brK{\norm{\rmL\MelX\indkco{k}}_{\IL^1}^2+\frac{1}{N} \rho_N(k|\MelUrez)\norm{\MelUrez \rmL\indk{k}}_{\IL^2}^2}\\
			&\hspace{4cm}+ \frac{\const_2}{\magenta{\kappa^2}N}\br{\Constfeps\ConstfY+128^\rmp\sigmabar_{Y,2\rmq}^2}\\
			&\hspace{4cm}+ \const_3\times \br{\sigmabar_{\varepsilon,2}^2+\sigmabar_{\varepsilon,4}^2}\norm{\rmL\MelX\br{1\lor (M\abs{\MelU}^2)}^{-\frac{1}{2}}}_{\IL^1}^2\\
			&\hspace{1cm} +  \const_4\times\Constfeps\ConstfY\br{\Mo^2\lor \Const{event}\eta_\varepsilon}\br{1\lor\norm{\MelU}_{\IL^\infty}^6}\norm{\rmL\MelX\indkco{K_M}}_{\IL^1}^2
	\end{align*} 
	for some numerical constants $\const_1,\const_2,\const_3,\const_4\in\pRz$.
\end{theorem}

\begin{proof}[Proof of \cref{lifumu:unknown:datadriven}]
	The proof is divided in several steps.
	We note that all inequalities are to be understood $\IP_\varepsilon^M$-almost surely, when we consider the conditional expectation given $(\varepsilon_j)_{j\in\nset M}$.
	\begin{itemize}
		\item \textbf{Step I:} We start by upper bounding the local data-driven risk conditioned on the 
		additional sample $(\varepsilon_j)_{j\in\nset{M}}$.
		Applying \cref{lifumu:unknown:dd:maindecomposition}, we obtain for any $\ko\in\magenta{\calK^N\cap(0,K_M]}$
		\begin{align*}
			&\IE_{f|\varepsilon}^{N,M}\brE{\abs{\linhat_{\khat}-\lin(f)}^2} 
				\leq  24\norm{\rmL\indkco{\ko}\MelXcheck}_{\IL^1}^2 + 48\IE_{f|\varepsilon}^{N,M}\brE{\penhat_{\ko}^{\rmwhatU}}\\ 
					&\hspace{2cm} + 4\IE_{f|\varepsilon}^{N,M}\brE{\abs{\linhat_{\ko}-\lincheck}^2}
					  + 48\IE_{f|\varepsilon}^{N,M}\brE{\max_{k\in\magenta{\calK^N}}\brK{\br{\abs{\rmZ_k^{\rmwhatU}}^2-\pen_k^{\rmwhatU}}_+}}\\
					&\hspace{2cm} + 96\IE_{f|\varepsilon}^{N,M}\brE{\max_{k\in\magenta{\calK^N}}\brK{\br{\penbar_k^{\rmwhatU}-\penhat_k^{\rmwhatU}}_+}}
					 + 2\abs{\lincheck-\lin(f)}^2\\
					&\hspace{1cm} \eqqcolon T_1+T_2+T_3+T_4+T_5+T_6.
		\end{align*}
		In the next step, we bound the conditional expectations from above, i.e. $T_2,T_3,T_4$ and $T_5$.
		
		\item \textbf{Step II:} Following \cref{lem:lifumu:penhat:unknown:estimate}  
				we obtain 
				\begin{align*}
					T_2+T_3  
					 &\leq 1924\times \frac{1}{N}\Constfeps\ConstfY\norm{\rmL\indk{\ko}\MelUhatrez\1_{\frakMl}}_{\IL^2}^2 \rho_N(\ko|\MelUhatrez\1_{\frakMl})\\
						&\hspace{1cm} + 4\times  \norm{\rmL\MelXcheck\indkco{\ko}}_{\IL^1}^2.
				\end{align*}
				Utilizing \cref{lem:lifumu:unknwon:rhoestimate}, the first summand can be further upper bounded, and we conclude 
				\begin{align*}
					T_2+T_3  
					 &\leq \const_1\times \Constfeps\ConstfY\frac{1}{N} \rho_N(\ko|\MelUrez)\norm{\MelUrez \rmL\indk{\ko}}_{\IL^2}^2\\
					 &\hspace{3cm}\times\br{\1_{\mho_{\ko}}+ \br{1\lor\norm{\MelU}_{\IL^\infty}^6}\br{1\lor\ln(M)}^2M^{1+\frac{1}{\rmp-1}}\1_{\mho_{\ko}^C}}\\
				 	&\hspace{1cm} + 4\times  \norm{\rmL\MelXcheck\indkco{\ko}}_{\IL^1}^2.
				\end{align*}
				Further, an application of \cref{lem:lifumu:unknown:concentration} yields
				\begin{equation*}
					T_4+T_5 \leq \frac{\const_2}{\magenta{\kappa^2}N}\br{\Constfeps\ConstfY+128^\rmp\sigmabar_{Y,2\rmq}^2}
				\end{equation*}
				with a finite constant $\const_2\in\pRz$.
				Note, the upper bounds of $T_4$ and $T_5$ are deterministic and thus do not need to be further considered.
				We have found upper bounds for all conditional expectations, which now only depend on the additional
				sample $(\varepsilon_j)_{j\in\nset{M}}$.
				In the next step we consider the expectations
					\begin{equation*}
						\IE_\varepsilon^M[T_1],\quad \IE_\varepsilon^M[T_2+T_3]\quad \and \quad \IE_\varepsilon^M[T_6].
					\end{equation*}

			\item \textbf{Step III:} We can deduce from the proof of \cref{lifumu:unknown:errorestimatesIII} that
				\begin{equation*}
					\IE_\varepsilon^M[T_1] = 24 \IE_{\varepsilon}^M\brE{\norm{\rmL\MelXcheck\indkco{\ko}}_{\IL^1}^2} \leq 48\sigmabar_{\varepsilon,2}^2\norm{\rmL\MelX\indkco{\ko}}_{\IL^1}^2.
				\end{equation*}
				Next, we have
				\begin{align*}
					\IE_\varepsilon^M[T_2+T_3]  
					 &\leq \const_1\times \Constfeps\ConstfY\frac{1}{N} \rho_N(\ko|\MelUrez)\norm{\MelUrez \rmL\indk{\ko}}_{\IL^2}^2\\
					 &\hspace{1cm}\times\IE_\varepsilon^M\brE{\br{\1_{\mho_{\ko}}+ \br{1\lor\norm{\MelU}_{\IL^\infty}^6}\br{1\lor\ln(M)}^2M^{1+\frac{1}{\rmp-1}}\1_{\mho_{\ko}^C}}}\\
				 	&\hspace{1cm} + 4\times  \IE_\varepsilon^M\brE{\norm{\rmL\MelXcheck\indkco{\ko}}_{\IL^1}^2}.
				\end{align*}
				We note that the second summand can be upper bounded as $\IE_\varepsilon^M[T_1]$.
				For the first term, we observe
				\begin{align*}
					&\IE_\varepsilon^M\brE{\br{\1_{\mho_{\ko}}+ \br{1\lor\norm{\MelU}_{\IL^\infty}^6}\br{1\lor\ln(M)}^2M^{1+\frac{1}{\rmp-1}}\1_{\mho_{\ko}^C}}}\\
						&\hspace{2cm} \leq \br{\IP_{\varepsilon}^M(\mho_{\ko})+ \br{1\lor\norm{\MelU}_{\IL^\infty}^6}\br{1\lor\ln(M)}^2M^{1+\frac{1}{\rmp-1}}\IP_{\varepsilon}^M(\mho_{\ko}^C)}.
				\end{align*}
				As \magenta{$\ko\in\calK^N \cap (0,K_M]$} with $K_M$ as defined in \cref{lifumu:bound:event}, we have 
				\begin{equation*}
					\IP_\varepsilon^M\br{\mho_{\ko}^C} \leq\IP_\varepsilon^M\br{\mho_{K_M}^C} \leq \br{\Mo^2\lor \Const{event}\eta_\varepsilon}\times M^{-2}
				\end{equation*}
				with $\Mo\geq 3$.
				Further, as $\rmp\geq 3$ 
				\begin{align*}
					\br{1\lor\ln(M)}^2M^{1+\frac{1}{\rmp-1}}M^{-2} & =\br{1\lor\ln(M)}^2 M^{-1 + \frac{1}{\rmp-1}}
						 \leq \br{\frac{4}{\rme}}^2\leq 3,
				\end{align*}
				where we have used $\frac{2(\rmp-1)}{(\rmp-2)}\in[2,4]$ and $\frac{\ln(x)^A}{x}\leq\br{\frac{A}{\rme}}^{A}$ for all $x\geq 1$.
				All together
				\begin{align*}
					&\br{\IP_\varepsilon^M\br{\mho_{\ko}}+ \br{1\lor\norm{\MelU}_{\IL^\infty}^6}\br{1\lor\ln(M)}^2M^{1+\frac{1}{\rmp-1}}\IP_\varepsilon^M\br{\mho_{\ko}^C}}\\
					&\hspace{1cm}\leq \br{1+ 3\br{\Mo^2\lor \Const{event}\eta_\varepsilon}\br{1\lor\norm{\MelU}_{\IL^\infty}^6}}\\
					&\hspace{1cm}\leq 6\br{\Mo^2\lor \Const{event}\eta_\varepsilon}\br{1\lor\norm{\MelU}_{\IL^\infty}^6}.
				\end{align*}
				We conclude
				\begin{align*}
					\IE_\varepsilon^M[T_2+T_3]  
					 &\leq \const_1\times \Constfeps\ConstfY\frac{1}{N} \rho_N(\ko|\MelUrez)\norm{\MelUrez \rmL\indk{\ko}}_{\IL^2}^2\\
					 &\hspace{1cm}\times6\br{\Mo^2\lor \Const{event}\eta_\varepsilon}\br{1\lor\norm{\MelU}_{\IL^\infty}^6}\\
				 	&\hspace{1cm} + 8\sigmabar_{\varepsilon,2}^2\norm{\rmL\MelX\indkco{\ko}}_{\IL^1}^2.
				\end{align*}
				Lastly, \cref{lifumu:unknown:errorestimatesIV} yields
				\begin{equation*}
					\IE_{\varepsilon}^M\brE{T_6} =2\IE_{\varepsilon}^M\brE{\abs{\lincheck-\lin(f)}^2} \leq \const_3\times \br{\sigmabar_{\varepsilon,2}^2+\sigmabar_{\varepsilon,4}^2}\norm{\rmL\MelX\br{1\lor (M\abs{\MelU}^2)}^{-\frac{1}{2}}}_{\IL^1}^2 
				\end{equation*}
				with a finite numerical constant $\const_3\in\pRz$.
			\item \textbf{Step IV:} Combining all bounds derived so far, there exists a numerical constant $\const_4\in\pRz$, such that for all $\ko\in\magenta{\calK^N\cap(0,K_M]}$
			\begin{align*}
				&\IE_{f,\varepsilon}^{N,M}\brE{\abs{\linhat_{\khat}-\lin(f)}^2} \leq \const_5\times\Constfeps\ConstfY\br{\Mo^2\lor \Const{event}\eta_\varepsilon}\br{1\lor\norm{\MelU}_{\IL^\infty}^6}\\
					&\hspace{4cm}\times \br{\norm{\rmL\MelX\indkco{\ko}}_{\IL^1}^2+\frac{1}{N} \rho_N(\ko|\MelUrez)\norm{\MelUrez \rmL\indk{\ko}}_{\IL^2}^2}\\
					&\hspace{4cm}+ \frac{\const_2}{N\magenta{\kappa^2}}\br{\Constfeps\ConstfY+128^\rmp\sigmabar_{Y,2\rmq}^2}\\
					&\hspace{4cm}+ \const_3\times \br{\sigmabar_{\varepsilon,2}^2+\sigmabar_{\varepsilon,4}^2}\norm{\rmL\MelX\br{1\lor (M\abs{\MelU}^2)}^{-\frac{1}{2}}}_{\IL^1}^2.
			\end{align*}
			
			\item \textbf{Step V:} Observe that the upper bound derived in \textbf{Step IV} holds for all $\ko\in\magenta{\calK^N\cap(0,K_M]}$. 
				Thus, by replacing $\ko$ with $\ko\land K_M$, the preceding inequality is also true for all $\ko\in\magenta{\calK^N}$.
				As $\magenta{\kappa\IN}\ni k \mapsto\rho_N(k|\MelUrez)\norm{\MelUrez \rmL\indk{k}}_{\IL^2}^2\in\pRzz$ is a non-decreasing function, 
				we can extend the last display for all $\ko\in\magenta{\calK^N}$ by
				\begin{align*}
					&\IE_{f,\varepsilon}^{N,M}\brE{\abs{\linhat_{\khat}-\lin(f)}^2} \leq \const_5\times\Constfeps\ConstfY\br{\Mo^2\lor \Const{event}\eta_\varepsilon}\br{1\lor\norm{\MelU}_{\IL^\infty}^6}\\
						&\hspace{4cm}\times \br{\norm{\rmL\MelX\indkco{\ko}}_{\IL^1}^2+\frac{1}{N} \rho_N(\ko|\MelUrez)\norm{\MelUrez \rmL\indk{\ko}}_{\IL^2}^2}\\
						&\hspace{4cm}+ \frac{\const_2}{N\magenta{\kappa^2}}\br{\Constfeps\ConstfY+128^\rmp\sigmabar_{Y,2\rmq}^2}\\
						&\hspace{4cm}+ \const_3\times \br{\sigmabar_{\varepsilon,2}^2+\sigmabar_{\varepsilon,4}^2}\norm{\rmL\MelX\br{1\lor (M\abs{\MelU}^2)}^{-\frac{1}{2}}}_{\IL^1}^2\\
						&\hspace{1cm} +  \const_5\times\Constfeps\ConstfY\br{\Mo^2\lor \Const{event}\eta_\varepsilon}\br{1\lor\norm{\MelU}_{\IL^\infty}^6}\norm{\rmL\MelX\indkco{K_M}}_{\IL^1}^2.
				\end{align*} 
				If we set
				\begin{equation*}
					\ko:\in\argmin_{k\in\magenta{\kappa}\IN}\brK{\norm{\rmL\MelX\mathds{1}_{[-k,k]^C}}_{\IL^1}^2 + \frac{1}{N}\norm{\rmL\MelUrez\indk{k}}_{\IL^2}^2\rho_N(k|\rmwU)},
				\end{equation*}
				as in the proof of \cref{prop:lifumu:neworacle:known}, we obtain similarly $\ko\in\magenta{\calK^N}$. The claim follows.
				\end{itemize}
\end{proof}
	Similar as in \cref{sec:lifumu:datadriven:known}, we are now deriving rates of convergence that
result from the local data-driven risk derived in \cref{lifumu:unknown:datadriven}. As before, we thus impose
some additional regularity assumptions and derive a corresponding maximal inequality
before illustrating resulting rates of convergence.

	\begin{corollary}[Local Estimation - Data-Driven Maximal Inequality - UED]\label{lem:lifumuu:nknown:maximialinequality:datadriven}
		\mbox{}\\
		Let $N,M\in\IN$. Grant \cref{ass:regularity:lifumu}.
		Let $\khat$, $\magenta{\calK^N}$ and $K_M$ be chosen as in and below the conditions of \cref{lifumu:unknown:datadriven}.
		Then, for all $N$ sufficiently large, we have
		\begin{align*}
			&\IE_{f,\varepsilon}^{N,M}\brE{\abs{\linhat_{\khat}-\lin(f)}^2} \leq \const_1\times\Constfeps\ConstfY\br{\Mo^2\lor \Const{event}\eta_\varepsilon}\br{1\lor\norm{\MelU}_{\IL^\infty}^6}\\
				&\hspace{5cm}\times \br{\rmr_2^2\lor\rmr_\infty^2}\times \min_{k\in\IN}\ddrisklocal(k|\rms,\rmb,\rmL)\\
				&\hspace{4cm}+ \frac{\const_2}{N\magenta{\kappa^2}}\br{\Constfeps\ConstfY+128^\rmp\sigmabar_{Y,2\rmq}^2}\\
				&\hspace{4cm}+ \const_3 \br{\sigmabar_{\varepsilon,2}^2+\sigmabar_{\varepsilon,4}^2}\rmr_2^2\rmr_\infty^2\times\norm{\rmL\rms^\dagger\br{1\lor (M(\rmb^\dagger)^2)}^{-\frac{1}{2}}}_{\IL^2}^2\\
				&\hspace{1cm} + \const_4\Constfeps\ConstfY\br{\Mo^2\lor \Const{event}\eta_\varepsilon}\br{1\lor\norm{\MelU}_{\IL^\infty}^6}\rmr_2^2\times\norm{\rmL\rms^\dagger\indkco{K_M}}_{\IL^2}^2
		\end{align*} 
		for some numerical constants $\const_1,\const_2,\const_3,\const_4\in\pRz$.
\end{corollary}

\begin{proof}[Proof of \cref{lem:lifumuu:nknown:maximialinequality:datadriven}]
	The claim is a direct consequence of \cref{lifumu:unknown:datadriven} and thus the details of the proof are omitted.
\end{proof}

In the following illustration we summarize rates of convergence for the local data-driven maximal risk as derived in 
\cref{lem:lifumuu:nknown:maximialinequality:datadriven}.
To that end, we consider scenarios in which $f$ and $f^\varepsilon$ are \textit{ordinary smooth} or \textit{super smooth} for parameters $b,s\in\pRzz$, where we refer again to \cref{def:ordinaryandexponentially}.

\begin{illustration}[Rates of Convergence - Data-Driven Maximal Inequality - UED]\label{example:lifumu:unknown:datadriven:rates}
	\mbox{}\\
	Grant the assumptions of \cref{lem:lifumuu:nknown:maximialinequality:datadriven} and additionally \cref{lifumu:representant:decay:ass}.
	\Cref{tab:lifumu:datadriven:unknown} provides the order of 	
	\begin{equation*}
		\norm{\rmL\rms^\dagger\br{1\lor \br{M(\rmb^\dagger)^2}}^{-\frac{1}{2}}}_{\IL^2}^2\quad\and\quad\norm{\rmL\rms^\dagger\indkco{K_M}}_{\IL^2}^2
	\end{equation*}
	in the asymptotic regime $M\to\infty$. 
	Keep in mind that $\ddrisklocal$ remains the same as derived in \cref{example:lifumu:known:datadriven:rates}.
	\begin{table}[H]
		\centering
		\footnotesize
		\setlength{\tabcolsep}{6pt}        
		\renewcommand{\arraystretch}{1} 
		\begin{tabular}{@{}
			C{0.125cm}
			C{2.5cm}!{\vrule width 0.6pt}
			C{0.3cm}
			C{0.3cm}!{\vrule width 0.6pt}
			C{3.125cm}
			C{4.125cm}
		  @{}}
	
		\toprule
		& \ref{lifumu:eq:decay:L} & $f$ & $f_\varepsilon$ &
		$\norm{\rmL\rms^\dagger\br{1\lor (M\abs{\rmb^\dagger}^2)}^{-\frac{1}{2}}}_{\IL^2}^2$ & $\norm{\rmL\rms^\dagger\indkco{K_M}}_{\IL^2}^2$\\
		\midrule
		\end{tabular}
		\renewcommand{\arraystretch}{3} 
		\begin{tabular}{@{}
			C{0.125cm}
			C{2.5cm}!{\vrule width 0.6pt}
			C{0.3cm}
			C{0.3cm}!{\vrule width 0.6pt}
			C{3.125cm}
			C{4.125cm}
		  @{}}
		\midrule
		(1)
		&$s > p+\frac{1}{2}$
		& \textsc{os} & \textsc{os}
		& $\displaystyle M^{-(1\land \frac{2(s-p)-1}{2b})}$
		& $\displaystyle \br{\frac{\ln(M)^{2(1+\rho)}}{M}}^{\frac{2(s-p)-1}{2b}}$
		\\
		
		(2)
		& $s > p+\frac{1}{2}$
		& \textsc{os} & \textsc{es}
		& $\displaystyle \ln(M)^{-\frac{2(s-p)-1}{2b}}$
		& $\displaystyle \ln(M)^{-\frac{2(s-p)-1}{2b}}$
		\\
		
		(3)
		& $s > p+\frac{1}{2}$
		& \textsc{es} & \textsc{os}
		& $\displaystyle \frac{1}{M}$
		& $\displaystyle \exp\br{-\ln(M)^{-(1+\rho)\frac{2s}{b}}M^\frac{s}{b}}$
		\\
		
		\bottomrule
		\end{tabular}
		
		\captionof{table}{Order of the additional terms in the local data-driven maximal risk in the asymptotic regime $M\to\infty$ when the error distribution is unknown.}
		\label{tab:lifumu:datadriven:unknown}
		\end{table}
		As we have seen in \cref{example:lifumu:known:oracle:rates}
		that 
		\begin{equation*}
			\norm{\rmL\rms^\dagger\indkco{K_M}}_{\IL^2}^2 \asymp K_M^{-2(s-p)+1}
		\end{equation*}
		in case of \textsc{os}, and 
		\begin{equation*}
			\norm{\rmL\rms^\dagger\indkco{K_M}}_{\IL^2}^2 \asymp \exp\br{-K_M^{-2s}}
		\end{equation*}
		in case of \textsc{es},
		it suffices to determine the order of $K_M$ from \cref{bound:event:km} in \cref{lifumu:bound:event}.
		We note that there are finite constants $\const_1=\const_1(\eta_\varepsilon,\rho)\in\pRz$ and $\const_2=\const_2(\rho)\in\pRz$, such that $\tau_M\leq \const_1 (\ln M)^{\frac{1}{2}}$ for all $M\in\IN_{\geq M_\circ}$ and 
				\begin{equation*}
					\inf_{t\in[-k,k]}|\rmwbar(t)|(\ln k)^{\frac{1}{2}+\rho} \geq \const_2,
				\end{equation*}
				for all $k\in\IR_{\geq3}$. Assume that $f^\varepsilon$ is \textit{ordinary smooth}. 
				Then there exists a finite constant $\const_3\in\pRz$, such that
				\begin{equation*}
					\inf_{t\in[-k,k]}|\MelU(t)|k^{b} \geq \const_3,
				\end{equation*}
				for all $k\in\IR_{\geq1}$.
				Consequently, there exists a constant $\const_4 \in\pRz$ small enough, such that for $k=\left\lfloor\const_4 M^{\frac{1}{2b}}(\ln M)^{\frac{-1-\rho}{b}}\right\rfloor$
				\begin{equation*}
					\inf_{t\in[-k,k]}|\rmwbar(t)||\MelU(t)|\geq  \const_2\const_3 k^{-b}(\ln k)^{-\frac{1}{2}-\rho}
					\geq 6 M^{-\frac{1}{2}} \const_1 (\ln M)^{\frac{1}{2}} \geq 
					6\tau_M M^{-\frac{1}{2}}.
				\end{equation*}
				Therefore, exploiting \cref{bound:event:km}, i.e. $K_M\geq k$ we obtain 
				\begin{equation*}
					\norm{\rmL\rmsos^\dagger\indkco{K_M}}_{\IL^2}^2 \asymp \left(\frac{(\ln M)^{2+2\rho}}{M}\right)^{\frac{2(s-p)-1}{2b}}.
				\end{equation*}
				This shows (1). (3) follows analogously.
				Now assume that $\MelU$ is \textit{super smooth}.
				Then
				\begin{equation*}
					\inf_{t\in[-k,k]}|\MelU(t)|\exp(k^{2b})\gtrsim 1,
				\end{equation*}
				for all $k\in\IR_{\geq1}$. Consequently, there exists a constant $\const_4 \in\pRz$ small enough, such that for $ k = \left\lfloor \left(\ln\frac{\const_4\cdot M}{(\ln M)^{2+2\rho}}\right)^{\frac{1}{2b}}\right\rfloor$
				\begin{equation*}
					\inf_{t\in[-k,k]}|\rmwbar(t)||\MelU(t)|\geq  \const_2 \exp(-k^{2b})(\ln k)^{-\frac{1}{2}-\rho}
					\geq 6 M^{-\frac{1}{2}} \const_1 (\ln M)^{\frac{1}{2}} \geq 
					6\tau_M M^{-\frac{1}{2}}.
				\end{equation*}
				Therefore, exploiting \cref{bound:event:km}, i.e. $K_M\geq k$, we obtain 
				\begin{equation*}
					\norm{\rmL\rmsos^\dagger\indkco{K_M}}_{\IL^2}^2 \asymp (\ln M)^{-\frac{2(s-p)-1}{2b}}.
				\end{equation*}
				This shows (2).
	\end{illustration}

\section{Simulation experiments}
\label{sec:simulation-experiments}

In this section, we complement the theoretical results of this work by a
simulation study. To that end, we restrict our attention to the estimation of the
probability density function
\(
f_X
\)
and the survival function
\(
S_X
\).
The experiments cover two qualitatively different inverse problems: an
ordinary-smooth setting and an exponentially smooth setting. In each case,
the Mellin transforms of the signal and the multiplicative error belong to the
same broad smoothness regime. 

\subsection{Simulation design}
\label{sec:simulation-design}

\subsubsection{Ordinary-smooth and exponentially smooth models}

We study the following two data-generating mechanisms.
\begin{itemize}
    \item \textbf{Ordinary-smooth setting (OS).}
    The signal follows a beta distribution and the error follows a Pareto
    distribution,
    \begin{equation*}
        X\sim\operatorname{Beta}(2,3),
        \qquad
        \varepsilon\sim\operatorname{Pareto}(3.1,2).
    \end{equation*}

    \item \textbf{Exponentially smooth setting (ES).}
    The signal follows a gamma distribution and the error is lognormal,
    \begin{equation*}
        X\sim\operatorname{Gamma}\br{3,\frac{7}{20}},
        \qquad
        \varepsilon\sim\operatorname{Lognormal}\br{\frac{1}{2},\frac{1}{16}}.
    \end{equation*}
    In this model, the
    Mellin transforms exhibit exponential or faster-than-polynomial decay.
    The inverse problem is consequently substantially more ill-posed than in
    the ordinary-smooth setting.
\end{itemize}

If we write shorthand
\(z=c-1+2\pi\iota t\), the Mellin transforms involved in the two models are
\begin{align*}
    \mathcal{M}_c[f_{\operatorname{Beta}(2,3)}](t)
        &=\frac{B(2+z,3)}{B(2,3)},
        &
    \mathcal{M}_c[f_{\operatorname{Pareto}(3.1,2)}](t)
        &=\frac{3.1\times 2^z}{3.1-z},\\
    \mathcal{M}_c[f_{\operatorname{Gamma}(3,7/20)}](t)
        &=(7/20)^{-z}\frac{\Gamma(3+z)}{\Gamma(3)},
        &
    \mathcal{M}_c[f_{\operatorname{LN}(1/2,1/16)}](t)
        &=\exp\!\left(\frac{z}{2}+\frac{z^2}{32}\right),
\end{align*}
where we refer to \cite{BrennerMiguel_2023} for computations.
Throughout, we fix \(c=1.5\) for both target functions $f_X$ and $\SX$, respectively.
In particular, this value satisfies the condition \(c\in\IR_{>1}\) required
for survival-function inversion, \cref{example:Illustration1}.\\

\begin{figure}[H]
    \centering
    \begin{subfigure}[t]{0.48\textwidth}
        \centering
        \includegraphics[width=\linewidth]{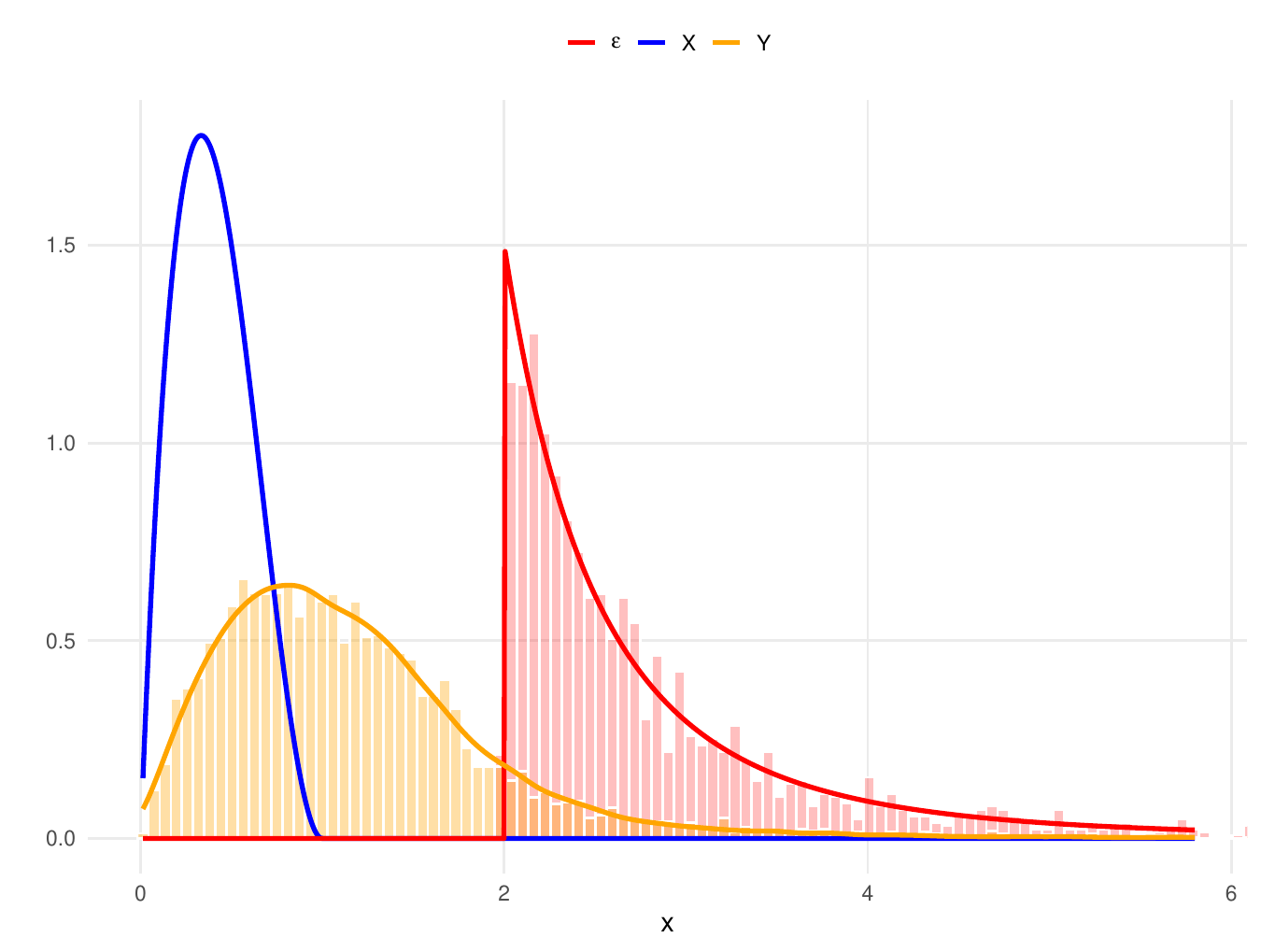}
        \caption{Ordinary-smooth model: $X\sim f_{\operatorname{Beta}(2,3)}$ and $\varepsilon\sim f_{\operatorname{Pareto}(3.1,2)}$ .}
        \label{fig:simulation-model-os}
    \end{subfigure}
    \hfill
    \begin{subfigure}[t]{0.48\textwidth}
        \centering
        \includegraphics[width=\linewidth]{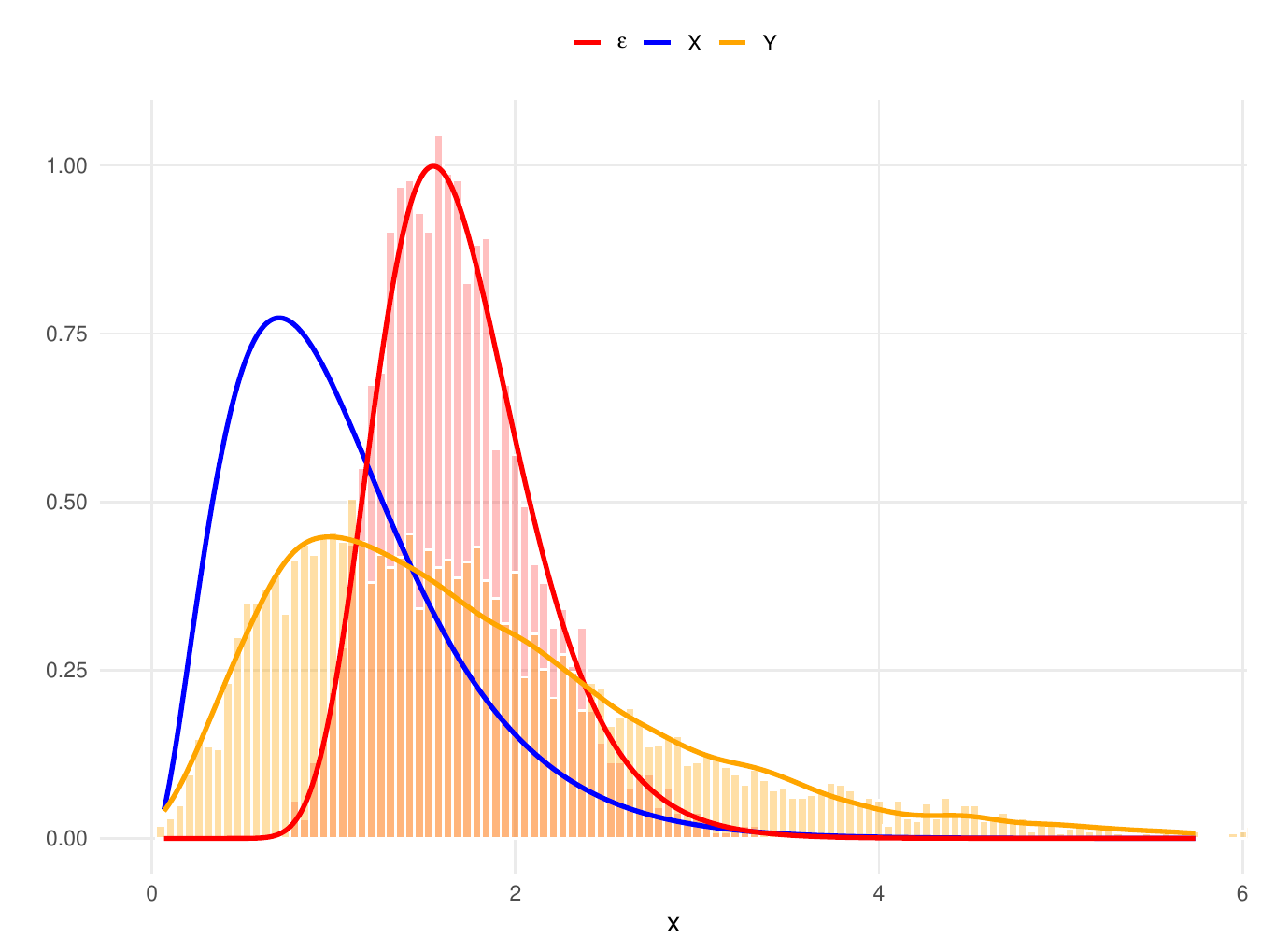}
        \caption{Exponentially smooth model: $X\sim f_{\operatorname{Gamma}(3,7/20)}$ and $\varepsilon\sim f_{\operatorname{LN}(1/2,1/16)}$.}
        \label{fig:simulation-model-es}
    \end{subfigure}
    \caption{Statistical models used in the simulation study. The blue and red
    curves show the densities of \(X\) and \(\varepsilon\), respectively, while the orange curve represents the density of
    \(Y=X\varepsilon\). The orange and red histograms correspond to independent
    samples of \(Y\) and \(\varepsilon\).}
    \label{fig:simulation-models}
\end{figure}

\Cref{fig:simulation-models} illustrates the two statistical models.
The density of \(X\) and the density of \(\varepsilon\) are evaluated exactly.
Since the density of the product \(Y=X\varepsilon\) does not have a convenient
closed form for all model combinations, it is represented by standard kernel estimates
based on a large independent sample. The displayed histograms use
separate samples: \(N= 5000\) observations of \(Y\) and an auxiliary sample of
\(M=2000\) observations from the error distribution.

\subsection{Known error distribution}
\label{sec:simulation-known-error}

\subsubsection{Effect of the sample size for a fixed cut-off}\label{subsec:known:sample_size_effect}

We first fix the regularisation parameter at \(k=0.5\) and consider the sample
sizes
\begin{equation}
    N\in\{100, 500, 1000, 5000\}.
\end{equation}
This experiment isolates the effect of the sample size, 
while holding the
deterministic regularisation bias fixed. We expect
the stochastic error to decrease with \(N\), with a slower and less stable improvement in
the exponentially smooth model.\\

Each experiment consists of \(R=250\) independent Monte Carlo repetitions, where we investigate both
targets
\(
\vartheta(x)\in\{f_X(x),S_X(x)\}.
\)
To quantify the estimation accuracy, we evaluate the pointwise risk (MSE) at \(x_0=0.8\).
Let
\(
\widehat\vartheta_k^{(r)}(x_0)
\)
denote the estimate obtained in repetition \(r\). 
The Monte Carlo MSE is estimated directly by
\begin{equation*}
    \widehat{\operatorname{MSE}}(x_0;k,N)
    =\frac{1}{R}\sum_{r=1}^R
      \left(
        \widehat\vartheta_k^{(r)}(x_0)-\vartheta(x_0)
      \right)^2.
\end{equation*}

For the graphical comparison, all \(R\) estimated curves are shown in grey,
the true target is shown in black, and the pointwise median of the Monte Carlo
estimates is shown in red. Density and survival-function estimates are stacked
vertically within the same figure. This common layout is used throughout the
study.

\begin{figure}[H]
    \begin{subfigure}[t]{0.45\textwidth}
        \centering
        \includegraphics[width=\linewidth]{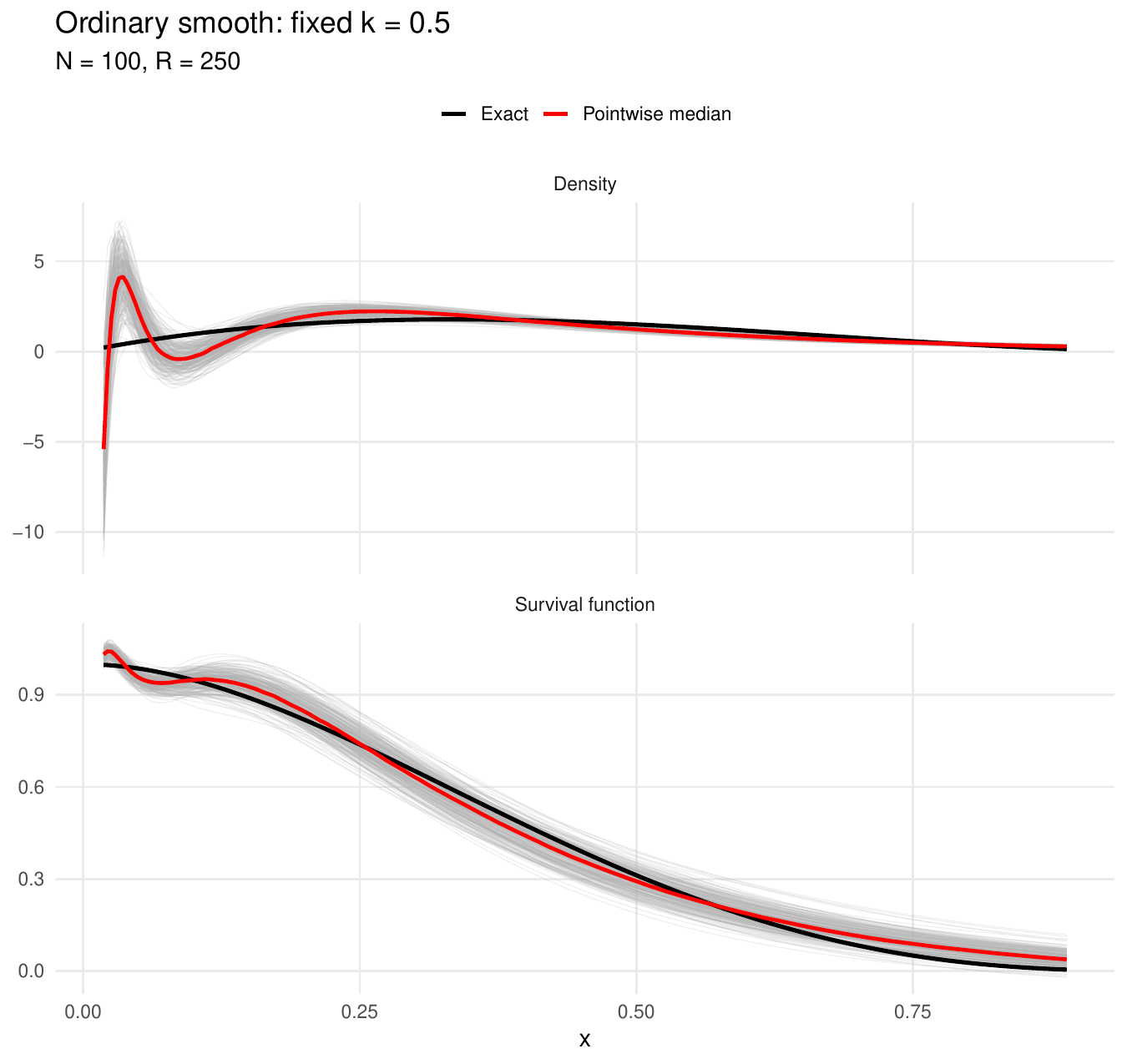}
        \caption{OS, \(N=100\).}
    \end{subfigure}
    \hfill
    \begin{subfigure}[t]{0.45\textwidth}
        \centering
        \includegraphics[width=\linewidth]{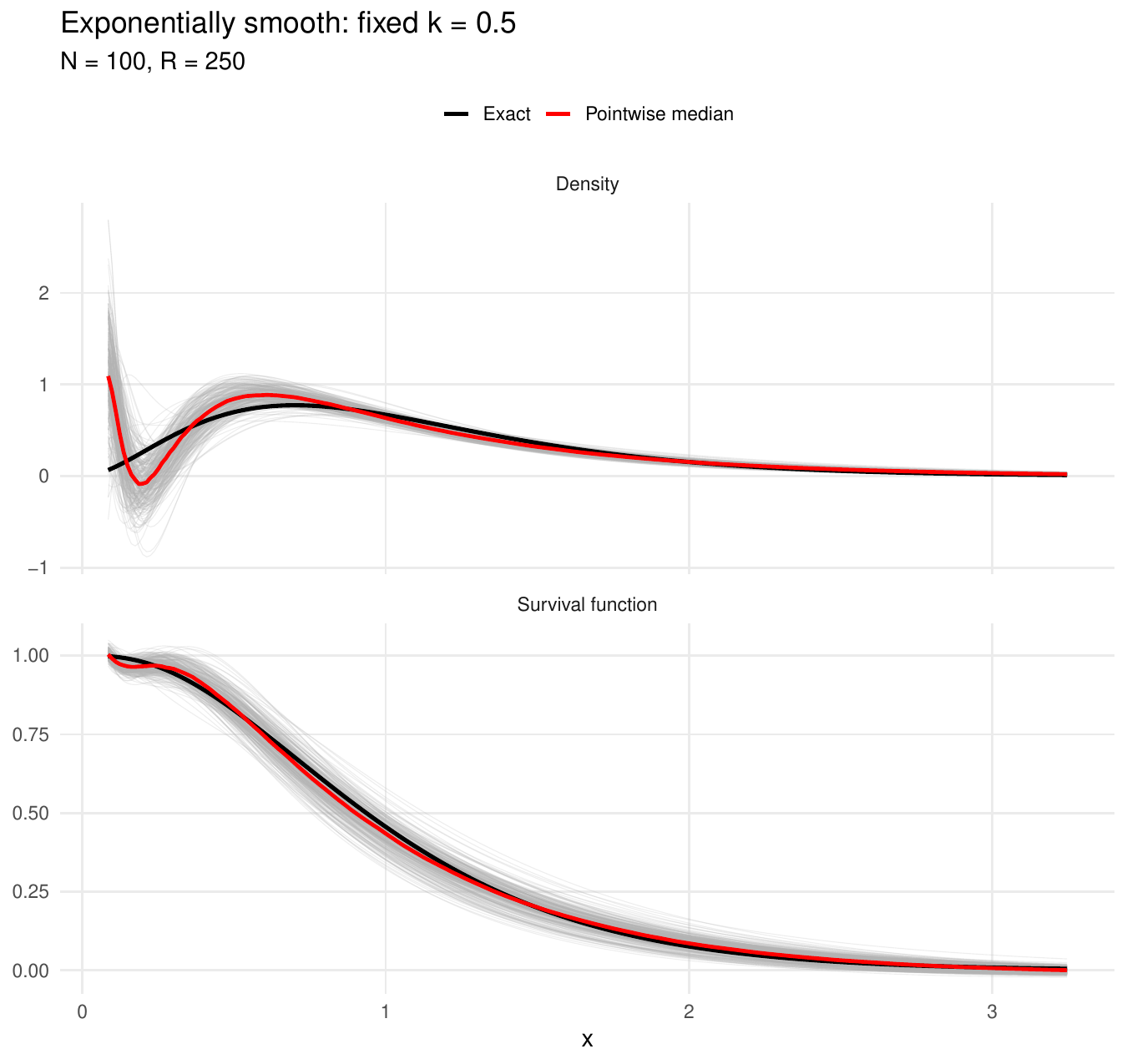}
        \caption{ES, \(N=100\).}
    \end{subfigure}

    \medskip
    \begin{subfigure}[t]{0.45\textwidth}
        \centering
        \includegraphics[width=\linewidth]{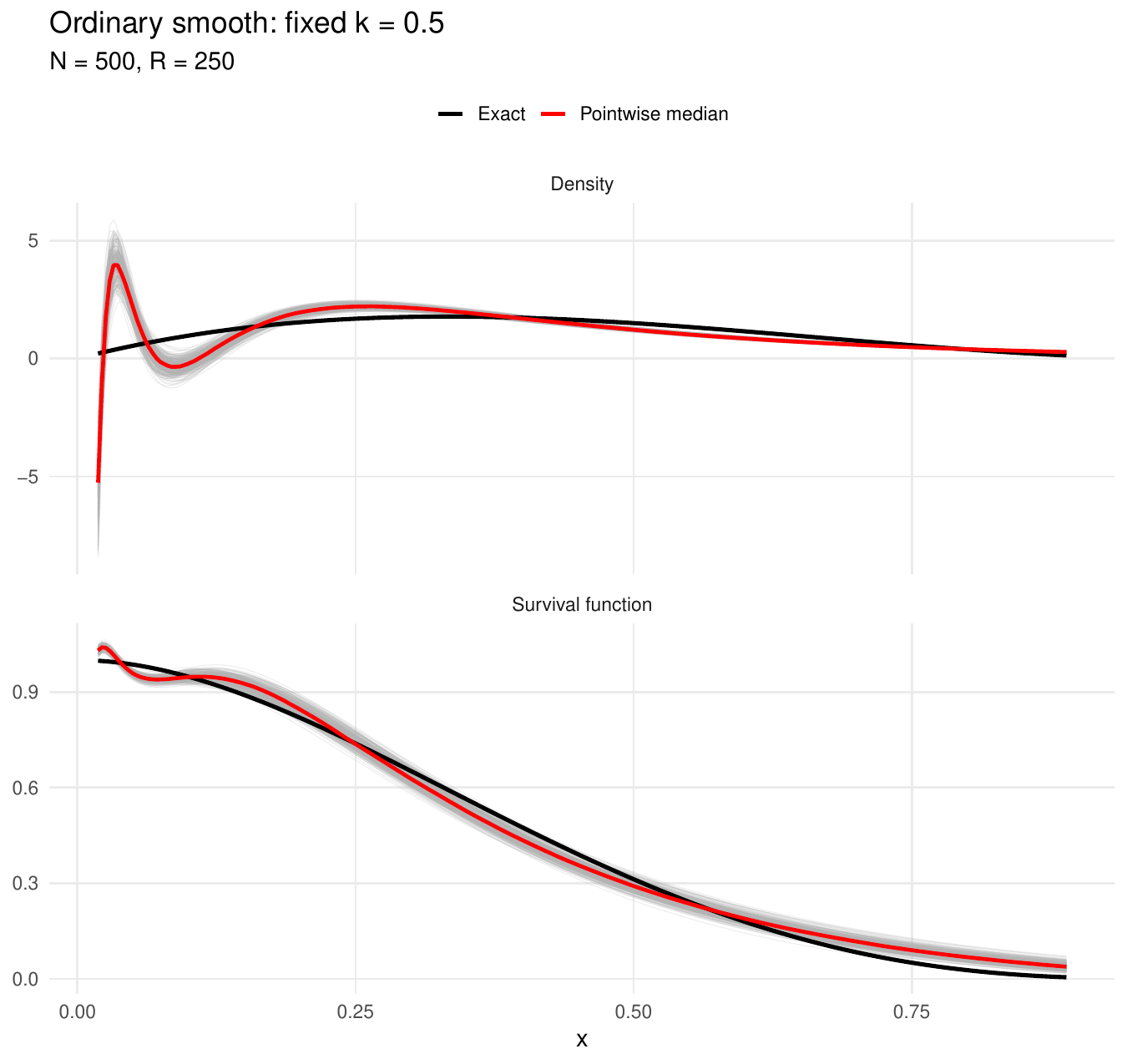}
        \caption{OS, \(N=500\).}
    \end{subfigure}
    \hfill
    \begin{subfigure}[t]{0.45\textwidth}
        \centering
        \includegraphics[width=\linewidth]{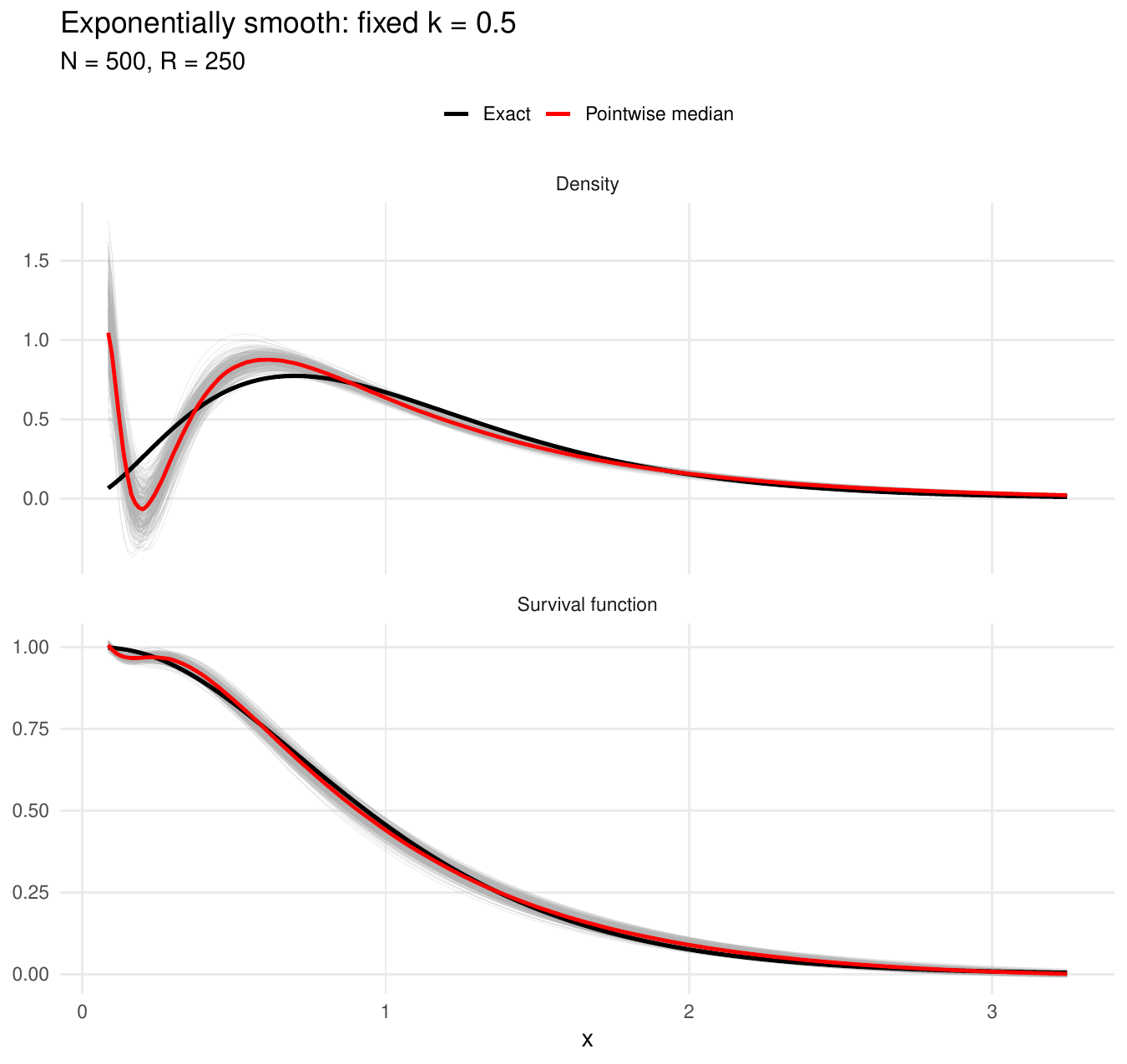}
        \caption{ES, \(N=500\).}
    \end{subfigure}
\end{figure}
\begin{figure}[H]
	\ContinuedFloat
    \begin{subfigure}[t]{0.45\textwidth}
        \centering
        \includegraphics[width=\linewidth]{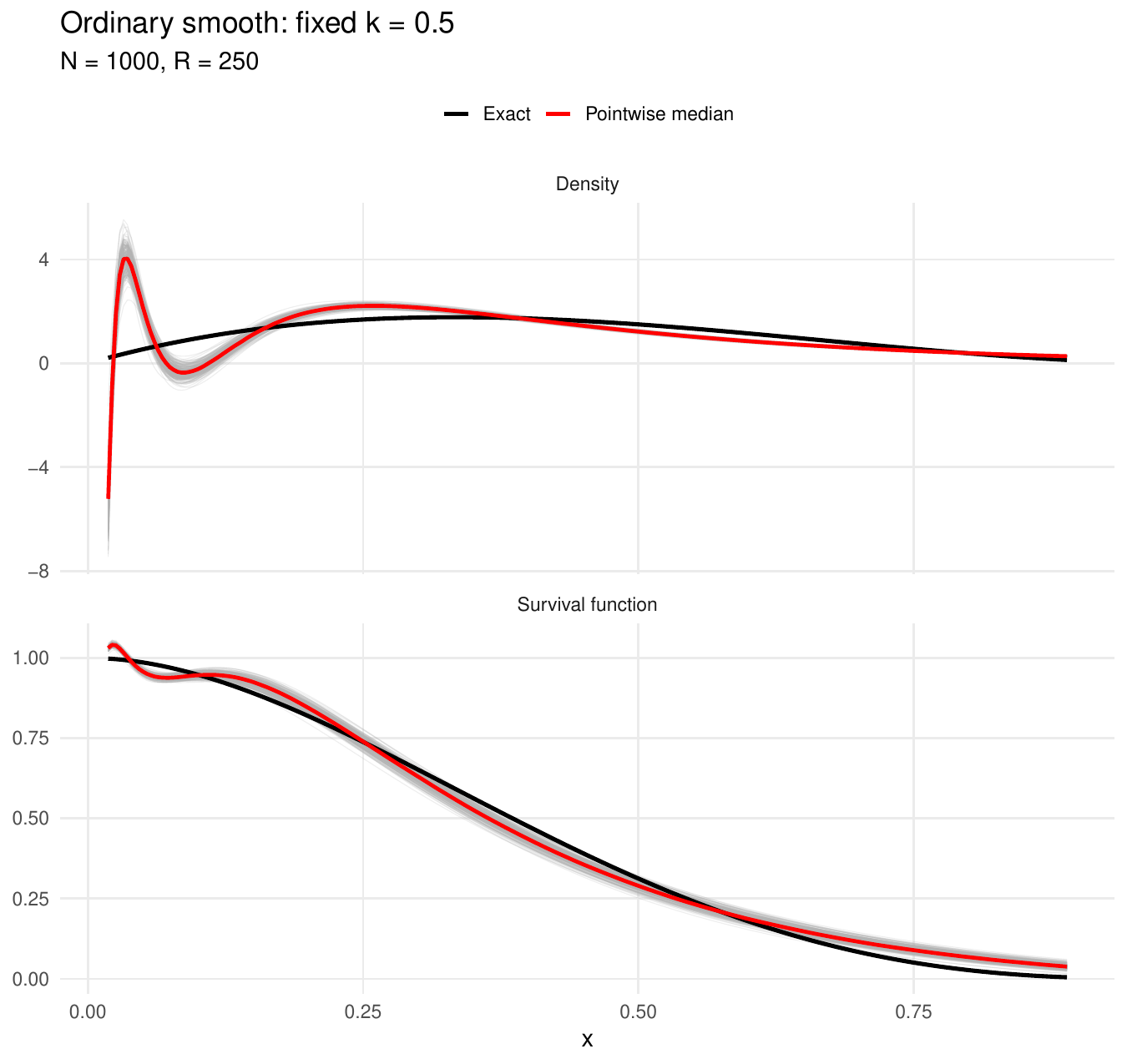}
        \caption{OS, \(N=1000\).}
    \end{subfigure}
    \hfill
    \begin{subfigure}[t]{0.45\textwidth}
        \centering
        \includegraphics[width=\linewidth]{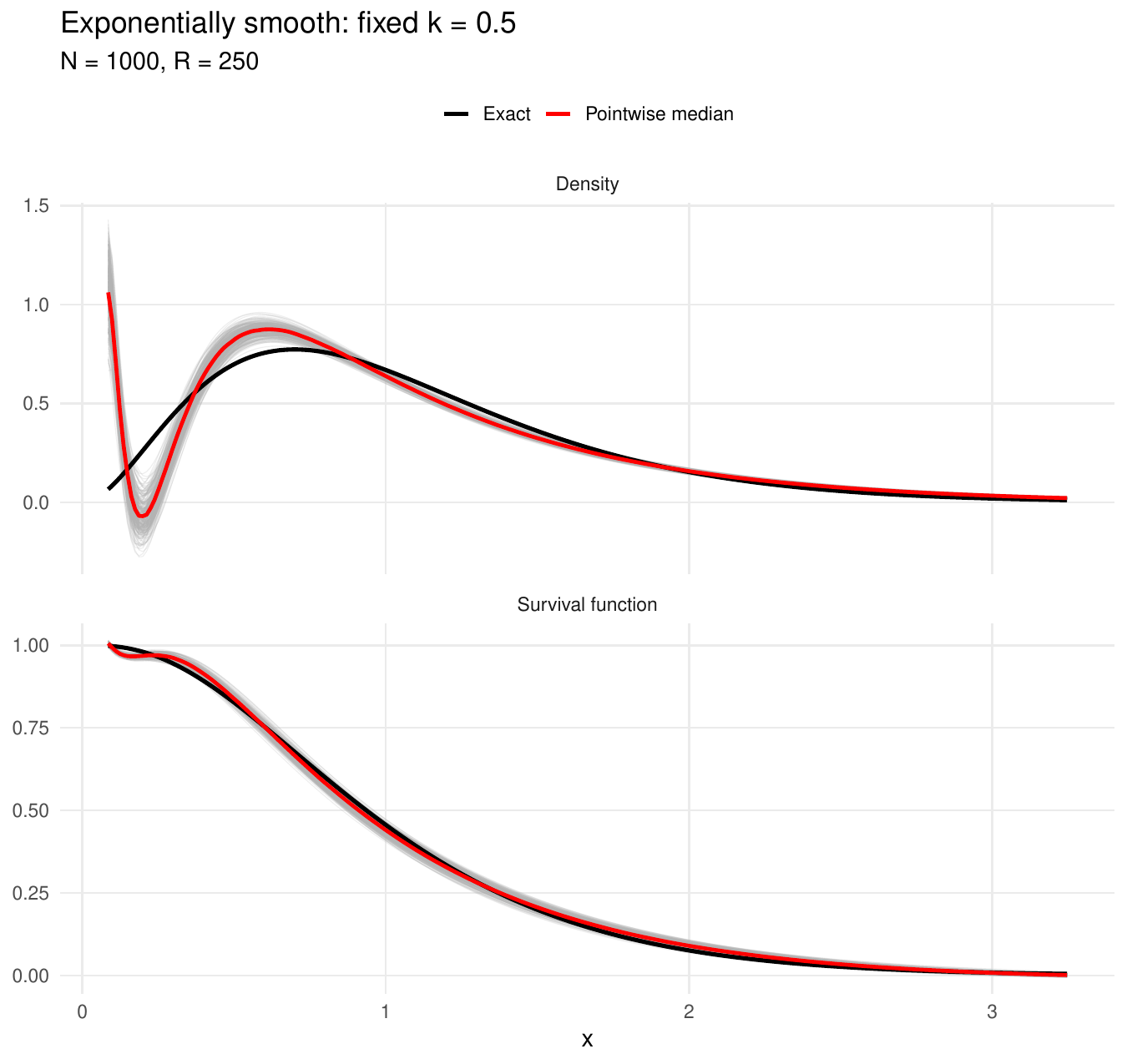}
        \caption{ES, \(N=1000\).}
    \end{subfigure}
	\end{figure}

\begin{figure}[H]
	\ContinuedFloat
	\begin{subfigure}[t]{0.45\textwidth}
        \centering
        \includegraphics[width=\linewidth]{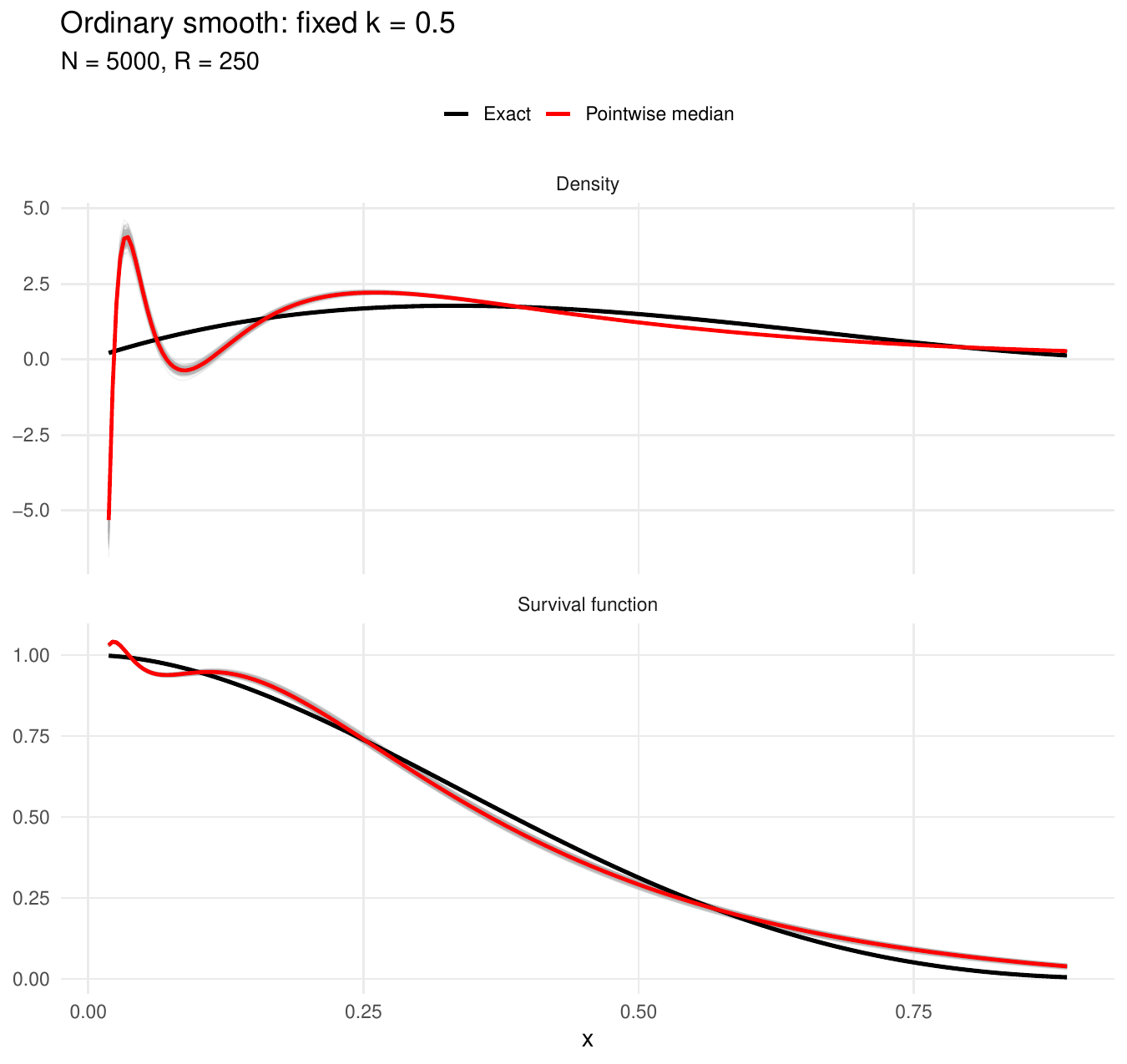}
        \caption{OS, \(N=5000\).}
    \end{subfigure}
    \hfill
    \begin{subfigure}[t]{0.45\textwidth}
        \centering
        \includegraphics[width=\linewidth]{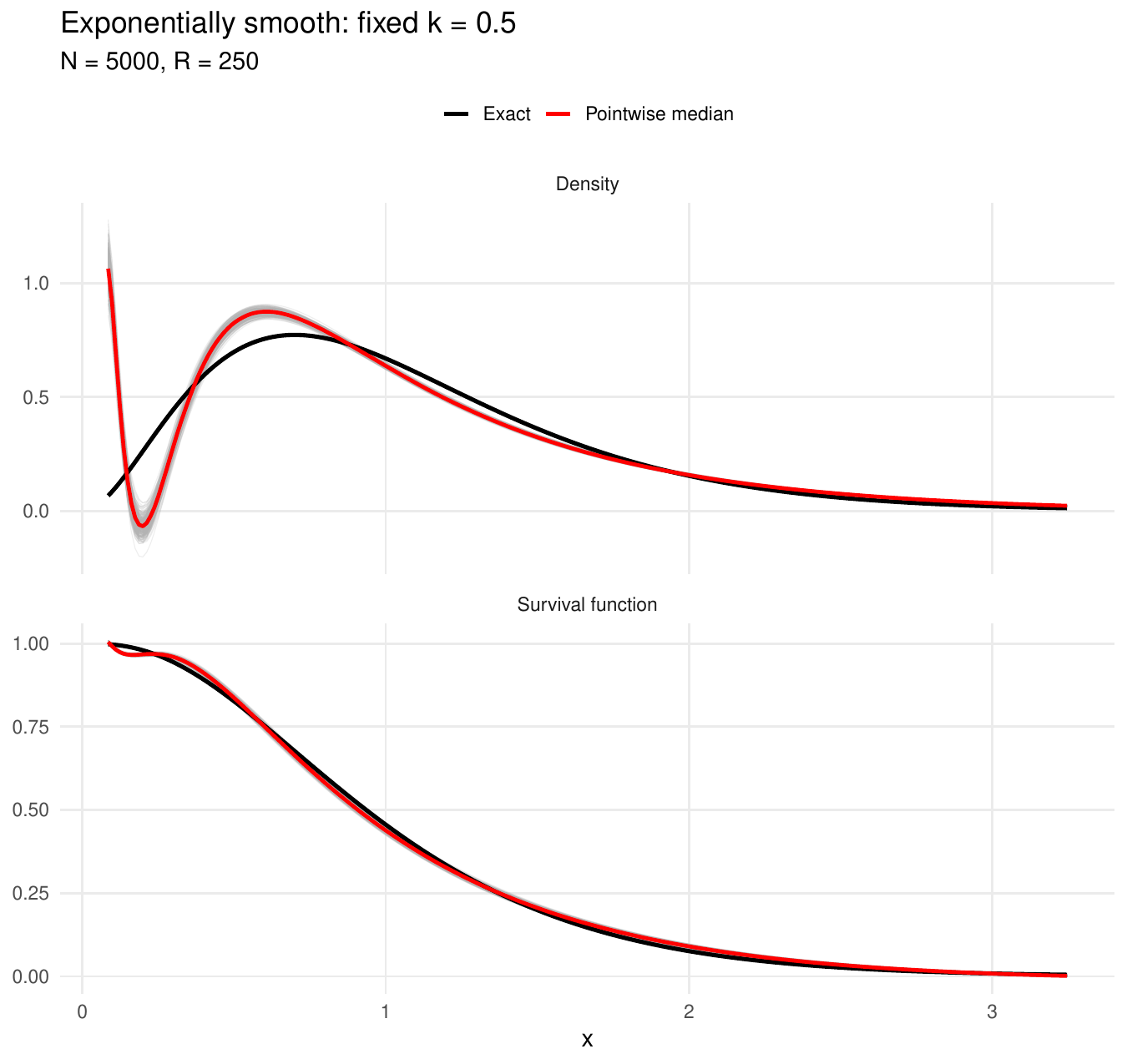}
        \caption{ES, \(N=5000\).}
    \end{subfigure}
    \caption{Monte Carlo estimates for a fixed cut-off \(k=0.5\). Each panel
    contains the density estimate above the survival-function estimate. Grey
    lines represent the \(R=250\) repetitions, the black line is the true
    target, and the red line is the pointwise Monte Carlo median.}
\end{figure}

\begin{table}[htbp]
    \centering
    \small
    \caption{Pointwise Monte Carlo MSE at \(x_0=0.8\) for the
             fixed cut-off \(k=0.5\).}
    \label{tab:pointwise-mse-known}

    \begin{tabular}{llcccc}
        \toprule
        & & \multicolumn{4}{c}{
            \(\widehat{\mathrm{MSE}}(0.8;0.5,\cdot)\;\times 10^{-3}\)
        } \\
        \cmidrule(lr){3-6}
        \textbf{Smoothness}
        & \textbf{Target}
        & \(N=100\)
        & \(N=500\)
        & \(N=1000\)
        & \(N=5000\) \\
        \midrule
        OS & Density
           & 2.98 & 0.765 & 0.449 & 0.234 \\
        OS & Survival function
           & 2.34 & 1.87  & 1.74  & 1.70  \\
        \addlinespace
        ES & Density
           & 4.61 & 1.76  & 1.43  & 1.06  \\
        ES & Survival function
           & 2.60 & 0.726 & 0.476 & 0.361 \\
        \bottomrule
    \end{tabular}
\end{table}

\subsubsection{Effect of the cut-off for a fixed sample size}

Next, we fix \(N=1000\) and vary the cut-off over
\begin{equation}
    k\in\{0.1,0.5,1.0,2.0\}.
    \label{eq:k-grid-fixed-n}
\end{equation}
This experiment directly illustrates the regularization trade-off.
For simplicity, we restrict ourselves here to the ordinary smoothness case. 

\begin{figure}[H]
    \centering
    \begin{subfigure}[t]{0.45\textwidth}
        \centering
        \includegraphics[width=\linewidth]{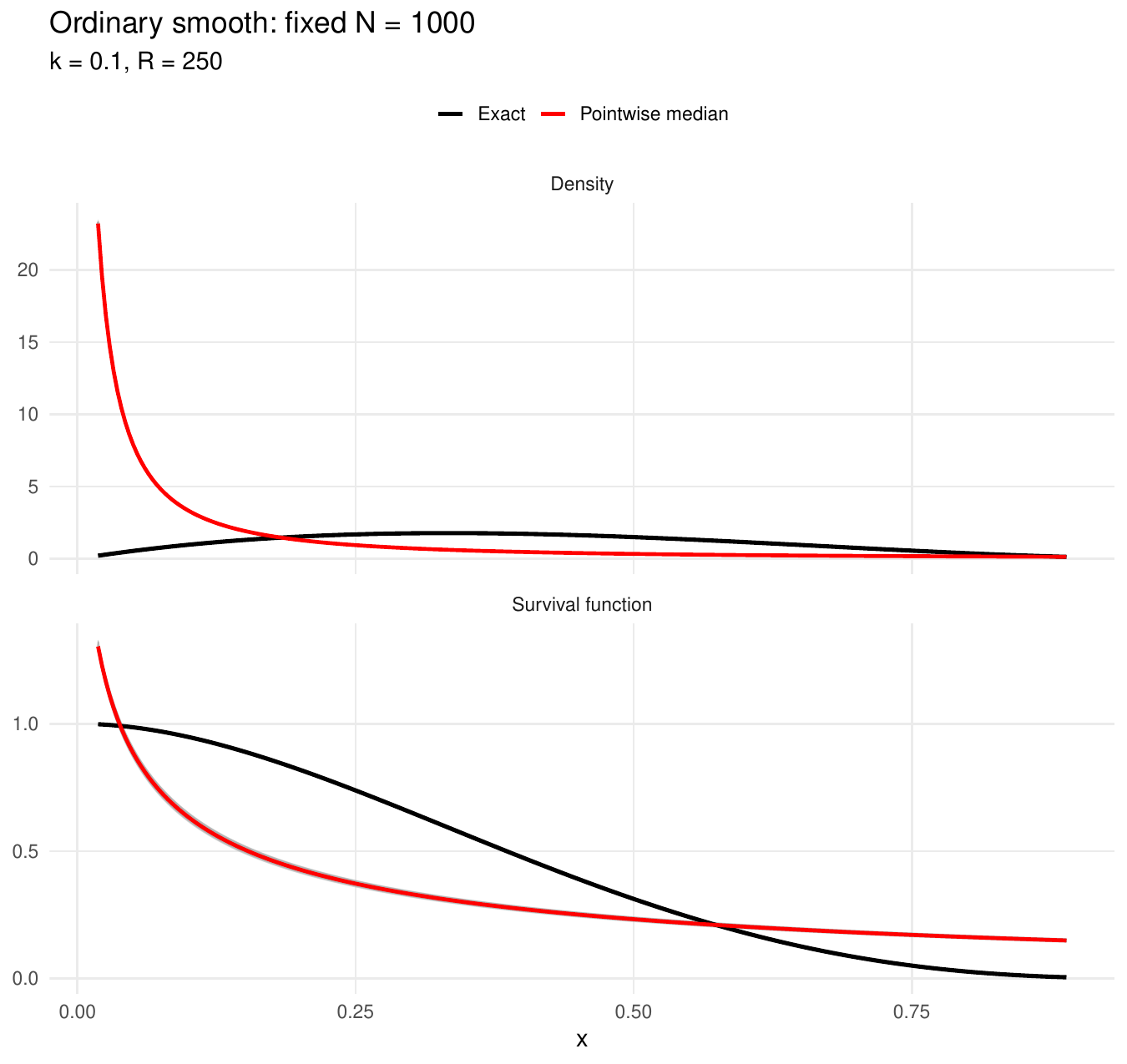}
        \caption{OS, \(N=1000\), $k = 0.1$.}
    \end{subfigure}
    \hfill
    \begin{subfigure}[t]{0.45\textwidth}
        \centering
        \includegraphics[width=\linewidth]{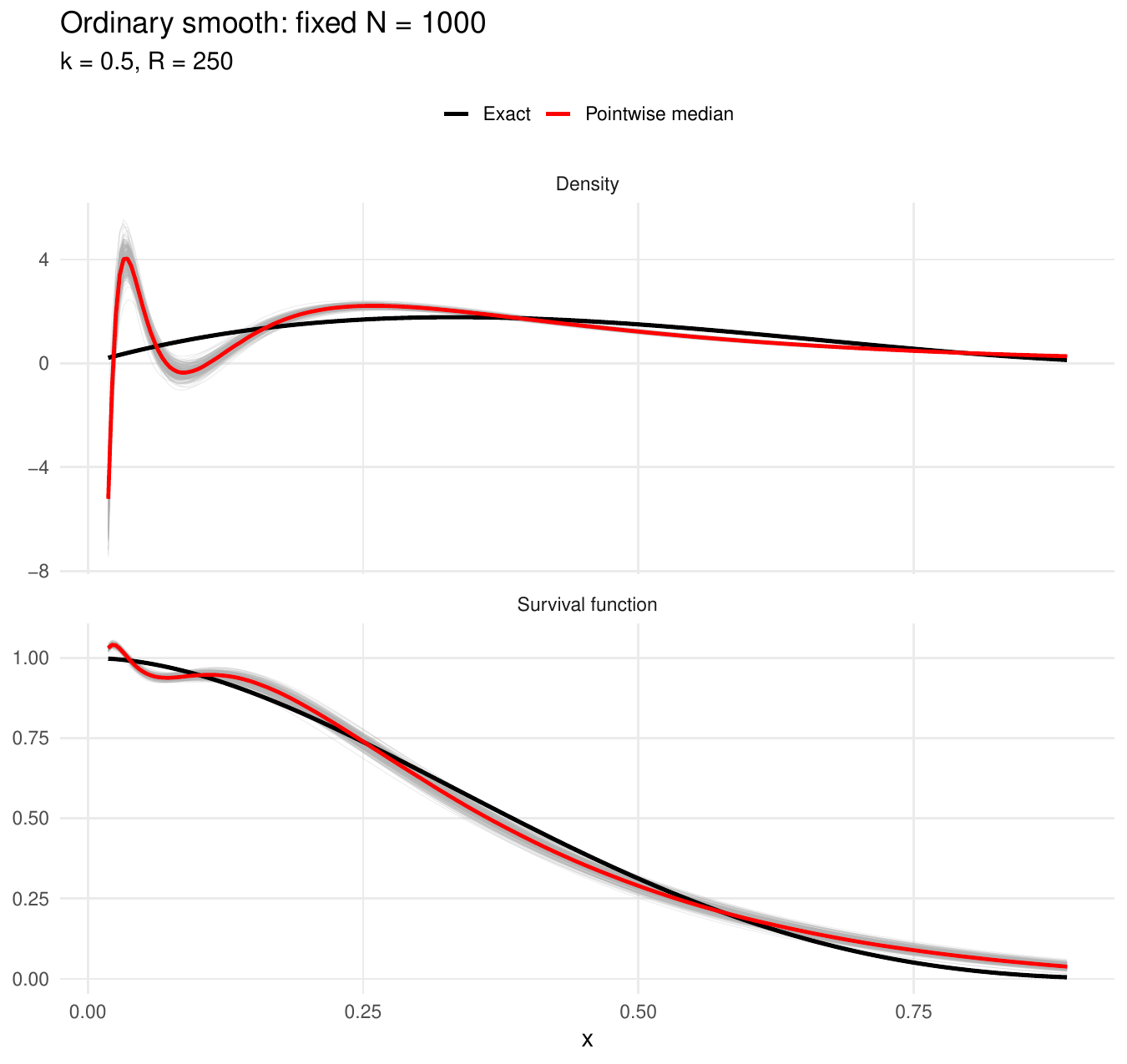}
       \caption{OS, \(N=1000\), $k = 0.5$.}
    \end{subfigure}

    \medskip
    \begin{subfigure}[t]{0.45\textwidth}
        \centering
        \includegraphics[width=\linewidth]{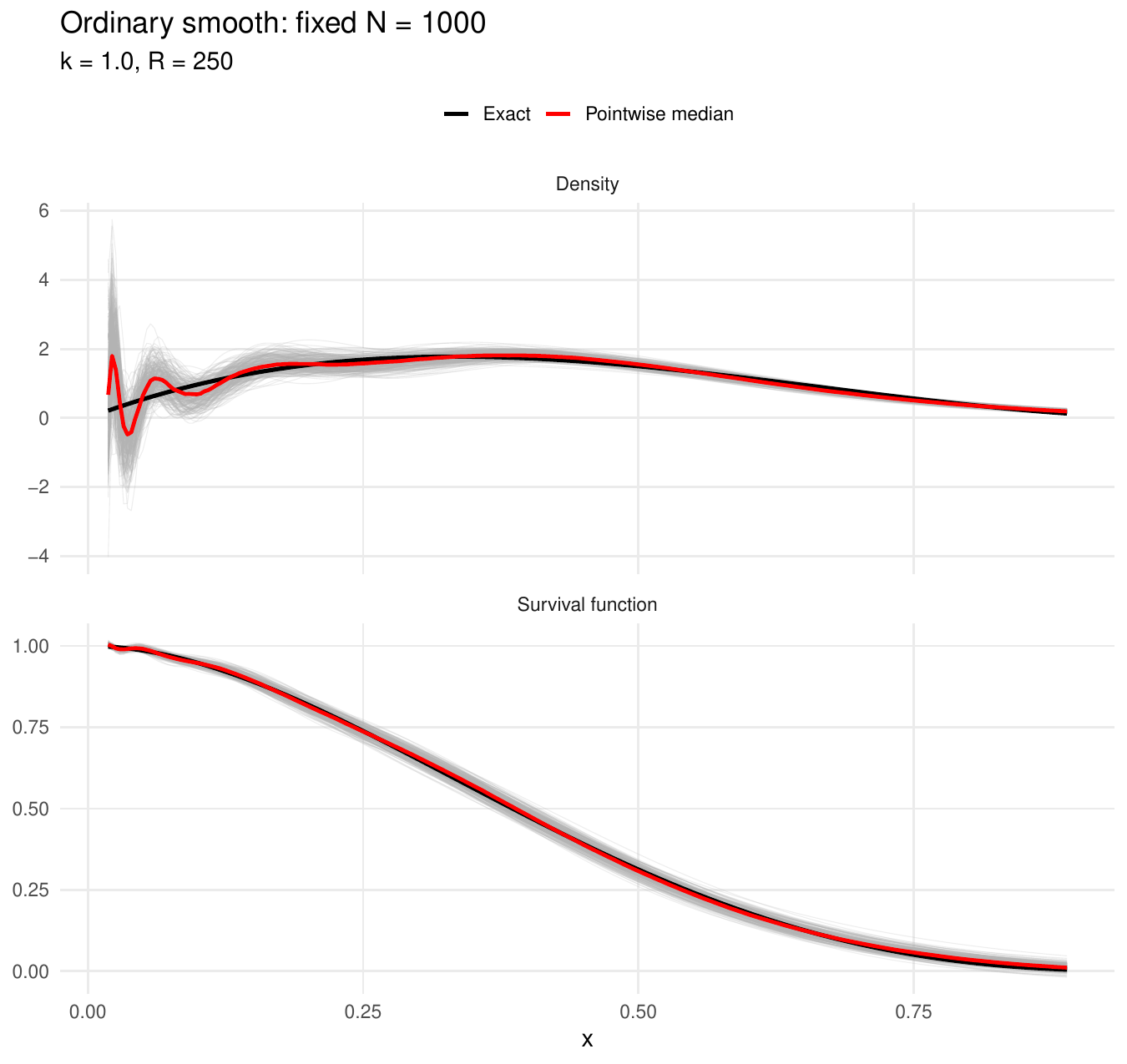}
       \caption{OS, \(N=1000\), $k = 1.0$.}
    \end{subfigure}
    \hfill
    \begin{subfigure}[t]{0.45\textwidth}
        \centering
        \includegraphics[width=\linewidth]{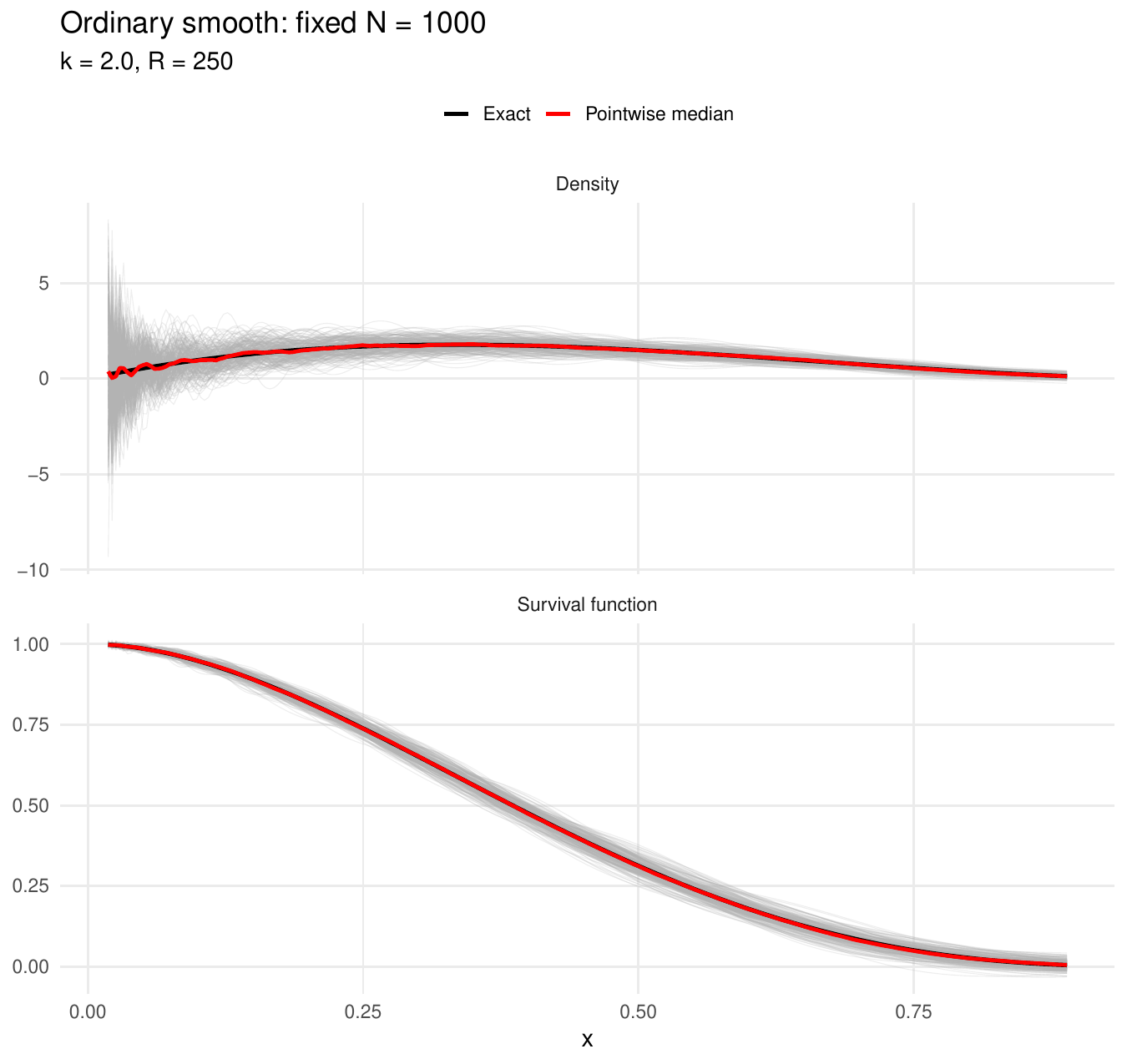}
        \caption{OS, \(N=1000\), $k = 2.0$.}
    \end{subfigure}
    \caption{Representative cut-off values for fixed \(N=1000\).}
    \label{fig:fixed-n-by-k}
\end{figure}

\begin{table}[htbp]
    \centering
    \small
    \caption{Pointwise Monte Carlo MSE at \(x_0=0.8\) for the
             fixed cut-off \(N=1000\).}
    \label{tab:pointwise-mse-cutoff}

    \begin{tabular}{llcccc}
        \toprule
        & & \multicolumn{4}{c}{
            \(\widehat{\mathrm{MSE}}(0.8;\,\cdot,1000)
            \;(\times 10^{-3})\)
        } \\
        \cmidrule(lr){3-6}
        \textbf{Smoothness}
        & \textbf{Target}
        & \(k=0.1\)
        & \(k=0.5\)
        & \(k=1\)
        & \(k=2\) \\
        \midrule
        OS & Density
           & 49.1 & 0.449 & 1.98  & 19.2  \\
        OS & Survival function
           & 18.4 & 1.74  & 0.197 & 0.232 \\
        \bottomrule
    \end{tabular}
\end{table}

\subsubsection{Data-driven cut-off selection}
\label{sec:data-driven-simulation}

We simulate the data-driven selection rule
developed in \cref{sec:lifumu:datadriven:known} (particuarly \cref{prop:lifumu:neworacle:known}), by choosing $\kappa = 1/100$ and replacing the penalty constant $6$ in \cref{eq:khat:known:final} by the factor $0.6$.
The selection is performed for the evaluation point \(x_0=0.8\)
and is carried out separately for the probability density function and the survival function.\\

\begin{figure}[H]
    \centering
    \begin{subfigure}[t]{0.48\textwidth}
        \centering
        \includegraphics[width=\linewidth]{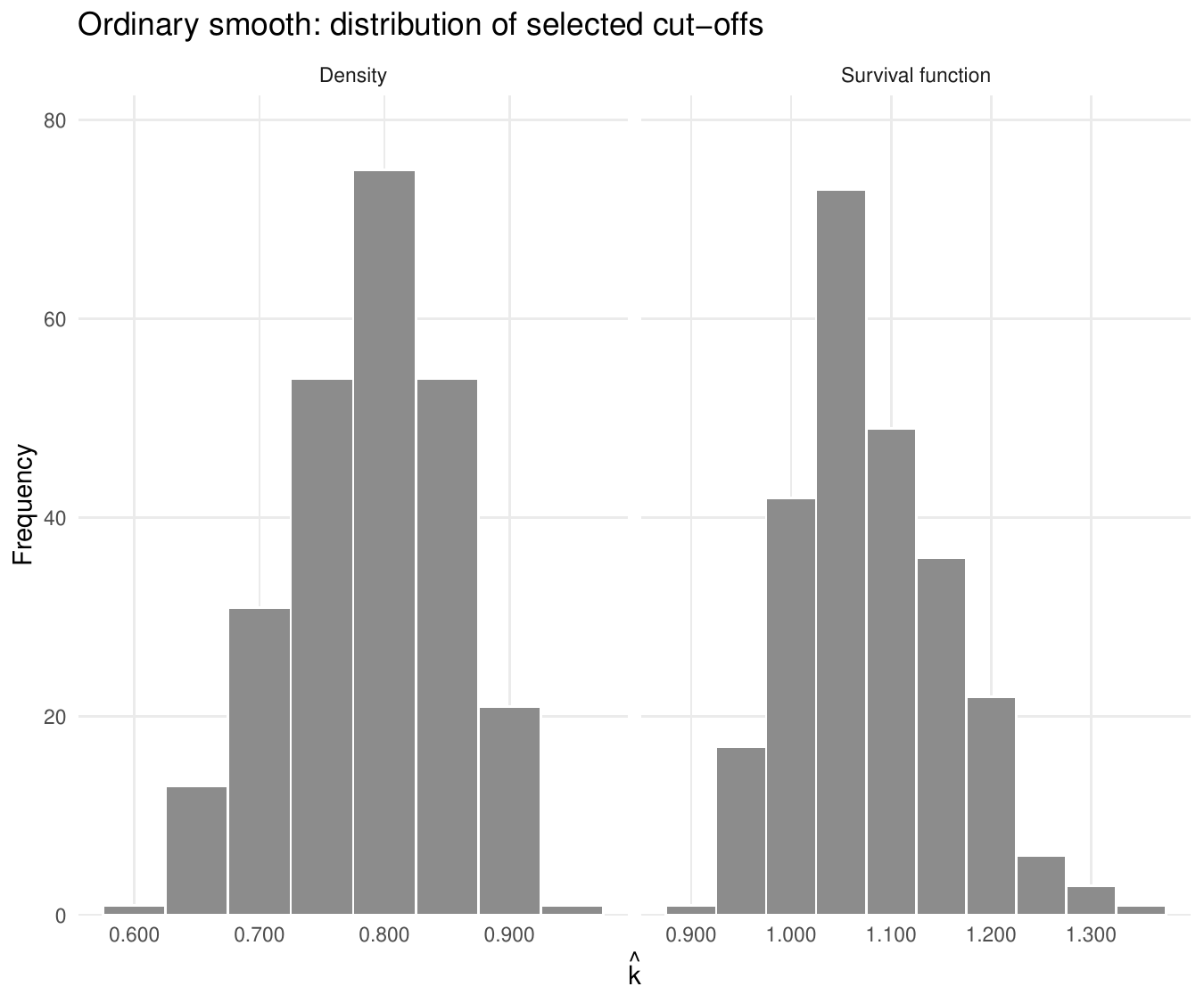}
        \caption{Data-driven cut-off selection distribution: Ordinary Smoothness.}
        \label{fig:selected-curves-placeholder}
    \end{subfigure}
    \hfill
    \begin{subfigure}[t]{0.48\textwidth}
        \centering
        \includegraphics[width=\linewidth]{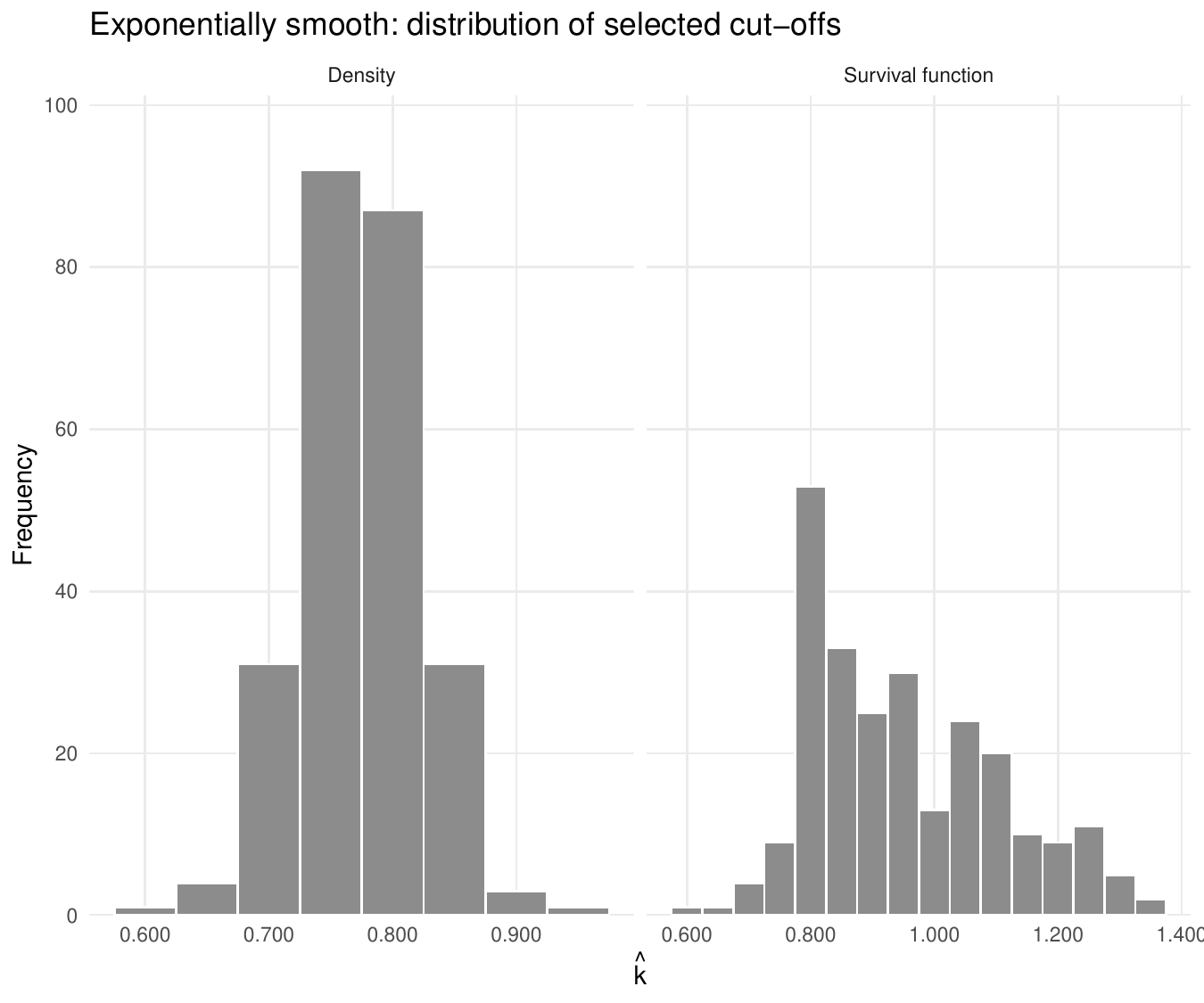}
        \caption{Data-driven cut-off selection distribution: Exponentially Smoothness.}
        \label{fig:selected-k-placeholder}
    \end{subfigure}
    \caption{Data-driven cut-off selection distributions illustrated by corresponding histograms for $N = 5000$.}
    \label{fig:data-driven-selection}
\end{figure}

For $R = 250$ Monte Carlo repetitions, we record the selected cut-off \(\hat{k}\).
The resulting distributions are illustrated by histograms in \cref{fig:data-driven-selection}.

\begin{table}[htbp]
    \centering
    \small
    \caption{Pointwise Monte Carlo MSE at \(x_0=0.8\) for data-driven cut-off selection.}
    \label{tab:pointwise-mse}

    \begin{tabular}{llcccc}
        \toprule
        & & \multicolumn{4}{c}{
            \(\widehat{\mathrm{MSE}}(0.8; \hat{k},\cdot)\;(\times 10^{-3})\)
        } \\
        \cmidrule(lr){3-6}
        \textbf{Smoothness}
        & \textbf{Target}
        & \(N=100\)
        & \(N=500\)
        & \(N=1000\)
        & \(N=5000\) \\
        \midrule
        OS & Density
           & 7.29 & 1.80  & 0.917 & 0.281 \\
        OS & Survival function
           & 1.28 & 0.380 & 0.285 & 0.0667 \\
        \addlinespace
        ES & Density
           & 8.65 & 3.08  & 1.85  & 0.586 \\
        ES & Survival function
           & 2.48 & 0.531 & 0.317 & 0.0620 \\
        \bottomrule
    \end{tabular}
\end{table}

\subsection{Unknown error distribution}
\label{sec:simulation-unknown-error}
We now consider the second regime \assuRef{SRII} when not knowing the error distribution, but observing 
an independent auxiliary sample
\begin{equation*}
    \varepsilon_1,\ldots,\varepsilon_M
    \overset{\mathrm{i.i.d.}}{\sim}f_\varepsilon.
\end{equation*}

\subsubsection{Effect of the additional sample size for a fixed cut-off}
As in \cref{subsec:known:sample_size_effect}, we first fix the regularization parameter at \(k=0.5\).
We fix the sample size $N= 5000$ and consider additional sample sizes
\begin{equation}
    M\in\{100, 500, 1000, 5000\}.
    \label{eq:fixed-k-sample-sizes}
\end{equation}
As before we simulate $\hat{f}_k$ and $\SXhat_k$ using a Monte-Carlo simulation.

\begin{figure}[H]
    \centering
    \begin{subfigure}[t]{0.45\textwidth}
        \centering
        \includegraphics[width=\linewidth]{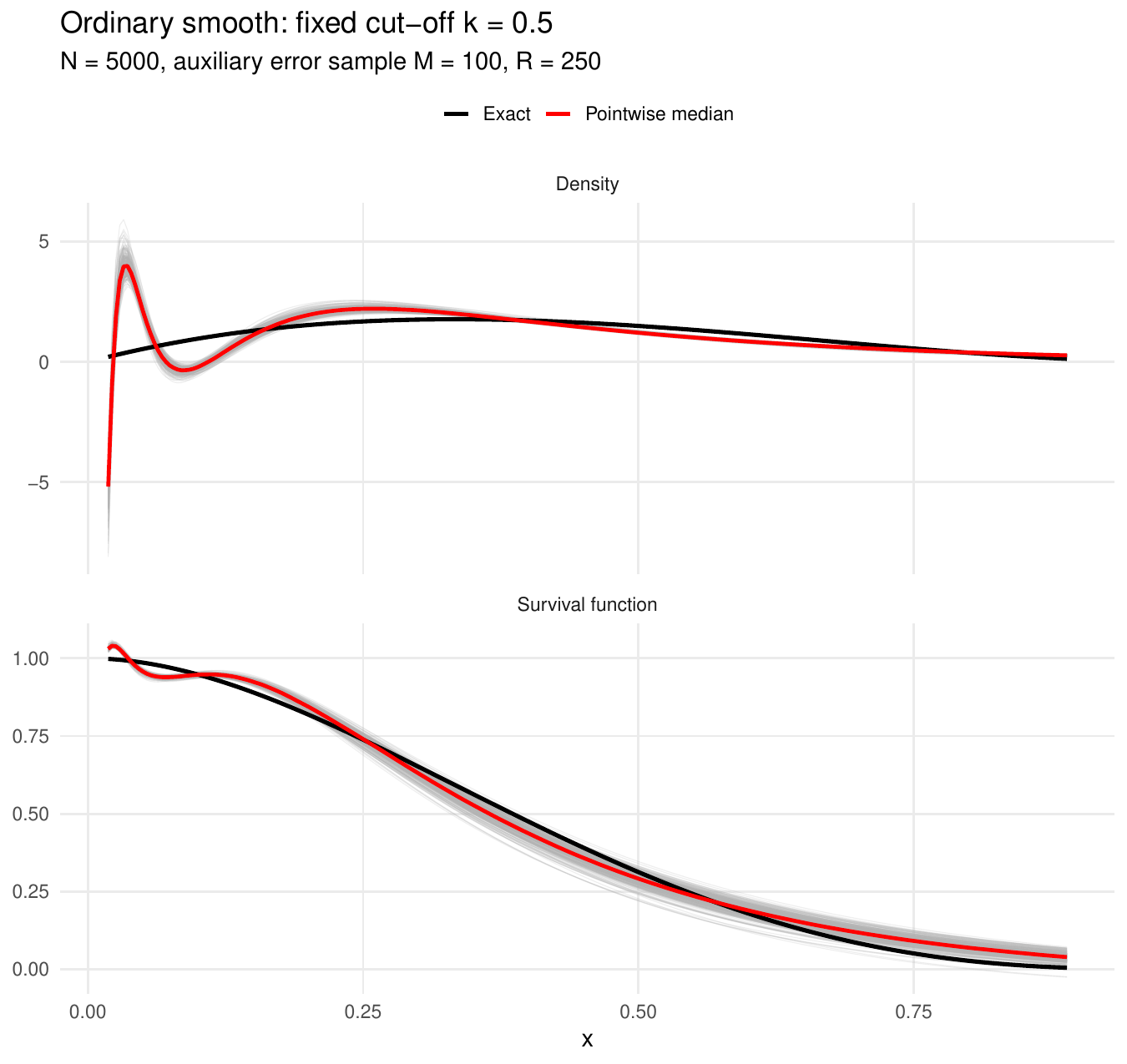}
        \caption{OS, $N=5000$, \(M=100\).}
    \end{subfigure}
    \hfill
    \begin{subfigure}[t]{0.45\textwidth}
        \centering
        \includegraphics[width=\linewidth]{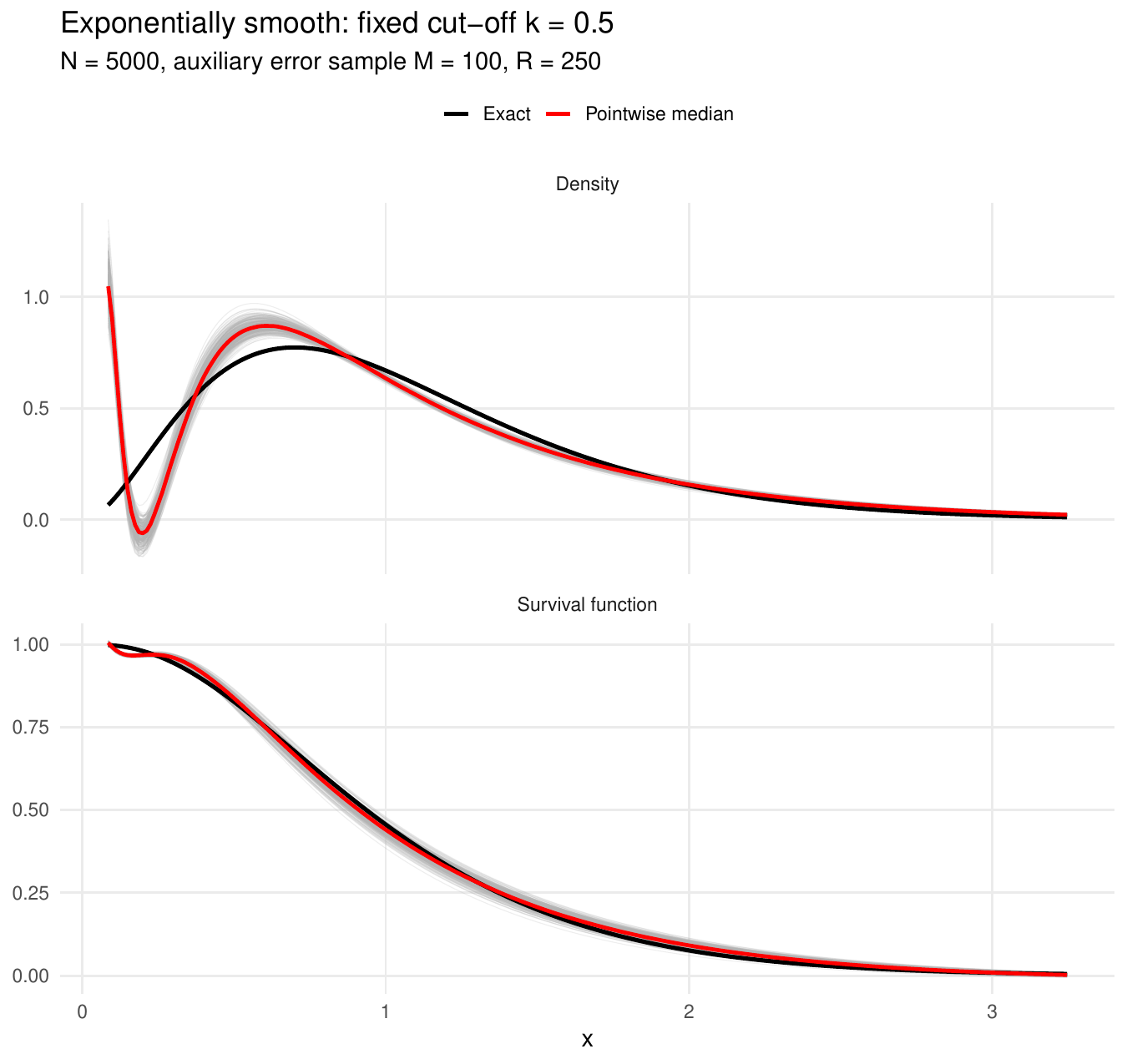}
        \caption{ES,  $N=5000$, \(M=100\).}
    \end{subfigure}

    \medskip
    \begin{subfigure}[t]{0.45\textwidth}
        \centering
        \includegraphics[width=\linewidth]{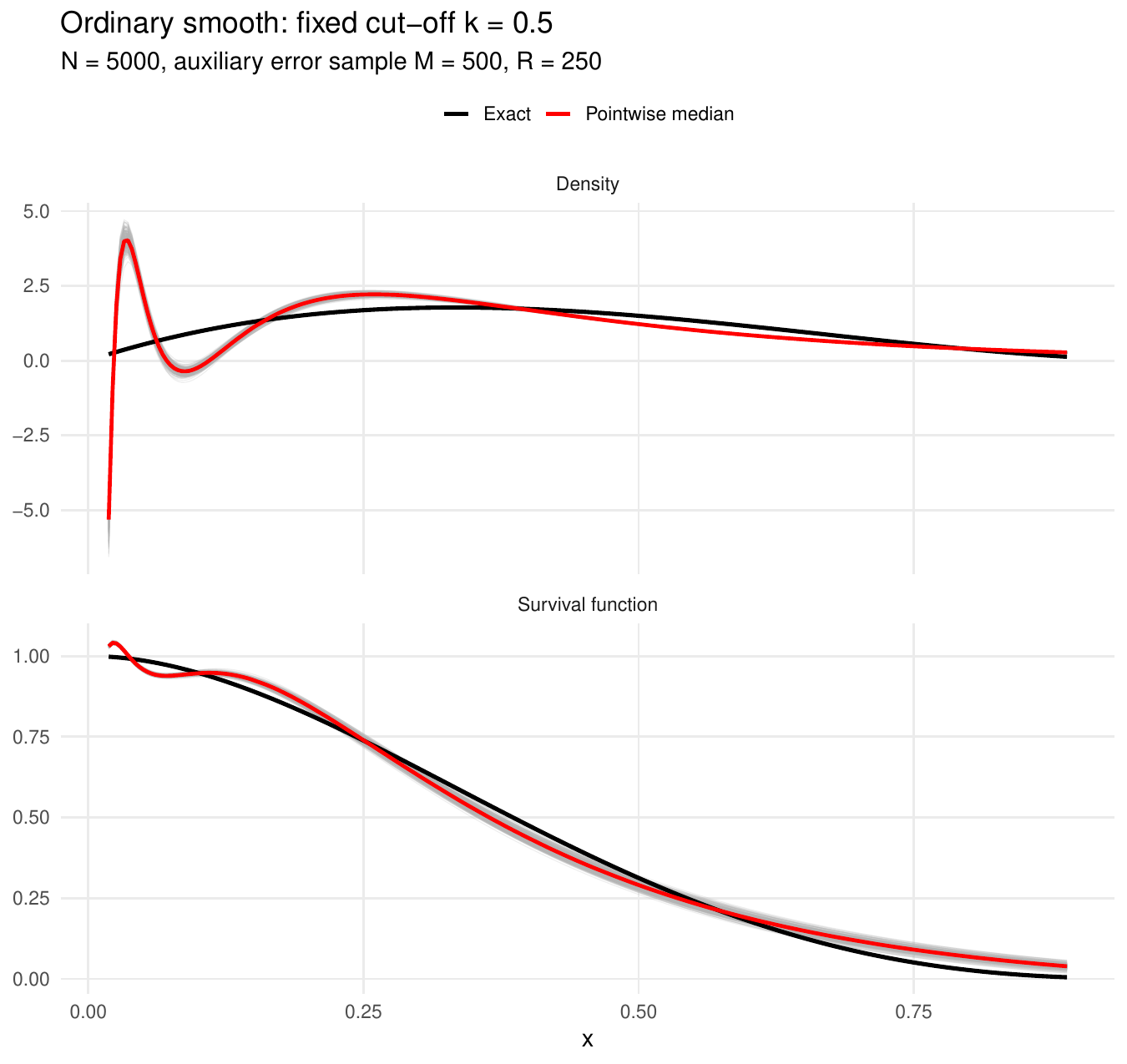}
        \caption{OS, $N=5000$, \(M=500\).}
    \end{subfigure}
    \hfill
    \begin{subfigure}[t]{0.45\textwidth}
        \centering
        \includegraphics[width=\linewidth]{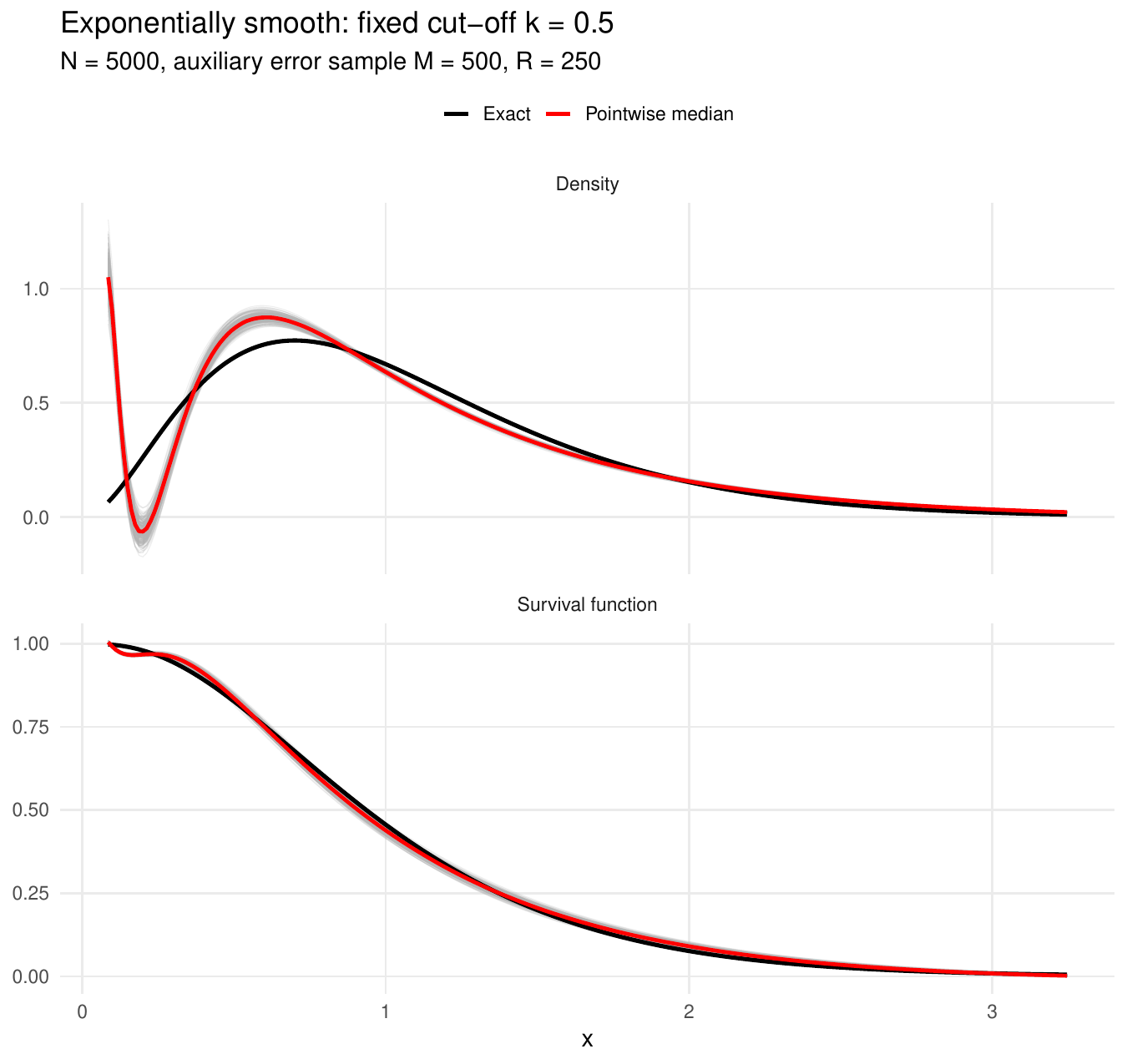}
        \caption{ES, $N=5000$, \(M=500\).}
    \end{subfigure}

    	\end{figure}

\begin{figure}[H]
	\ContinuedFloat
    \begin{subfigure}[t]{0.45\textwidth}
        \centering
        \includegraphics[width=\linewidth]{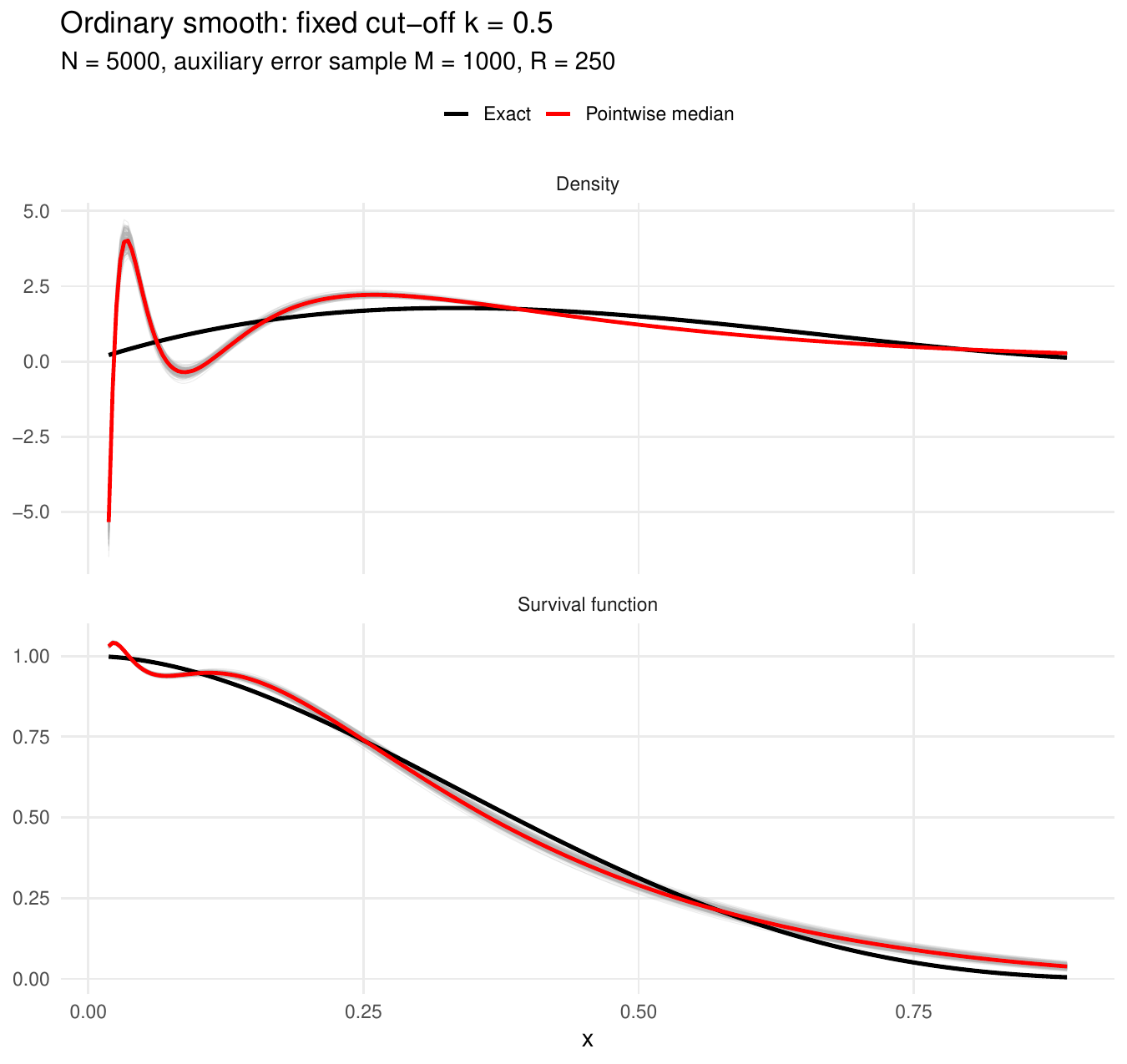}
        \caption{OS, $N=5000$, \(M=1000\).}
    \end{subfigure}
    \hfill
    \begin{subfigure}[t]{0.45\textwidth}
        \centering
        \includegraphics[width=\linewidth]{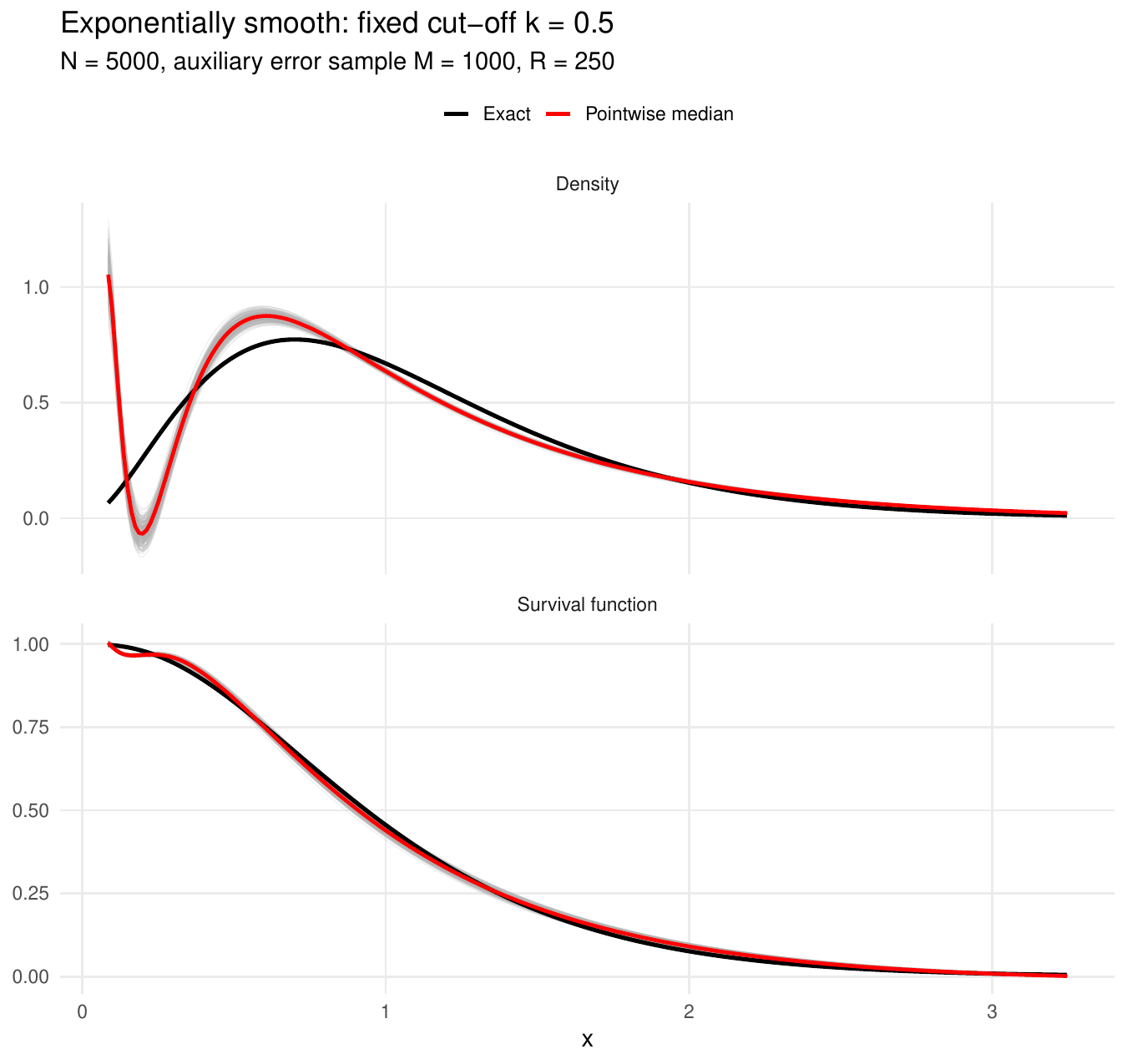}
        \caption{ES, $N=5000$, \(M=1000\).}
    \end{subfigure}
	\medskip
	\begin{subfigure}[t]{0.45\textwidth}
        \centering
        \includegraphics[width=\linewidth]{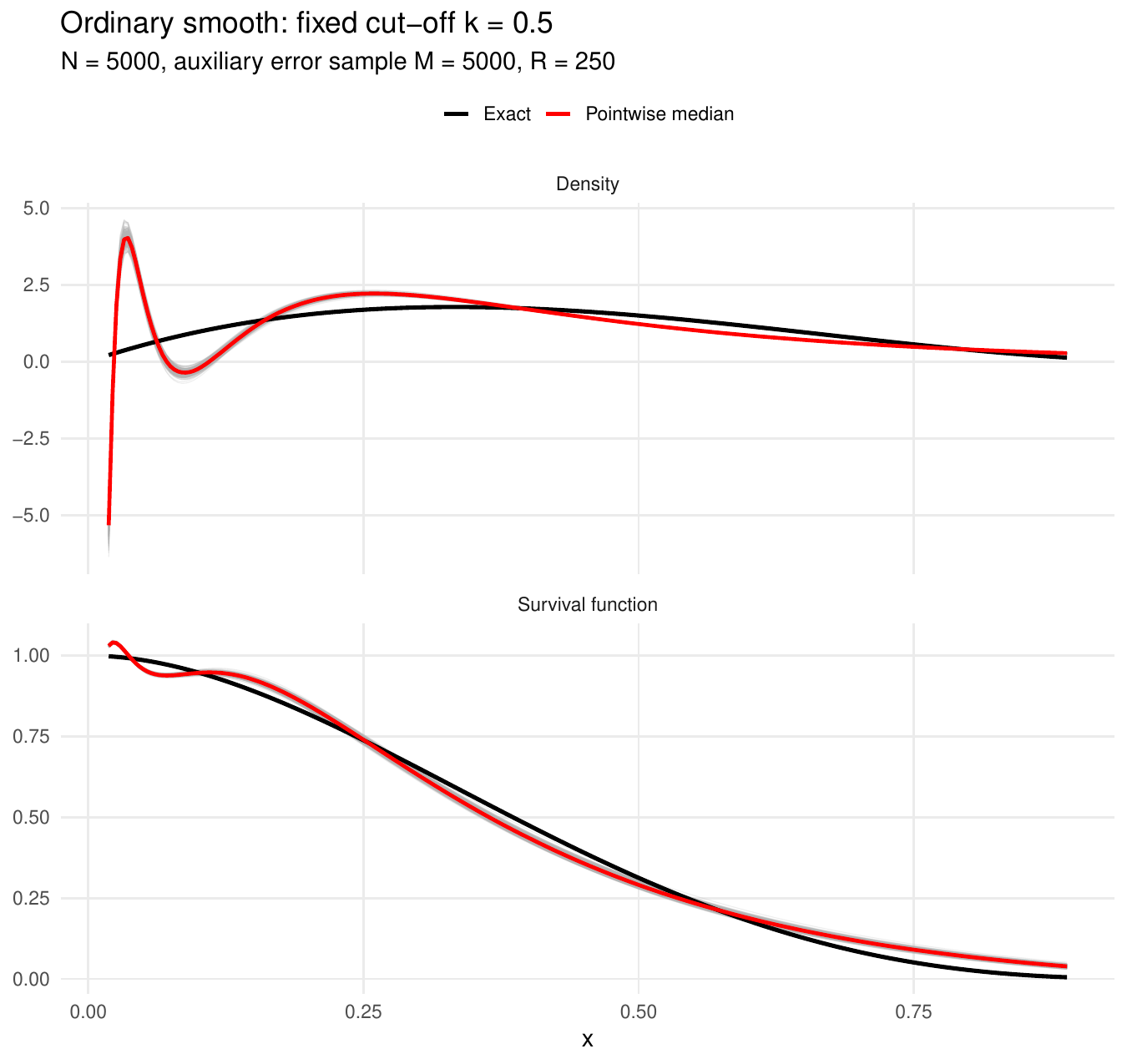}
        \caption{OS, $N=5000$, \(M=5000\).}
    \end{subfigure}
    \hfill
    \begin{subfigure}[t]{0.45\textwidth}
        \centering
        \includegraphics[width=\linewidth]{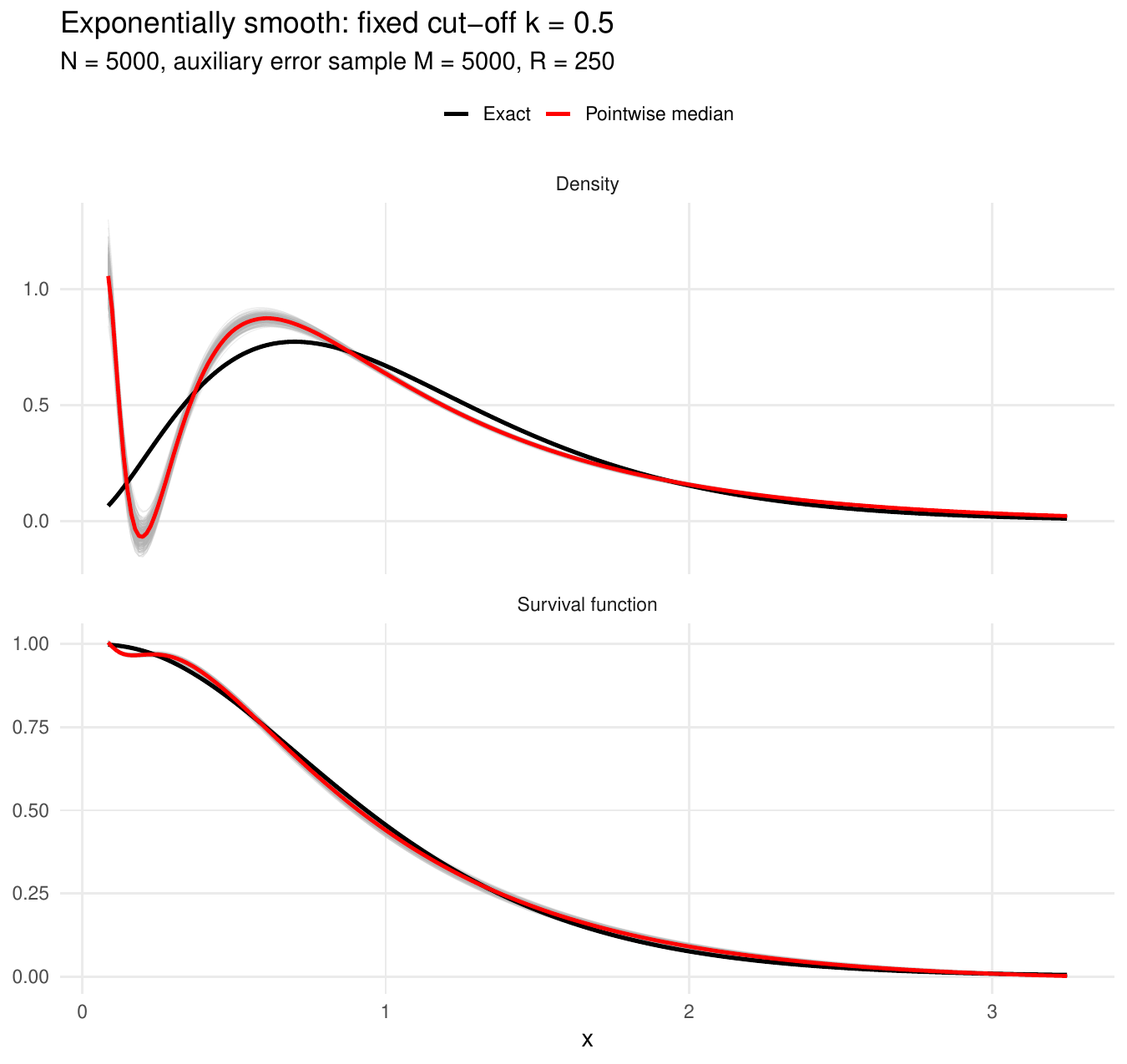}
        \caption{ES, $N=5000$, \(M=5000\).}
    \end{subfigure}
    \caption{Monte Carlo estimates for a fixed cut-off \(k=0.5\). Each panel
    contains the density estimate above the survival-function estimate. Grey
    lines represent the \(R=250\) repetitions, the black line is the true
    target, and the red line is the pointwise Monte Carlo median.}
    \label{fig:fixed-k-by-n}
\end{figure}

\begin{table}[htbp]
    \centering
    \small
    \caption{Pointwise Monte Carlo MSE at \(x_0=0.8\) for
             \(N=5000\) with unknown error distribution.}
    \label{tab:pointwise-mse-unknown}

    \begin{tabular}{llcccc}
        \toprule
        & & \multicolumn{4}{c}{
            \(\widehat{\mathrm{MSE}}\;(\times 10^{-3})\)
        } \\
        \cmidrule(lr){3-6}
        \textbf{Smoothness}
        & \textbf{Target}
        & \(M=100\)
        & \(M=500\)
        & \(M=1000\)
        & \(M=5000\) \\
        \midrule
        OS & Density
           & 0.855 & 0.321 & 0.295 & 0.242 \\
        OS & Survival function
           & 1.97  & 1.75  & 1.724  & 1.719  \\
        \addlinespace
        ES & Density
           & 1.029  & 1.028  & 1.028  & 1.027  \\
        ES & Survival function
           & 0.540 & 0.422 & 0.392 & 0.357 \\
        \bottomrule
    \end{tabular}
\end{table}

\subsubsection{Data-driven cut-off selection}
As before, we simulate the data-driven selection rule
developed in \cref{sec:lifumu:unknown:datadriven} (particuarly \cref{lifumu:unknown:datadriven}), again with $\kappa = 100$.

In every Monte Carlo repetition, \(\widehat k\) is computed from the current
sample.
The experiment is repeated for $N = 5000$ and
\(M\in\{100, 500,1000,5000\}\) and both smoothness regimes.

\begin{figure}[H]
    \centering
    \begin{subfigure}[t]{0.48\textwidth}
        \centering
        \includegraphics[width=\linewidth]{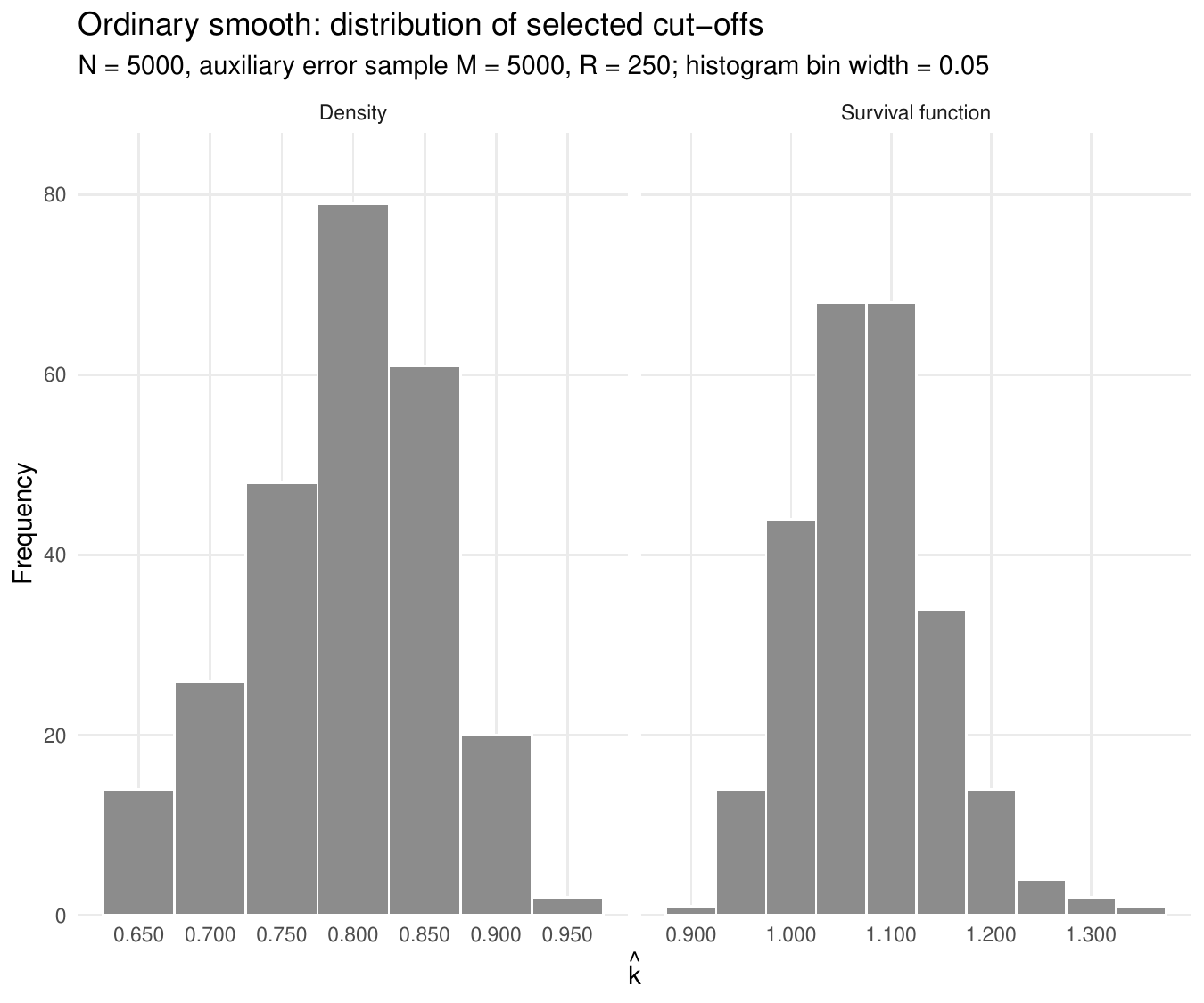}
        \caption{Data-driven cut-off selection distribution: Ordinary Smoothness.}
    \end{subfigure}
    \hfill
    \begin{subfigure}[t]{0.48\textwidth}
        \centering
        \includegraphics[width=\linewidth]{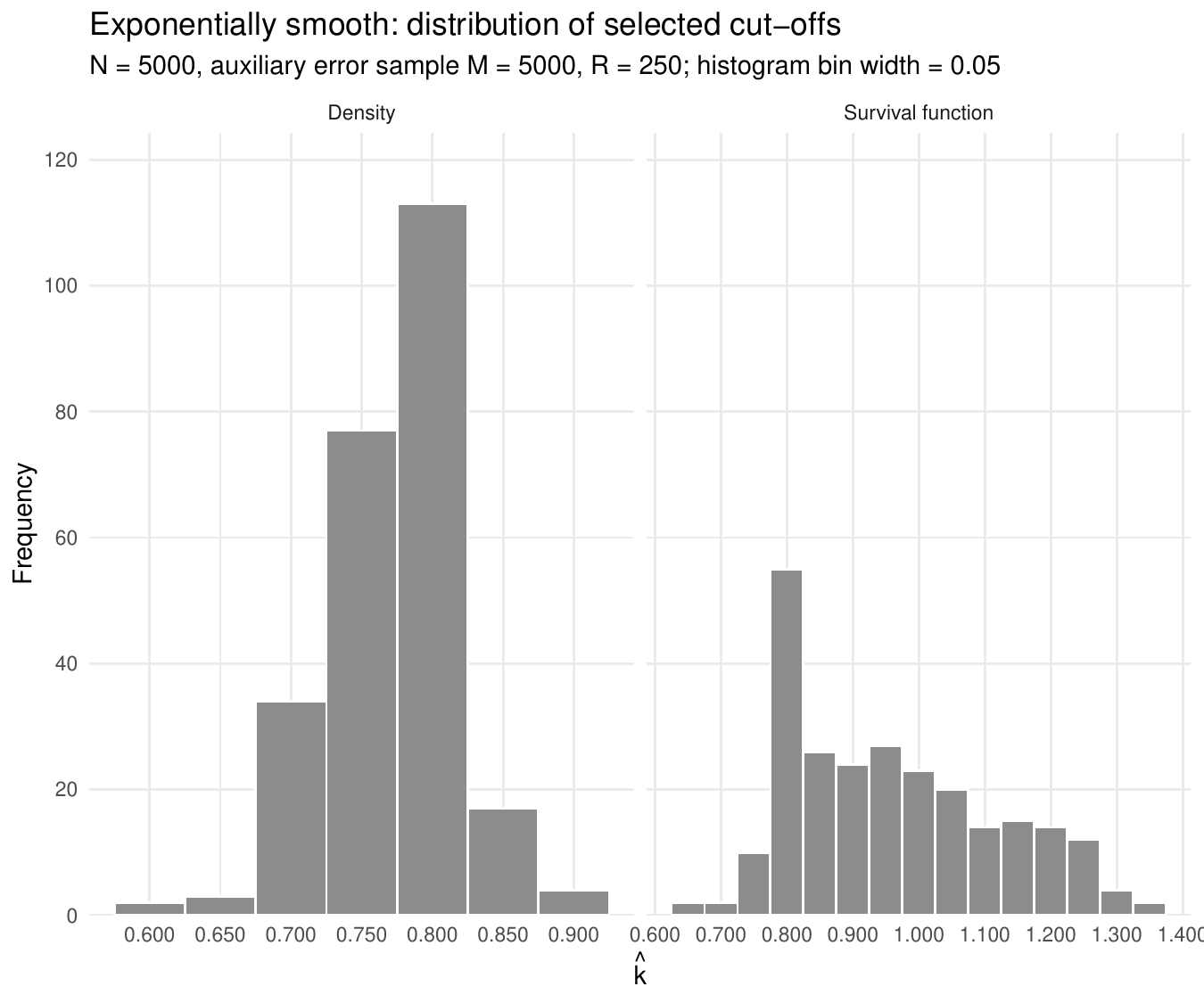}
        \caption{Data-driven cut-off selection distribution: Exponentially Smoothness.}
    \end{subfigure}
    \caption{Data-driven cut-off selection distributions illustrated by corresponding histograms for $N = M = 5000$.}
    \label{fig:data-driven-selection:unknown}
\end{figure}

\begin{table}[htbp]
    \centering
    \small
    \caption{Pointwise Monte Carlo MSE at \(x_0=0.8\) for
             \(N=5000\) with unknown error distribution and
             data-driven cut-off selection.}
    \label{tab:pointwise-mse-unknown-data-driven}

    \begin{tabular}{llcccc}
        \toprule
        & & \multicolumn{4}{c}{
            \(\widehat{\mathrm{MSE}}\;(\times 10^{-3})\)
        } \\
        \cmidrule(lr){3-6}
        \textbf{Smoothness}
        & \textbf{Target}
        & \(M=100\)
        & \(M=500\)
        & \(M=1000\)
        & \(M=5000\) \\
        \midrule
        OS & Density
           & 2.14  & 0.807 & 0.584  & 0.421 \\
        OS & Survival function
           & 0.468 & 0.133 & 0.0999 & 0.0639 \\
        \addlinespace
        ES & Density
           & 0.804 & 0.611 & 0.570  & 0.573 \\
        ES & Survival function
           & 0.272 & 0.114 & 0.0898 & 0.0743 \\
        \bottomrule
    \end{tabular}
\end{table}

\section*{Funding}
\addcontentsline{toc}{section}{Funding}

MS acknowledges partial funding from the Deutsche Forschungsgemeinschaft (DFG, German Research Foundation) through SFB 1294 (project ID 318763901),
as well as under Germany’s Excellence Strategy—EXC 2181/1—390900948 (the Heidelberg STRUCTURES Cluster of Excellence).

\section*{AI Disclosure}
\addcontentsline{toc}{section}{AI Disclosure} 
ChatGPT 6.1 was used for proof-reading and for the simulations of \cref{sec:simulation-experiments}. 
All AI-assisted output was carefully reviewed and verified by the authors, who take full responsibility for the content and results presented in this work.

\printbibliography[heading=bibintoc,title={References}]

@article{AndersenHansen_2001,
  title = {Multiplicative Censoring: Density Estimation by a Series Expansion Approach},
  shorttitle = {Multiplicative Censoring},
  author = {Andersen, Kim E. and Hansen, Martin B.},
  year = 2001,
  month = oct,
  journal = {Journal of Statistical Planning and Inference},
  volume = {98},
  number = {1},
  pages = {137--155},
  issn = {0378-3758},
  doi = {10.1016/S0378-3758(00)00237-8},
  urldate = {2026-04-06}
}

@article{AsgharianWolfson_2005,
  title = {Asymptotic Behavior of the Unconditional {{NPMLE}} of the Length-Biased Survivor Function from Right Censored Prevalent Cohort Data},
  author = {Asgharian, Masoud and Wolfson, David B.},
  year = 2005,
  month = oct,
  journal = {The Annals of Statistics},
  volume = {33},
  number = {5},
  pages = {2109--2131},
  publisher = {Institute of Mathematical Statistics},
  issn = {0090-5364, 2168-8966},
  doi = {10.1214/009053605000000372},
  urldate = {2026-02-06}
}

@article{Barronetal_1999,
  title = {Risk Bounds for Model Selection via Penalization},
  author = {Barron, Andrew and Birg{\'e}, Lucien and Massart, Pascal},
  year = 1999,
  month = feb,
  journal = {Probability Theory and Related Fields},
  volume = {113},
  number = {3},
  pages = {301--413},
  issn = {1432-2064},
  doi = {10.1007/s004400050210},
  urldate = {2026-02-06},
  langid = {english}
}

@article{Belomestnyetal_2016,
  title = {Nonparametric {{Laguerre}} Estimation in the Multiplicative Censoring Model},
  author = {Belomestny, Denis and Comte, Fabienne and {Genon-Catalot}, Valentine},
  year = 2016,
  month = jan,
  journal = {Electronic Journal of Statistics},
  volume = {10},
  number = {2},
  pages = {3114--3152},
  publisher = {{Institute of Mathematical Statistics and Bernoulli Society}},
  issn = {1935-7524, 1935-7524},
  doi = {10.1214/16-EJS1203},
  urldate = {2026-02-06}
}

@article{BelomestnyGoldenshluger_2020,
  title = {Nonparametric Density Estimation from Observations with Multiplicative Measurement Errors},
  author = {Belomestny, Denis and Goldenshluger, Alexander},
  year = 2020,
  month = feb,
  journal = {Annales de l'Institut Henri Poincar\'e, Probabilit\'es et Statistiques},
  volume = {56},
  number = {1},
  pages = {36--67},
  publisher = {Institut Henri Poincar\'e},
  issn = {0246-0203},
  doi = {10.1214/18-AIHP954},
  urldate = {2026-02-06}
}

@article{BirgeMassart_1998,
  title = {Minimum Contrast Estimators on Sieves: Exponential Bounds and Rates of Convergence},
  shorttitle = {Minimum Contrast Estimators on Sieves},
  author = {Birg{\'e}, Lucien and Massart, Pascal},
  year = 1998,
  month = sep,
  journal = {Bernoulli},
  volume = {4},
  number = {3},
  pages = {329--375},
  publisher = {{Bernoulli Society for Mathematical Statistics and Probability}},
  issn = {1350-7265},
  urldate = {2026-02-19},
  langid = {english}
}

@article{BrennerMiguel_2022,
  title = {Anisotropic Spectral Cut-off Estimation under Multiplicative Measurement Errors},
  author = {Brenner Miguel, Sergio Filipe},
  year = 2022,
  month = jul,
  journal = {Journal of Multivariate Analysis},
  volume = {190},
  pages = {104990},
  issn = {0047-259X},
  doi = {10.1016/j.jmva.2022.104990},
  urldate = {2026-02-06}
}

@book{BrennerMiguel_2023,
  title = {The {{Mellin}} Transform in Nonparametric Statistics},
  author = {Brenner Miguel, Sergio Filipe},
  year = 2023,
  publisher = {Dissertation, Ruprecht-Karls-Universit\"at Heidelberg, 2023},
  address = {Heidelberg},
  doi = {10.11588/heidok.00033855},
  urldate = {2025-11-09},
  langid = {english}
}

@article{BrennerMiguelPhandoidaen_2022,
  title = {Multiplicative Deconvolution in Survival Analysis under Dependency},
  author = {Brenner Miguel, Sergio Filipe and Phandoidaen, Nathawut},
  year = 2022,
  month = mar,
  journal = {Statistics},
  volume = {56},
  number = {2},
  pages = {297--328},
  publisher = {Taylor \& Francis},
  issn = {0233-1888},
  doi = {10.1080/02331888.2022.2058503},
  urldate = {2026-02-06}
}

@article{Bruneletal_2016,
  title = {Nonparametric Density and Survival Function Estimation in the Multiplicative Censoring Model},
  author = {Brunel, Elodie and Comte, Fabienne and {Genon-Catalot}, Valentine},
  year = 2016,
  month = sep,
  journal = {TEST},
  volume = {25},
  number = {3},
  pages = {570--590},
  issn = {1863-8260},
  doi = {10.1007/s11749-016-0479-1},
  urldate = {2026-02-06},
  langid = {english}
}

@article{ComteDion_2016,
  title = {Nonparametric Estimation in a Multiplicative Censoring Model with Symmetric Noise},
  author = {Comte, F. and Dion, C.},
  year = 2016,
  month = oct,
  journal = {Journal of Nonparametric Statistics},
  volume = {28},
  number = {4},
  pages = {768--801},
  publisher = {Taylor \& Francis},
  issn = {1048-5252},
  doi = {10.1080/10485252.2016.1225737},
  urldate = {2026-02-06}
}

@article{ComteLacour_2011,
  title = {Data-{{Driven Density Estimation}} in the {{Presence}} of {{Additive Noise}} with Unknown {{Distribution}}},
  author = {Comte, F. and Lacour, C.},
  year = 2011,
  month = sep,
  journal = {Journal of the Royal Statistical Society Series B: Statistical Methodology},
  volume = {73},
  number = {4},
  pages = {601--627},
  issn = {1369-7412},
  doi = {10.1111/j.1467-9868.2011.00775.x},
  urldate = {2026-02-06}
}

@article{Johannes_2009,
  title = {Deconvolution with Unknown Error Distribution},
  author = {Johannes, Jan},
  year = 2009,
  month = oct,
  journal = {The Annals of Statistics},
  volume = {37},
  number = {5A},
  pages = {2301--2323},
  publisher = {Institute of Mathematical Statistics},
  issn = {0090-5364, 2168-8966},
  doi = {10.1214/08-AOS652},
  urldate = {2026-02-06}
}

@article{JohannesSchwarz_2013,
  title = {Adaptive Circular Deconvolution by Model Selection under Unknown Error Distribution},
  author = {Johannes, Jan and Schwarz, Maik},
  year = 2013,
  month = nov,
  journal = {Bernoulli},
  volume = {19},
  number = {5A},
  pages = {1576--1611},
  publisher = {{Bernoulli Society for Mathematical Statistics and Probability}},
  issn = {1350-7265},
  doi = {10.3150/12-BEJ422},
  urldate = {2026-02-06}
}

@book{Massart_2007,
  title = {Concentration {{Inequalities}} and {{Model Selection}}},
  author = {Massart, Pascal},
  editor = {Picard, Jean},
  year = 2007,
  series = {Lecture {{Notes}} in {{Mathematics}}},
  volume = {1896},
  publisher = {Springer},
  address = {Berlin, Heidelberg},
  doi = {10.1007/978-3-540-48503-2},
  urldate = {2026-02-06},
  copyright = {http://www.springer.com/tdm},
  isbn = {978-3-540-48497-4},
  langid = {english}
}

@book{Meister_2009,
  title = {Deconvolution {{Problems}} in {{Nonparametric Statistics}}},
  author = {Meister, Alexander},
  year = 2009,
  series = {Lecture {{Notes}} in {{Statistics}}},
  volume = {193},
  publisher = {Springer},
  address = {Berlin, Heidelberg},
  doi = {10.1007/978-3-540-87557-4},
  urldate = {2026-02-03},
  copyright = {http://www.springer.com/tdm},
  isbn = {978-3-540-87556-7 978-3-540-87557-4}
}

@article{Migueletal_2021,
  title = {Spectral Cut-off Regularisation for Density Estimation under Multiplicative Measurement Errors},
  author = {Brenner Miguel, Sergio Filipe and Comte, Fabienne and Johannes, Jan},
  year = 2021,
  month = jan,
  journal = {Electronic Journal of Statistics},
  volume = {15},
  number = {1},
  pages = {3551--3573},
  publisher = {{Institute of Mathematical Statistics and Bernoulli Society}},
  issn = {1935-7524, 1935-7524},
  doi = {10.1214/21-EJS1870},
  urldate = {2026-02-06}
}

@article{Migueletal_2023,
  title = {Linear Functional Estimation under Multiplicative Measurement Error},
  author = {Brenner Miguel, Sergio Filipe and Comte, Fabienne and Johannes, Jan},
  year = 2023,
  month = aug,
  journal = {Bernoulli},
  volume = {29},
  number = {3},
  pages = {2247--2271},
  publisher = {{Bernoulli Society for Mathematical Statistics and Probability}},
  issn = {1350-7265},
  doi = {10.3150/22-BEJ1540},
  urldate = {2025-11-20}
}

@article{Migueletal_2024,
  title = {Multiplicative Deconvolution under Unknown Error Distribution},
  author = {Brenner Miguel, Sergio Filipe and Johannes, Jan and Siebel, Maximilian},
  year = 2024,
  month = jan,
  journal = {Electronic Journal of Statistics},
  volume = {18},
  number = {2},
  pages = {4795--4850},
  publisher = {{Institute of Mathematical Statistics and Bernoulli Society}},
  issn = {1935-7524, 1935-7524},
  doi = {10.1214/24-EJS2314},
  urldate = {2025-07-11}
}

@misc{Neubert_2025,
  type = {Dissertation},
  title = {Nonparametric Inference in Convolution Models},
  author = {Neubert, Bianca Marie},
  year = 2025,
  address = {Heidelberg},
  doi = {10.11588/heidok.00037748},
  urldate = {2026-02-05},
  copyright = {info:eu-repo/semantics/openAccess},
  langid = {english}
}

@article{NeumannHossjer_1997,
  title = {On the Effect of Estimating the Error Density in Nonparametric Deconvolution},
  author = {Neumann, Michael H.},
  year = 1997,
  month = jan,
  journal = {Journal of Nonparametric Statistics},
  volume = {7},
  number = {4},
  pages = {307--330},
  publisher = {Taylor \& Francis},
  issn = {1048-5252},
  doi = {10.1080/10485259708832708},
  urldate = {2026-02-06}
}

@book{Petrov_1995,
  title = {Oxford {{Studies In Probability}} 4: {{Limit Theorems}} of {{Probability Theory Sequences}} of {{Independent Random Variables}}},
  shorttitle = {Oxford {{Studies In Probability}} 4},
  author = {Petrov, Valentin V},
  year = 1995,
  month = apr,
  publisher = {Oxford University Press},
  doi = {10.1093/oso/9780198534990.001.0001},
  urldate = {2026-02-28},
  isbn = {978-0-19-853499-0}
}

@article{vanEsetal_2000,
  title = {Survival Analysis under Cross-Sectional Sampling: Length Bias and Multiplicative Censoring},
  shorttitle = {Survival Analysis under Cross-Sectional Sampling},
  author = {{van Es}, Bert and Klaassen, Chris A. J. and Oudshoorn, Karin},
  year = 2000,
  month = dec,
  journal = {Journal of Statistical Planning and Inference},
  series = {Perspectives in Modern Statistical Inference: Parametric Semiparametrics, Nonparametrics},
  volume = {91},
  number = {2},
  pages = {295--312},
  issn = {0378-3758},
  doi = {10.1016/S0378-3758(00)00183-X},
  urldate = {2026-02-04}
}

@article{Vardi_1989,
  title = {Multiplicative Censoring, Renewal Processes, Deconvolution and Decreasing Density: {{Nonparametric}} Estimation},
  shorttitle = {Multiplicative Censoring, Renewal Processes, Deconvolution and Decreasing Density},
  author = {Vardi, Y.},
  year = 1989,
  month = dec,
  journal = {Biometrika},
  volume = {76},
  number = {4},
  pages = {751--761},
  issn = {0006-3444},
  doi = {10.1093/biomet/76.4.751},
  urldate = {2026-02-04}
}

@article{VardiZhang_1992,
  title = {Large {{Sample Study}} of {{Empirical Distributions}} in a {{Random-Multiplicative Censoring Model}}},
  author = {Vardi, Y. and Zhang, Cun-Hui},
  year = 1992,
  month = jun,
  journal = {The Annals of Statistics},
  volume = {20},
  number = {2},
  pages = {1022--1039},
  publisher = {Institute of Mathematical Statistics},
  issn = {0090-5364, 2168-8966},
  doi = {10.1214/aos/1176348668},
  urldate = {2026-02-04}
}

@article{GoldenshlugerLepski_2011,
  title = {Bandwidth Selection in Kernel Density Estimation: {{Oracle}} Inequalities and Adaptive Minimax Optimality},
  shorttitle = {Bandwidth Selection in Kernel Density Estimation},
  author = {Goldenshluger, Alexander and Lepski, Oleg},
  year = 2011,
  month = jun,
  journal = {The Annals of Statistics},
  volume = {39},
  number = {3},
  pages = {1608--1632},
  publisher = {Institute of Mathematical Statistics},
  issn = {0090-5364, 2168-8966},
  doi = {10.1214/11-AOS883},
  urldate = {2026-11-06},
  langid = {english}
}

\newpage
\appendix

\section{Proofs of \cref{sec:lifumu:known}}

\subsection{\texorpdfstring{Proofs of \cref{sec:lifumu:oracle:known}}{Proofs of Section 3.1}}\label{sec:app:lifumu:known:oracle}

\begin{proof}[Proof of \cref{example:lifumu:known:oracle:rates}] 
	As we study the asymptotic regime $N\to\infty$, we assume without loss of generality that $\kdiamond \asymp \kdiamond(N)\in\IR_{\geq 1}$.
		\begin{itemize}
			\item  Proof of (\textsc{os}): For $k\in\IR_{\geq 1}$ and $f$ being ordinary smooth, we obtain
					\begin{align*}
						\norm{\rmL\indkco{k}\rmsos^\dagger}_{\IL^2}^2&
							 = \int_\IR \abs{\rmL(t)}^2\rmsos^\dagger(t)^2\indkco{k}(t)\rmd\lebesgue(t)\\
							&\leq \ConstLupper^2\times\int_\IR \br{1+t^2}^{-(s-p)}\indkco{k}(t)\rmd\lebesgue(t)\\
							&\leq 2\ConstLupper^2\times\int_k^\infty t^{-2(s-p)} \rmd\lebesgue(t)
							 =\const(\rmL,s,p) \times k^{-2(s-p)+1}.
					\end{align*}
				
			\item Proof of (\textsc{es}): 
			For $k\in\IR_{\geq 1}$ and $f$ being exponentially smooth, we obtain analogously
				\begin{align*}
 					 \norm{\rmL\indkco{k}\rmses^\dagger}_{\IL^2}^2
					&\leq 2\ConstLupper^2\times\int_k^{\infty} \br{1+t^2}^{p}\exp(-t^{2s})\rmd\lebesgue(t)\\
					& = 2\ConstLupper^2\times\int_k^{\infty} \br{1+t^2}^{p} t^{2s-1}t^{-2s+1}\exp(-t^{2s})\rmd\lebesgue(t)\\
					& \leq \const(\rmL,p)\times \int_k^{\infty}  t^{2s-1}t^{-2s+2p+1}\exp(-t^{2s})\rmd\lebesgue(t)\\
					& \leq \const(\rmL,p)\times \int_k^{\infty}  t^{2s-1}\exp(-t^{2s})\rmd\lebesgue(t)\\
					& =  \const(\rmL,p,s)\times \exp\br{-k^{2s}}.
				\end{align*}
		\end{itemize}
		After deriving the order of the bias term, we focus on the cases $(1)$-$(7)$ separately by firstly deriving the order of the variance term and lastly balancing both terms.
		\begin{itemize}
			\item (1):
				Let $p+b>-\frac{1}{2}$. We then have for $k\in\pRz$
				\begin{align*}
					\norm{\rmL\indk{k}\rmbos}_{\IL^2}^2
					 	 &= \int_{[-k,k]}\abs{\rmL(t)}^2\rmbos(t)^2\rmd\lebesgue(t)\\
						& \leq \const(\rmL,p,b)\int_0^k t^{2(p+b)}\rmd t
						\leq \const(\rmL,p,b)\times k^{2(p+b)+1}.
				\end{align*}
				Thus, by balancing the bias and the variance term, we obtain 
				\begin{align*}
					&\frac{1}{N} k^{2(p+b)+1} \asymp k^{-2(s-p)+1}\quad
					\Longleftrightarrow\quad \frac{1}{N} \asymp k^{-2(s+b)}\quad
					\Longleftrightarrow \quad N^{\frac{1}{2(s+b)}}\asymp k.
				\end{align*}
				Computing $\maxrisklocal \br{\kdiamond|\rmsos,\rmbos,\rmL}$ with $\kdiamond \asymp N^{\frac{1}{2(s+b)}}$ yields the claim.

			\item (2): 
				Let $p+b=-\frac{1}{2}$. Similarly, we have for $k$ sufficiently large
					\begin{equation*}
						\norm{\rmL\indk{k}\rmbos}_{\IL^2}^2 = \int_{-k}^k \frac{1}{\sqrt{1+t^2}}\rmd t = 2\arsinh(k) \leq \const \times \ln(k).
					\end{equation*}
					Thus minimizing the function
						\begin{equation*}
							\pRz\ni k\mapsto \frac{1}{N}\ln(k) +k^{-2(s-p)+1} \in\pRz
						\end{equation*}
					leads to an extremum 
					\begin{equation*}
						\kdiamond \asymp N^{\frac{1}{2(s-p)-1}}.
					\end{equation*}
					Computing $\maxrisklocal \br{\kdiamond|\rmsos,\rmbos,\rmL}$ with $\kdiamond \asymp N^{\frac{1}{2(s-p)-1}}$ yields the claim.
			\item (3):
				 Let $p+b<-\frac{1}{2}$. We then have for $k\in\pRz$
					\begin{equation*}
                        \norm{\rmL\indk{k}\rmbos}_{\IL^2}^2\leq \const(\rmL)\int_{-k}^k \frac{1}{(1+t^2)^{-(p+b)}}\rmd t \leq \const(\rmL) \int_{-\infty}^\infty \frac{1}{(1+t^2)^{-(p+b)}} \rmd t < \infty
					\end{equation*}
					
					Thus by balancing both terms, we obtain 
					\begin{equation*}
						\kdiamond \asymp N^\frac{1}{2(s-p)-1}.
					\end{equation*}
					Computing $\maxrisklocal \br{\kdiamond|\rmsos,\rmbos,\rmL}$ with $\kdiamond \asymp N^{\frac{1}{2(s-p)-1}}$ yields the claim.

			\item (4): For $p\in\IR$, and $b\in\pRzz$, we have 
				\begin{align*}
					\norm{\rmL\indk{k}\rmbes}_{\IL^2}^2 &\leq \const(\rmL)\int_{0}^k \br{1+t^2}^p\exp\br{t^{2b}}\rmd t\\
						& \leq \const(\rmL,p)\times\br{1+\int_{1}^k t^{2p+2b}\exp\br{t^{2b}}t^{-2b}\rmd t}\\
						& \leq  \const(\rmL,p)\times\br{1+k^{(2p+2b+1)_+}\exp\br{k^{2b}}}\\
						& \leq  \const(\rmL,p)\times k^{(2p+2b+1)_+}\exp\br{k^{2b}}.
				\end{align*}
				Thus, balancing both terms, we see that the exponential in the second term is the dominant one.
				We obtain for large $N$
				\begin{equation*}
					\exp\br{k^{2b}} \asymp N\times k^{-2(s-p)+1} \quad \Longleftrightarrow \quad k^{2b} \asymp \ln(N) - \ln(k)
				\end{equation*}
				and thus 
				\begin{equation*}
					\kdiamond \asymp \ln(N)^{\frac{1}{2b}}.
				\end{equation*}
				Computing $\maxrisklocal \br{\kdiamond|\rmsos,\rmbes,\rmL}$ with $\kdiamond \asymp \ln(N)^{\frac{1}{2b}}$ yields the claim.
			
			\item $(5)$ - $(7)$: Similarly as derived in $(4)$, we see that 
				the exponential in the bias term is the dominant expression. The optimal choice is thus given by 
				$$\kdiamond \asymp \ln(N)^{\frac{1}{2s}}.$$ Computing $\maxrisklocal \br{\kdiamond|\rmses,\rmbos,\rmL}$ with this choice then yields the claims.
		\end{itemize}
\end{proof}

\subsection{\texorpdfstring{Proofs of \cref{sec:lifumu:datadriven:known}}{Proofs of Section 3.2}}\label{sec:lifumu:appendix:datadriven:known}

In this section we derive general results for the data-driven selection methods, which not only provide the 
proofs of the statements formulated in \cref{sec:lifumu:datadriven:known}, but also in \cref{sec:lifumu:unknown:datadriven}.
To that end, let $\rmw\colon\IR\to\pRzz$ be a weight function, such that $\rmw\indk{k}\in\IL^2$ for all $k\in\pRz$. 
Let $(Y_i)_{i\in\nset{N}}$, $N\in\IN$, be i.i.d. copies of $Y\sim\IP^Y$, a strictly positive real-valued random variable, such that
$f^Y\in\IL^1_{\rmx^{c-1}}$ for some $c\in\pRz$.
We define $\phi\colon\pRz\times\IR\to\IC$ by $\phi(y,t)\coloneqq y^{c-1+2\pi\iota t}$, noting that $\phi$ is jointly 
$\scrB_{\pRz}\otimes\scrB_{\IR}-\scrB_{\IC}$-measurable.
For every fixed $y\in\pRz$ we set $\phi_y\coloneqq \phi(y,\cdot)$.
Evidently, $\phi$ does not need to be a bounded function, and therefore we decompose $\phi = \phib + \phiu$ into a bounded part $\phib\colon\pRz\times\IR\to\IC$
and in an unbounded part $\phiu\colon\pRz\times\IR\to\IC$ by defining    
\begin{equation*}
    \phi(y,t) = \phib(y,t) + \phiu(y,t) \coloneqq \1_{(0,\nu)}(y^{c-1}) y^{c-1+2\pi\iota t} +\1_{[\nu,\infty)}(y^{c-1}) y^{c-1+2\pi\iota t}
\end{equation*}
for some $\nu\in\pRz$.
We define for each $k\in\IR_{>0}$ the empirical process
\begin{align*}
    \rmZ_k^{\rmw}&\coloneqq \frac{1}{N}\sum_{i\in\nset{N}}\int_{\IR}\br{\phi(Y_i,t)-\IE\brE{\phi(Y_i,t)}}\rmw(t)\indk{k}(t)\rmd\lebesgue(t)\\
    & = \int_{\IR}\br{\MelYhat(t)-\MelY(t)}\rmw(t)\indk{k}(t)\rmd\lebesgue(t).
\end{align*}
For
	$\mathrm{j}\in\{\mathrm{b},\mathrm{u}\}$ we set
	$\MelYhat^{\mathrm{j}}(t)=\tfrac{1}{N}\sum_{i\in\nset{N}}\phi^{\mathrm{j}}(Y_1,t)$ and
	$\MelY^{\mathrm{j}}(t)=\E\big[\phi^{\mathrm{j}}(Y_1,t)\big]$ for
	each $t\in\IR$, respectively, such that $\MelYhat-\MelY=\MelYhat^{\mathrm{b}}-\MelY^{\mathrm{b}}+\MelYhat^{\mathrm{u}}-\MelY^{\mathrm{u}}$.
Thus, we can set
\begin{align*}
    \IWb_k&\coloneqq 
     \int_{\IR}\br{\MelYhatb(t)-\MelYb(t)}\rmw(t)\indk{k}(t)\rmd\lebesgue(t)
\end{align*}
as well as
\begin{align*}
    \IWu_k&\coloneqq
     \int_{\IR}\br{\MelYhatu(t)-\MelYu(t)}\rmw(t)\indk{k}(t)\rmd\lebesgue(t),
\end{align*}
and observe $\rmZ_k^{\rmw} = \IWb_k + \IWu_k$.\\

\begin{lemma}[Variance of $\IWu_k$]\label{lifumu:app:varinace:process}
    \mbox{}\\
    Assume, there exists $\rmp\in\IR_{\geq 2}$, such that $\IE\brE{Y^{2\rmp(c-1)}}<\infty$.
    For any $N\in\IN$ and $k\in\IR_{>0}$, we have
    \begin{equation*}
        \IE_Y^N\brE{\abs{\IWu_k}^2} \leq \frac{2k}{N}\norm{\indk{k}\rmw}_{\IL^2}^2\nu^{2(1-\rmp)}\IE\brE{Y^{2\rmp(c-1)}}.
    \end{equation*}
\end{lemma}

\begin{proof}[Proof of \cref{lifumu:app:varinace:process}]
    Observing that $\IWu_k$ defines a centered random variable, we have 
    \begin{align*}
        N\times\Var_Y^N\br{\IWu_k} &= \IE\brE{\abs{\int_{\IR}\br{\phiu(Y,t)-\IE\brE{\phiu(Y,t)}}\rmw(t)\indk{k}(t)\rmd\lebesgue(t)}^2}\\
        & \leq \IE\brE{\abs{\int_{\IR}{\phiu(Y,t)\rmw(t)\indk{k}(t)\rmd\lebesgue(t)}}^2}\\
        & \leq \norm{\rmw\indk{k}}_{\IL^2}^2\cdot \IE\brE{\norm{\phiu(Y,\cdot)\indk{k}}_{\IL^2}^2},
    \end{align*}
    where the last line follows from the Cauchy-Schwarz inequality.
    Evidently, for all $t\in\IR$, we have 
    \begin{align*}
      \abs{\phiu(Y,t)}^2 &= Y^{2(c-1)}\1_{[\nu,\infty)}(Y^{c-1}) = Y^{2\rmp(c-1)-2(\rmp-1)(c-1)}\1_{[\nu,\infty)}(Y^{c-1}) \leq Y^{2\rmp(c-1)}\nu^{-2(\rmp-1)},
    \end{align*}
    and thus
    \begin{equation*}
        \IE\brE{\norm{\phiu(Y,\cdot)\indk{k}}_{\IL^2}^2} \leq 2k\nu^{-2(\rmp-1)}\IE\brE{Y^{2\rmp(c-1)}}.
    \end{equation*}
    This shows the claim.\\
\end{proof}

{
  \renewcommand{\rmp}{\mathrm{p}}
  \renewcommand{\penbar}{\mathfrak{Pen}}
    \begin{lemma}[Concentration of $\IWb_k$]\label{lem:lifumu:app:concentration:process:bounded}
        \mbox{}\\
        Let $k\in\IR_{>0}$ be fixed. Let
        \begin{equation*}
            \rmV_k^{\rmw}\coloneqq \IE\brE{\abs{\int_\IR\phib(Y_1,t)\rmw(t)\indk{k}(t)\rmd\lebesgue(t)}^2}.
        \end{equation*}
        For any $\rmp,\Delta_k\in\IR_{\geq 1}$ set
        \begin{equation*}
            \penbar_k^{\rmw}\coloneqq \frac{\ln(\Delta_k)}{N} \br{4\rmp\rmV_k^{\rmw} + \frac{512k\rmp^2}{N}\norm{\rmw\indk{k}}_{\IL^2}^2\nu^2\ln(\Delta_k)}.
        \end{equation*}
        We then have 
        \begin{equation*}
            \forall k\in\IR_{>0}:\;\IE_{Y}^N\brE{\br{\abs{\IWb_k}^2-\penbar_k^{\rmw}}_+} \leq \frac{8\Delta_k^{-\rmp}}{N}\norm{\rmw\indk{k}}_{\IL^2}^2\br{\ConstfY\Constfeps+\frac{256k\nu^2}{N}}
        \end{equation*}
    \end{lemma}

\begin{proof}[Proof of \cref{lem:lifumu:app:concentration:process:bounded}]
    We aim to apply Bernstein's inequality (precisely, \cref{eq:Bern:4} of \cref{lem:BernsteinInequality_general}) on $\IWb_k$. 
    As $\IWb_k$ defines a $\IC$-valued random variable, we consider its real and imaginary part separately.
    To that end, we define the family $(\Phib(Y_i))_{i\in\nset{N}}$ of mean-zero random variables by
    \begin{equation*}
        \forall i\in\nset{N}:\;\Phib(Y_i) \coloneqq \int_\IR \br{\phib(Y_i,t)-\IE[\phib(Y_i,t)]}\rmw(t)\indk{k}(t)\rmd\lebesgue(t).
    \end{equation*}
    Note, $\IWb_k = \frac{1}{N}\sum_{i\in\nset{N}} \Phib (Y_i)$.
    We show that $\br{\Re{\Phib(Y_i)}}_{i\in\nset{N}}$ and $\br{\Im{\Phib(Y_i)}}_{i\in\nset{N}}$ both satisfy \cref{eq:condition:general:Bernstein}.
    Indeed, by definition, we have $\abs{\phib(t,y)}\leq\nu$ for all $(t,y)\in\IR\times\pRz$, and thus
    \begin{align*}
        \forall i\in\nset{N}:\;\abs{\Phib(Y_i)} 
         \leq2 \nu\norm{\rmw\indk{k}}_{\IL^2}\norm{\indk{k}}_{\IL^2}
            \leq 2\sqrt{2k}\nu\norm{\rmw\indk{k}}_{\IL^2}\eqqcolon \rmb.
    \end{align*}
    In particular, we have 
    \begin{equation*}
        \forall i\in\nset{N}:\; \abs{\Re{\Phib(Y_i)}},\abs{\Im{\Phib(Y_i)}} \leq \rmb.
    \end{equation*}
    Further, for all $i\in\nset{N}$
    \begin{align*}
        \IE\brE{\Re{\Phib(Y_i)}^2} &= \IE\brE{\abs{ \int_\IR \Re{\phib(Y_i,t)-\IE[\phib(Y_i,t)]}\rmw(t)\indk{k}(t)\rmd\lebesgue(t)}^2}\\
            & \leq \IE\brE{\abs{ \int_\IR \Re{\phib(Y_i,t)}\rmw(t)\indk{k}(t)\rmd\lebesgue(t)}^2} \eqqcolon \rmv_{re}^2,
    \end{align*}
    and analogously
    \begin{align*}
        \IE\brE{\Im{\Phib(Y_i)}^2} & \leq \IE\brE{\abs{ \int_\IR \Im{\phib(Y_i,t)}\rmw(t)\indk{k}(t)\rmd\lebesgue(t)}^2} \eqqcolon \rmv_{im}^2.
    \end{align*}
    We note that 
    \begin{equation*}
        \rmv_{re}^2+\rmv_{im}^2 =  \IE\brE{\abs{ \int_\IR \phi(Y_i,t)\rmw(t)\indk{k}(t)\rmd\lebesgue(t)}^2} = \rmV_k^{\rmw} \eqqcolon \rmv^2 \leq\Constfeps\ConstfY \norm{\rmw\indk{k}}_{\IL^2}^2,
    \end{equation*}
    where in the last step we have applied \cref{lem:MellinTrick}.
    We conclude, for all $m\in\IN_{\geq 2}$
    \begin{align*}
        &\frac{1}{N}\sum_{i\in\nset{N}}\IE\brE{\abs{\Re{\Phib(Y_i)}}^m}  \leq \frac{m!}{2}\rmv_{re}^2 \rmb^{m-2}\\
        \quad\and\quad&\frac{1}{N}\sum_{i\in\nset{N}}\IE\brE{\abs{\Im{\Phib(Y_i)}}^m} \leq \frac{m!}{2}\rmv_{im}^2 \rmb^{m-2}.
    \end{align*}
    Thus, by an application \cref{lem:BernsteinInequality_general} \cref{eq:Bern:4} we derive for any $K\in\IR_{\geq 1}$
    \begin{align*}
        \IE_{Y}^N\brE{\br{N\abs{\Re{\IWb_k}}^2-\br{4\rmv_{re}^2+32\rmb^2\frac{\ln(K)}{N}}\ln(K)}_+}  \leq\frac{8}{K}\br{\rmv_{re}^2+\frac{16\rmb^2}{N}},
    \end{align*}
    and an analogous bound for $\Im{\IWb_k}$. We deduce, 
    \begin{align}
        &\IE_{Y}^N\brE{\br{N\abs{\IWb_k}^2-\br{4\rmv^2+64\rmb^2\frac{\ln(K)}{N}}\ln(K)}_+}\nonumber\\
            &\hspace{1cm}\leq \IE_{Y}^N\brE{\br{N\abs{\Re{\IWb_k}}^2-\br{4\rmv_{re}^2+32\rmb^2\frac{\ln(K)}{N}}\ln(K)}_+}\nonumber\\
                &\hspace{2cm}+ \IE_{Y}^N\brE{\br{N\abs{\Im{\IWb_k}}^2-\br{4\rmv_{im}^2+32\rmb^2\frac{\ln(K)}{N}}\ln(K)}_+}\nonumber\\
            &\hspace{1cm} \leq\frac{8}{K}\br{\rmv_{re}^2+\frac{16\rmb^2}{N}} +\frac{8}{K}\br{\rmv_{im}^2+\frac{16\rmb^2}{N}}\nonumber\\
            &\hspace{1cm} \leq \frac{8}{K}\br{\rmv^2+\frac{32\rmb^2}{N}}.\label{eq:proof:lifumu:concentration:bound}
    \end{align}
    Choosing $K=\Delta_k^\rmp\in\IR_{\geq 1}$ with $\Delta_k\in\IR_{\geq 1}$ and $\rmp\in\IR_{\geq 2}$ as in the assumption,
    we obtain 
    \begin{align*}
        \br{4\rmv^2+64\rmb^2\frac{\ln(K)}{N}}\ln(K) &= \br{4\rmV_k+64\times8k\nu^2\norm{\rmw\indk{k}}_{\IL^2}^2\frac{\ln(\Delta_k^\rmp)}{N}} \ln(\Delta_k^\rmp)\\
        & = \br{4\rmp\rmV_k+512k\nu^2\rmp^2\norm{\rmw\indk{k}}_{\IL^2}^2\frac{\ln(\Delta_k)}{N}} \ln(\Delta_k)\\
        & = N\times\penbar_k^{\rmw}.
    \end{align*}
    Further
    \begin{align*}
        \frac{8}{K}\br{\rmv^2+\frac{32\rmb^2}{N}} 
        &\leq \frac{8}{\Delta_k^\rmp}\br{\Constfeps\ConstfY\norm{\rmw\indk{k}}_{\IL^2}^2+\frac{256k\nu^2\norm{\rmw\indk{k}}_{\IL^2}^2 }{N}}\\
        & = \frac{8}{\Delta_k^\rmp}\norm{\rmw\indk{k}}_{\IL^2}^2\br{\Constfeps\ConstfY+\frac{256k\nu^2 }{N}}.
    \end{align*}
    Plugging in these computations into \cref{eq:proof:lifumu:concentration:bound} and dividing the whole inequality by $N$ yields the claim.
\end{proof}}

\begin{proposition}[Concentration of $\rmZ_k^{\rmw}$]\label{propo:lifumu:known:concentration:process:app}
    \mbox{}\\
    Assume that there exists $\rmp\in\IR_{\geq 2}$, such that $\IE\brE{Y^{2\rmp(c-1)}}<\infty$.
    Let $\magenta{\kappa}\in\pRz$. For all $k\in\IR_{>\magenta{\kappa}}$ define
    \begin{equation*}
        \Delta_k\coloneqq  \Delta_k^{\rmw}\coloneqq \magenta{2\br{\frac{k}{\kappa}}^2\br{1\lor k\norm{\rmw\indk{k}}_{\IL^2}^2}} \quad\text{and}\quad \nu^2 \coloneqq \nu_k^2 \coloneqq (\nu_k^{\rmw})^2 \coloneqq \frac{1}{128}\Delta_k^{\frac{1}{\rmp-1}},
    \end{equation*}
    and further 
    \begin{equation*}
        \penbar_k^{\rmw}\coloneqq 8\br{\rmV_k^{\rmw}+2k\norm{\rmw\indk{k}}_{\IL^2}^2\frac{\Delta_k^{\frac{1}{\rmp-1}}\ln(\Delta_k)}{N}} \frac{\ln(\Delta_k)}{N}.
    \end{equation*}
    Then for all $\magenta{\calK\subseteq\kappa\IN}$, we have 
    \begin{equation*}
        \IE_{Y}^N\brE{\max_{k\in\magenta{\calK}}\brK{\br{\abs{\rmZ_k^{\rmw}}^2-2\penbar_k^{\rmw}}_+}} \leq \frac{1}{N\magenta{\kappa^2}}\times \br{40\ConstfY\Constfeps+\frac{128^\rmp}{32}\IE\brE{Y^{2\rmp(c-1)}}}.
    \end{equation*}
\end{proposition}

\begin{proof}[Proof of \cref{propo:lifumu:known:concentration:process:app}]
    We start by decomposing the process $(\rmZ_k^{\rmw})_{k\in\calK}$ again by its bounded and its unbounded part, i.e. $\rmZ_k^{\rmw} = \IWb_k + \IWu_k$, and obtain
    \begin{align}
        \IE_{Y}^N\brE{\max_{k\in\magenta{\calK}}\brK{\br{\abs{\rmZ_k^{\rmw}}^2-2\penbar_k^{\rmw}}_+}} & \leq 2\sum_{k\in\magenta{\calK}}\IE_f^N\brE{\br{\abs{\IWb_k}^2-\penbar_k^{\rmw}}_+}\nonumber\\
            &\hspace{1cm}+2\sum_{k\in\magenta{\calK}}\IE_Y^N\brE{\abs{\IWu_k}^2}\label{eq:lifumu:proof:concentration:decompos},
    \end{align}
    where we have used 
    \begin{equation*}
        \max\brK{(a+b)^2-2c,0} \leq \max\brK{2a^2-2c + 2b^2,0} \leq 2\max\brK{a^2-c,0} +2b^2,
    \end{equation*}
    which holds for $a,b,c\in\IR$.
    Set $\mathrm{p} \coloneqq 1+\frac{1}{\rmq-1}\in(1,2]$ and to shorten notation $\rma \coloneqq \rmq-1$.
    Recall $\mathfrak{Pen}_k^{\rmw}$ from \cref{lem:lifumu:app:concentration:process:bounded}.
    We then have
    \begin{align*}
        N\times \mathfrak{Pen}_k^{\rmw} &= \br{4\mathrm{p}\rmV_k+512k\nu_k^2\mathrm{p}^2\norm{\rmw\indk{k}}_{\IL^2}^2\frac{\ln(\Delta_k)}{N}} \ln(\Delta_k)\\
        & \leq 8\times \br{\rmV_k+2k\norm{\rmw\indk{k}}_{\IL^2}^2\frac{\Delta_k^{\frac{1}{\rma}}\ln(\Delta_k)}{N}} \ln(\Delta_k) = N\times\penbar_k^{\rmw}.
    \end{align*} 
    Thus, exploiting  \cref{lem:lifumu:app:concentration:process:bounded} as well as \cref{lifumu:app:varinace:process},
    we can upper bound \cref{eq:lifumu:proof:concentration:decompos} by 
    \begin{align*}
        &\IE_{Y}^N\brE{\max_{k\in\magenta{\calK}}\brK{\br{\abs{\rmZ_k^{\rmw}}^2-2\penbar_k^{\rmw}}_+}}\\
        & \hspace{1cm}  \leq 2\sum_{k\in\magenta{\calK}}\br{\IE_Y^N\brE{\br{\abs{\IWb_k}^2-\penbar_k^{\rmw}}_+}+\IE_Y^N\brE{\abs{\IWu_k}^2}}\\
        & \hspace{1cm}  \leq 16\sum_{k\in\magenta{\calK}}\frac{\Delta_k^{-\mathrm{p}}}{N}\norm{\rmw\indk{k}}_{\IL^2}^2\br{\ConstfY\Constfeps+\frac{256k\nu_k^2}{N}}\\
            & \hspace{2cm}+ 4\sum_{k\in\magenta{\calK}}\frac{k}{N}\norm{\indk{k}\rmw}_{\IL^2}^2\nu_k^{2(1-\rmp)}\IE\brE{Y^{2\rmp(c-1)}}.
    \end{align*}
    Now observe, \magenta{as $\calK\subseteq\kappa\IN$, we have for all $k\in\calK$}
    \begin{align*}
       4\frac{k}{N}\norm{\indk{k}\rmw}_{\IL^2}^2\nu_k^{2(1-\rmp)}\IE\brE{Y^{2\rmp(c-1)}} & =  128^{\rmp-1}\times4\times\frac{k}{N}\norm{\indk{k}\rmw}_{\IL^2}^2\Delta_k^{-1}\IE\brE{Y^{2\rmp(c-1)}} \\
       &\leq128^{\rmp}\times64^{-1}\times\IE\brE{Y^{2\rmp(c-1)}}\frac{\magenta{\kappa^2}}{k^2}N^{-1}.
    \end{align*}
    Similarly, 
    \begin{align*}
        &16\frac{\Delta_k^{-\mathrm{p}}}{N}\norm{\rmw\indk{k}}_{\IL^2}^2\br{\ConstfY\Constfeps+\frac{256k\nu_k^2}{N}}
          \leq 24\times \frac{1}{N}\frac{1}{k^2}\ConstfY\Constfeps.
    \end{align*}
    We deduce, 
    \begin{align*}
        &\IE_{Y}^N\brE{\max_{k\in\magenta{\calK}}\brK{\br{\abs{\rmZ_k^{\rmw}}^2-2\penbar_k^{\rmw}}_+}}\\
        &\hspace{1cm} \leq \sum_{k\in\magenta{\calK}}\br{24\times \frac{1}{N}\frac{1}{k^2}\ConstfY\Constfeps+128^{\rmp}\times64^{-1}\times\IE\brE{Y^{2\rmp(c-1)}}\frac{1}{N}\frac{1}{k^2}}\\
        &\hspace{1cm} \leq \frac{1}{N}\times \br{24\ConstfY\Constfeps+\frac{128^\rmp}{64}\IE\brE{Y^{2\rmp(c-1)}}}\times \sum_{k\in\magenta{\calK}}\frac{1}{k^2}\\
        &\hspace{1cm} \leq \frac{1}{N}\times \br{24\ConstfY\Constfeps+\frac{128^\rmp}{64}\IE\brE{Y^{2\rmp(c-1)}}}\times \frac{\pi^2}{6\magenta{\kappa^2}}\\
       &\hspace{1cm} \leq \frac{1}{N\magenta{\kappa^2}}\times \br{40\ConstfY\Constfeps+\frac{128^\rmp}{32}\IE\brE{Y^{2\rmp(c-1)}}},
    \end{align*}
    which completes the proof.
\end{proof}

\begin{proposition}[Estimated Penalty]\label{lifumu:unknown:app:estimatedpenality}
    \mbox{}\\
    For any $N\in\IN$ and $k\in\pRz$ define
    \begin{equation*}
        \rmVhat_k^{\rmw}\coloneqq \frac{1}{N}\sum_{i\in\nset{N}}\abs{\int_{\IR}\phib(Y_i,t)\rmw(t)\indk{k}(t)\rmd\lebesgue(t)}^2.
    \end{equation*}
    With the notations of \cref{propo:lifumu:known:concentration:process:app}, we define the estimated penalty
    \begin{equation*}
        \forall k\in\pRz:\;\penhat_k^{\rmw}\coloneqq8\br{\orange{2\rmVhat_k^{\rmw}}+\orange{3}\times 2k\norm{\rmw\indk{k}}_{\IL^2}^2\frac{\Delta_k^{\frac{1}{\rmp-1}}\ln(\Delta_k)}{N}} \frac{\ln(\Delta_k)}{N}.
    \end{equation*}
    Then the following holds.
    \begin{enumerate}
        \item For all $k\in\pRz$, we have $\IE_Y^N\brE{\rmVhat_k^{\rmw}} \leq \Constfeps\ConstfY \norm{\rmw\indk{k}}_{\IL^2}^2$.
        \item For all $k\in\pRz$, we have $\IE_Y^N\brE{\penhat_k^{\rmw}} \leq \frac{40}{N}\times \norm{\rmw\indk{k}}_{\IL^2}^2\Constfeps\ConstfY\rho_N(k|\rmw)$, where $\rho(k|\rmw)$ is defined in \cref{not:lifumu:datadriven:stuff}. 
        \item For all $\magenta{\calK\subseteq\kappa\IN}$, we have 
            \begin{equation*}
                \IE_{Y}^{N}\brE{\max_{k\in\magenta{\calK}}\brK{\br{\penhat_k^{\rmw}-\penbar_k^{\rmw}}_+}} \leq  \frac{0.18}{\magenta{\kappa^2}}\times\Constfeps\ConstfY\frac{1}{N}.
        \end{equation*}
    \end{enumerate}
\end{proposition}

\begin{proof}[Proof of \cref{lifumu:unknown:app:estimatedpenality}]
    \mbox{}
    \begin{enumerate}
        \item Fix $k\in\pRz$. The claim follows directly as $Y_1,\dots,Y_N$ are drawn i.i.d. and by observing
             $\rmV_k^{\rmw}\leq\Constfeps\ConstfY \norm{\rmw\indk{k}}_{\IL^2}^2$ as shown in the proof of \cref{lem:lifumu:app:concentration:process:bounded}.
        
        \item  Fix $k\in\pRz$. As $\IE_Y^N\brE{\rmVhat_k^{\rmw}} = \rmV_k^{\rmw}$, we have 
        \begin{align*}
            \IE_{Y}^N\brE{\penhat_k^{\rmw}} &= \IE_{Y}^N\brE{8\br{\orange{2\rmVhat_k^{\rmw}}+\orange{3}\times 2k\norm{\rmw\indk{k}}_{\IL^2}^2\frac{\Delta_k^{\frac{1}{\rmp-1}}\ln(\Delta_k)}{N}} \frac{\ln(\Delta_k)}{N}} \\
                & = 8\br{\orange{2\rmV_k^{\rmw}}+\orange{3}\times 2k\norm{\rmw\indk{k}}_{\IL^2}^2\frac{\Delta_k^{\frac{1}{\rmp-1}}\ln(\Delta_k)}{N}} \frac{\ln(\Delta_k)}{N}.
        \end{align*}
        Now $\rmV_k^{\rmw}\leq\Constfeps\ConstfY \norm{\rmw\indk{k}}_{\IL^2}^2$ implies
        \begin{align*}
            \IE_{Y}^N\brE{\penhat_k^{\rmw}} 
                & = \frac{8}{N}\times \norm{\rmw\indk{k}}_{\IL^2}^2\br{2\Constfeps\ConstfY+3\times 2k\frac{\br{\Delta_k^{\rmw}}^{\frac{1}{\rmp-1}}\ln(\Delta_k^{\rmw})}{N}} \ln(\Delta_k^{\rmw})\\
                & \leq \frac{40}{N}\times \norm{\rmw\indk{k}}_{\IL^2}^2\Constfeps\ConstfY\rho_N(k|\rmw),
        \end{align*}
        which shows the claim.
        
        \item We fix some $k\in\calK$. We then have 
        \begin{align*}
            \penbar_k^{\rmw}-\penhat_k^\rmw & =  8\br{\rmV_k+2k\norm{\rmw\indk{k}}_{\IL^2}^2\frac{\Delta_k^{\frac{1}{\rmp-1}}\ln(\Delta_k)}{N}} \frac{\ln(\Delta_k)}{N}\\
                                &\hspace{1cm} - 8\br{\orange{2\rmVhat_k}+\orange{3}\times 2k\norm{\rmw\indk{k}}_{\IL^2}^2\frac{\Delta_k^{\frac{1}{\rmp-1}}\ln(\Delta_k)}{N}} \frac{\ln(\Delta_k)}{N}\\
                             & = 8\br{\rmV_k-2\rmVhat_k-2\times2k\norm{\rmw\indk{k}}_{\IL^2}^2\frac{128\nu_k^2\ln(\Delta_k)}{N}} \frac{\ln(\Delta_k)}{N}\\
                             & = 16\br{\frac{1}{2}\rmV_k-\rmVhat_k-2k\norm{\rmw\indk{k}}_{\IL^2}^2\frac{128\nu_k^2\ln(\Delta_k)}{N}} \frac{\ln(\Delta_k)}{N}\\
                             & = 16\br{\rmV_k-\rmVhat_k-\frac{1}{2}\rmV_k-2k\norm{\rmw\indk{k}}_{\IL^2}^2\frac{128\nu_k^2\ln(\Delta_k)}{N}} \frac{\ln(\Delta_k)}{N}.
        \end{align*}
        Defining $\rmb_k^2\coloneqq2k\nu_k^2\norm{\rmw\indk{k}}_{\IL^2}^2$ as well as $\rma_k\coloneqq 128\frac{\ln(\Delta_k)}{N}$ we obtain
        \begin{equation*}
            \br{ \penbar_k^{\rmw}-\penhat_k^\rmw }_+ = 16\frac{\ln(\Delta_k)}{N}\br{\rmV_k-\rmVhat_k-\frac{1}{2}\rmV_k-\rma_k\rmb_k^2}_+ .
        \end{equation*}
        Now we define for all $i\in\nset{N}$
        \begin{equation*}
            Z_i \coloneqq - \abs{\int_{\IR}\phib(Y_i,t)\rmw(t)\indk{k}(t)\rmd\lebesgue(t)}^2
        \end{equation*}
        and observe that $ \IE_{Y}^{N}\brE{Z_i} = -\rmV_k^\rmw$, as well as
        \begin{align*}
            \rmV_k-\rmVhat_k 
                & = \frac{1}{N}\sum_{i\in\nset{N}}\br{Z_i-\IE[Z_i]}.
        \end{align*}
        Similar as before, we have 
        \begin{align*}
            \sup_{y\in\pRz}\brK{\abs{\int_{\IR}\phib(y,t)\rmw(t)\indk{k}(t)\rmd\lebesgue(t)}^2} &\leq \nu_k^2\abs{\int_{\IR}\abs{\rmw(t)}\indk{k}(t)\rmd\lebesgue(t)}^2\\
                & \leq 2k\times\nu_k^2\times\norm{\rmw\indk{k}}_{\IL^2}^2  = \rmb_k^2 \eqqcolon \rmb,
        \end{align*}
        i.e. $\abs{Z_i}\leq \rmb$ for all $i\in\nset{N}$. Further,
        \begin{align*}
            \forall i\in\nset{N}:\;\IE\brE{\abs{Z_i}^2} &= \IE\brE{\abs{\int_{\IR}\phib(Y_i,t)\rmw(t)\indk{k}(t)\rmd\lebesgue(t)}^4}\\
                &\leq \rmb_k^2\times \IE\brE{\abs{\int_{\IR}\phib(Y_i,t)\rmw(t)\indk{k}(t)\rmd\lebesgue(t)}^2} = \rmb_k^2\rmV_k \eqqcolon \rmv^2.
        \end{align*}
        Particularly, we have 
        \begin{align*}
            \forall m\in\IN_{\geq 2}:\; \frac{1}{N}\sum_{i\in\nset{N}} \IE\brE{\abs{Z_i}^m} \leq\rmv^2\rmb^{m-2} \leq \frac{m!}{2}\rmv^2\rmb^{m-2},
        \end{align*}
        and an application of \cref{lem:BernsteinInequality_general} yields
        \begin{align*}
            & \IE_{Y}^{N}\brE{\br{\rmV_k-\rmVhat_k-\frac{1}{2}\rmV_k-\rma_k\rmb_k^2}_+}  = \int_\IR \IP_{Y}^N\br{\br{\rmV_k-\rmVhat_k-\frac{1}{2}\rmV_k-\rma_k\rmb_k^2}_+ \geq t}\rmd\lebesgue(t)\\
                & \hspace{1cm} = \int_\IR \IP_{Y}^N\br{\rmV_k-\rmVhat_k\geq t+\frac{1}{2}\rmV_k+\rma_k\rmb_k^2}\rmd\lebesgue(t)\\
                & \hspace{1cm} = \int_\IR \IP_{Y}^N\br{\frac{1}{N}\sum_{i\in\nset{N}}\br{Z_i-\IE[Z_i]}\geq t+\frac{1}{2}\rmV_k+\rma_k\rmb_k^2}\rmd\lebesgue(t)\\
                & \hspace{1cm} \leq \int_{0}^\infty \max\brK{\exp\br{-\frac{N\br{t+\frac{1}{2}\rmV_k+\rma_k\rmb_k^2}^2}{4\rmv^2}},\exp\br{-\frac{N\br{t+\frac{1}{2}\rmV_k+\rma_k\rmb_k^2}}{4\rmb}}}\rmd\lebesgue(t).
        \end{align*}
        Now noting that $(\abs{a}+\abs{b})^2\geq a^2+b^2$ and $(a^2/2+b^2)^2\geq a^2b^2$ both hold true for all $a,b\in\IR$, we obtain
        \begin{equation*}
            \br{t+\frac{1}{2}\rmV_k+\rma_k\rmb_k^2}^2\geq t^2 +\br{\frac{1}{2}\rmV_k+\rma_k\rmb_k^2}^2\geq t^2+\rmV_k\rma_k\rmb_k^2
        \end{equation*}
        and
        \begin{equation*}
            \br{t+\frac{1}{2}\rmV_k+\rma_k\rmb_k^2} \geq \br{t+\rma_k\rmb_k^2}.
        \end{equation*}
        Conclusively, 
        \begin{align*}
            & \int_{0}^\infty \max\brK{\exp\br{-\frac{N\br{t+\frac{1}{2}\rmV_k+\rma_k\rmb_k^2}^2}{4\rmv^2}},\exp\br{-\frac{N\br{t+\frac{1}{2}\rmV_k+\rma_k\rmb_k^2}}{4\rmb}}}\rmd t\\
            &\hspace{1cm} \leq \int_{0}^\infty \max\brK{\exp\br{-\frac{Nt^2+N\rmV_k\rma_k\rmb_k^2}{4\rmb_k^2\rmV_k}},\exp\br{-\frac{Nt+N\rma_k\rmb_k^2}{4\rmb_k^2}}}\rmd t\\
            &\hspace{1cm} \leq \exp\br{-\frac{N\rma_k}{4}}\times\int_{0}^\infty \max\brK{\exp\br{-\frac{Nt^2}{4\rmb_k^2\rmV_k}},\exp\br{-\frac{Nt}{4\rmb_k^2}}}\rmd t.
        \end{align*}
        Exploiting $(a^2b^2)^\frac{1}{2}\leq a^2+b^2/4$, as well as
        \begin{equation*}
           \frac{N\rma_k}{4} = 32\ln(\Delta_k) ,\quad \int_0^\infty\exp(-x/a)\rmd x = a ,\quad \int_0^\infty \exp(-x^2/(4a^2))\rmd x = (\pi a^2)^\frac{1}{2},
        \end{equation*}
        the last display can be upper bounded by 
        \begin{align*}
            \Delta_k^{-32}\br{\br{\pi\rmb_k^2\rmV_kN^{-1}}^\frac{1}{2}+4\rmb_k^2N^{-1}} &\leq  \Delta_k^{-32}\br{\frac{1}{4}\pi\rmb_k^2N^{-1} + \rmV_k+4\rmb_k^2N^{-1}}\\
                & \leq  \Delta_k^{-32}\br{\rmV_k+5\rmb_k^2N^{-1}}.
        \end{align*}
        As $\rmV_k\leq\Constfeps\ConstfY \norm{\rmw\indk{k}}_{\IL^2}^2$ and recalling the definition of $\rmb_k^2$, we obtain 
        \begin{align*}
             \IE_{Y}^{N}\brE{\br{\rmV_k-\rmVhat_k-\frac{1}{2}\rmV_k-\rma_k\rmb_k^2}_+} &\leq \Delta_k^{-32}\br{\rmV_k+5\rmb_k^2N^{-1}}\\
                &\leq \Delta_k^{-32}\norm{\rmw\indk{k}}_{\IL^2}^2\br{\Constfeps\ConstfY +5\times 2k\nu_k^2\times N^{-1}}\\
                &\leq \frac{1}{2}\Delta_k^{-30}k^{-2}\br{\Constfeps\ConstfY+10\times \frac{1}{128}\times N^{-1}}\\
                &\leq \frac{69}{128}\Constfeps\ConstfY\Delta_k^{-30}k^{-2}.
        \end{align*}
        Overall, 
        \begin{align*}
             \IE_{Y}^{N}\brE{\max_{k\in\magenta{\calK}}\brK{\br{\penbar_k^{\rmw}-\penhat_k^{\rmw}}_+}} & \leq \sum_{k\in\magenta{\calK}} \IE_{Y}^{N}\brE{  \br{ \penbar_k^{\rmw}-\penhat_k }_+ }\\
                & = \frac{16}{N}\sum_{k\in\magenta{\calK}}\ln(\Delta_k) \IE_{Y}^{N}\brE{\br{\rmV_k-\rmVhat_k-\frac{1}{2}\rmV_k-\rma_k\rmb_k^2}_+}\\
                & \leq  \frac{16}{N}\times \frac{69}{128}\times\Constfeps\ConstfY \times \sum_{k\in\magenta{\calK}}\ln(\Delta_k)\Delta_k^{-30}k^{-2}.
        \end{align*}
        The function $h(x)\coloneqq\frac{\ln(x)}{x^A}$, $x,A\in\IR_{\geq 1}$, has a global maximum at $x_o=\exp(\frac{1}{A})$ with value $h(x_o)=\frac{1}{A}\rme^{-1}$.
        As $\Delta_k\geq 1$, the above display is upper bounded by 
        \begin{equation*}
            \frac{16}{N}\times \frac{69}{128}\times\Constfeps\ConstfY \times\frac{1}{30\rme}\sum_{k\in\magenta{\calK}}k^{-2} \leq   \frac{16}{N\magenta{\kappa^2}}\times \frac{69}{128}\times\Constfeps\ConstfY \times\frac{1}{30\rme} \frac{\pi^2}{6}\leq \frac{0.18}{\magenta{\kappa^2}}\times\Constfeps\ConstfY\frac{1}{N}.
        \end{equation*}
        This shows the claim.
    \end{enumerate}
\end{proof}

\begin{proof}[Proof of \cref{example:lifumu:known:datadriven:rates}]
		Recall the asymptotic order derived for
		$$\norm{\rmL\indkco{k}\rms^\dagger}_{\IL^2}^2\quad\and\quad\norm{\rmL\indk{k}\rmb}_{\IL^2}^2$$ as well as		
		$\maxrisklocal$ in \cref{example:lifumu:known:oracle:rates}. 
		We need to find the order of $\rho_N(k|\rmb)$ in terms of the tuning parameter $k$.
			\begin{itemize}
				\item (1):
					Let $p+b>-\frac{1}{2}$. We then have for all $k\in\IR_{\geq 1}$
					\begin{align*}
						\Delta_k^{\rmbos} &= \magenta{2(k/\kappa)^2\br{1\lor k\norm{\rmL\rmbos\mathds{1}_{[-k,k]}}_{\IL^2}^2} }\leq 2k^{3+2(p+b)+1},
					\end{align*}
					and thus 
					\begin{equation*}
						\ln\br{\Delta_k^{\rmbos}} \leq \const(p,b)\ln(k).
					\end{equation*}
					Thus, for sufficiently large $k$, we have 
					\begin{align*}
						\rho_N(k|\rmbos) &\leq \br{1\lor\br{ \frac{2k}{N}\cdot \br{\Delta_k^{\rmbos}}^{\frac{1}{\rmp-1}}\br{1\lor\ln\br{\Delta_k^{\rmbos}}}}} \cdot \br{1\lor\ln\br{\Delta_k^{\rmbos}}}\\
							&\leq \br{ \frac{2k}{N}\cdot \br{\Delta_k^{\rmbos}}^{\frac{1}{\rmp-1}}\ln\br{\Delta_k^{\rmbos}}} \cdot \ln\br{\Delta_k^{\rmbos}}\\
							&\leq \const(p,b,\rmq)\times \frac{1}{N}k^{\frac{3+2(p+b)+1}{\rmq-1}+1}\ln(k)^2\\
							& = \const(p,b,\rmq)\times \frac{1}{N}k^{\frac{\rmq+2+2(p+b)+1}{\rmq-1}}\ln(k)^2\\
							& \leq \const(p,b,\rmq)\times \frac{1}{N}k^{2s-2p+(2(p+b)+1)-1}\ln(k)^2.
					\end{align*}
					As a consequence, we balance
					\begin{equation*}
						k^{-2(s-p)+1}\quad\and\quad \frac{1}{N^2}k^{2s-2p+2(2(p+b)+1)-1}\ln(k)^2,
					\end{equation*}
					which means we need to choose $k$, such that 
					\begin{equation*}
						1 \asymp \frac{1}{N^2}k^{4(s+b)}\ln(k)^2
					\end{equation*}
					in the asymptotic regime $N\to\infty$.
                    Choosing $\kdiamond \asymp \br{\frac{N}{\ln(N)}}^{\frac{1}{2(s+b)}}$
					and computing $\ddrisklocal$ yields the claim.
	
				\item (2): 
					Let $p+b=-\frac{1}{2}$. As before, for $k\in\IR_{\geq 1}$ sufficiently large
					\begin{align*}
						\Delta_k^{\rmbos} &= \magenta{2(k/\kappa)^2\br{1\lor k\norm{\rmL\rmbos\mathds{1}_{[-k,k]}}_{\IL^2}^2} }\leq \const\times k^3\ln(k),
					\end{align*}
					and thus 
					\begin{equation*}
						\ln\br{\Delta_k^{\rmbos}} \leq \const\ln(k).
					\end{equation*}
					In other words,
					\begin{align*}
						\rho_N(k|\rmbos) 
							&\leq \br{ \frac{2k}{N}\cdot \br{\Delta_k^{\rmbos}}^{\frac{1}{\rmp-1}}\ln\br{\Delta_k^{\rmbos}}} \cdot \ln\br{\Delta_k^{\rmbos}}\\
							&\leq \const(p,b,\rmq)\times \frac{1}{N}k^{\frac{3}{\rmq-1}+1}\ln(k)^{2+\frac{1}{\rmq-1}}\\
							& = \const(p,b,\rmq)\times \frac{1}{N}k^{\frac{\rmq+3-1}{\rmq-1}}\ln(k)^{2+\frac{1}{\rmq-1}}\\
							& \leq \const(p,b,\rmq)\times \frac{1}{N}k^{2s-2p-1}\ln(k)^{2+\frac{1}{\rmq-1}}.
					\end{align*}
					We need to balance
					\begin{equation*}
						k^{-2(s-p)+1}\quad\and\quad \frac{1}{N^2}k^{2s-2p-1}\ln(k)^{4},
					\end{equation*}
					which means we need to choose $k$, such that 
					\begin{equation*}
						1 \asymp \frac{1}{N^2}k^{4(s-p)-2}\ln(k)^4
					\end{equation*}
					in the asymptotic regime $N\to\infty$.
                    Choosing $\kdiamond \asymp \br{\frac{N}{\ln(N)^2}}^{\frac{1}{2(s-p)-1}}$ and 
					computing $\ddrisklocal$ yields the claim.

				\item (3): Let $p+b<-\frac{1}{2}$. As before, for $k\in\IR_{\geq 1}$ sufficiently large
					\begin{align*}
						\Delta_k^{\rmbos} &= \magenta{2(k/\kappa)^2\br{1\lor k\norm{\rmL\rmbos\mathds{1}_{[-k,k]}}_{\IL^2}^2} }\leq \const\times k^3,
					\end{align*}
					and thus 
					\begin{equation*}
						\ln\br{\Delta_k^{\rmbos}} \leq \const\ln(k).
					\end{equation*}
					Consequently, for sufficiently large $k$, we have 
					\begin{align*}
						\rho_N(k|\rmbos) 
							&\leq \br{ \frac{2k}{N}\cdot \br{\Delta_k^{\rmbos}}^{\frac{1}{\rmp-1}}\ln\br{\Delta_k^{\rmbos}}} \cdot \ln\br{\Delta_k^{\rmbos}}\\
							&\leq \const(p,b,\rmq)\times \frac{1}{N}k^{\frac{3}{\rmq-1}+1}\ln(k)^{2}\\
							& = \const(p,b,\rmq)\times \frac{1}{N}k^{\frac{\rmq+3-1}{\rmq-1}}\ln(k)^{2}\\
							& \leq \const(p,b,\rmq)\times \frac{1}{N}k^{2s-2p-1}\ln(k)^{2}.
					\end{align*}
					We need to balance
					\begin{equation*}
						k^{-2(s-p)+1}\quad\and\quad \frac{1}{N^2}k^{2s-2p-1}\ln(k)^{2},
					\end{equation*}
					which means we need to choose $k$, such that 
					\begin{equation*}
						1 \asymp \frac{1}{N^2}k^{4(s-p)-2}\ln(k)^2
					\end{equation*}
					in the asymptotic regime $N\to\infty$.
                    Choosing $\kdiamond \asymp \br{\frac{N}{\ln(N)}}^{\frac{1}{2(s-p)-1}}$ and 
					computing $\ddrisklocal$ yields the claim.

				\item (4): For $p\in\IR$ and $b\in\pRzz$, we obtain as before
					
				\begin{align*}
					\Delta_k^{\rmbes} &= \magenta{2(k/\kappa)^2\br{1\lor k\norm{\rmL\rmbes\mathds{1}_{[-k,k]}}_{\IL^2}^2} }\leq \const(\rmL)\times k^{3+(2p+2b+1)_+}\exp(k^{2b}),
				\end{align*}
				and thus 
				\begin{equation*}
					\ln\br{\Delta_k^{\rmbes}} \leq \const k^{2b}.
				\end{equation*}
				Consequently, for sufficiently large $k$, we have 
					\begin{align*}
						\rho_N(k|\rmbes) 
							&\leq \const(c,p,b,\rmL)\frac{1}{N}\exp(\const(\rmq)k^{2b}).
					\end{align*}
					We need to balance
					\begin{equation*}
						k^{-2(s-p)+1}\quad\and\quad \frac{1}{N^2}\exp(\const(\rmq)k^{2b})
					\end{equation*}
					which means we need to choose $k$, such that 
					\begin{equation*}
						k^{2b} \asymp \ln(N)-\ln(k)
					\end{equation*}
					in the asymptotic regime $N\to\infty$.
                    Choosing  $\kdiamond \asymp \ln(N)^{\frac{1}{2b}}$ and
					computing $\ddrisklocal$ yields the claim.
				
				\item $(5)$ - $(7)$: Analogously, as in \cref{example:lifumu:known:oracle:rates} the bias term is the dominant one, 
					which determines the order of the optimal choice $k^\diamond$.
                    Thus, as before, we choose
					$$\kdiamond \asymp \ln(N)^{\frac{1}{2s}},$$
					which leads to 
					\begin{equation*}
						\rho_N(\kdiamond|\rmb) \asymp \ln\ln(N).
					\end{equation*}
					Computing $\ddrisklocal$ yields the claim.
			\end{itemize}
\end{proof}

\section{Proofs of \cref{sec:lifumu:unknown}}
\subsection{\texorpdfstring{Proofs of \cref{sec:lifumu:unknown:oracle}}{Proofs of Section 3.3}}\label{sec:lifumu:app:unknown:oracle}

We start with basic statements for the risk conditioned on the additional data.
\begin{lemma}\label{lem:lifumu:basicIdentity}
    \mbox{}\\
    In the statistical regime \assuRef{SRII}, grant \cref{ass:lifumu:A1} and \cref{lifumu:ass:A3}. 
    Let
    \begin{equation*}
       \forall t\in\IR:\; \MelXcheck(t)\coloneqq \MelY(t)\MelUhatrez(t)\indMl(t),\quad \forall k\in\pRz:\;\MelXcheck^k(t)\coloneqq \MelXcheck(t)\indk{k}(t).
    \end{equation*}
    Define
	\begin{equation*}
		\lincheck \coloneqq \int_{\IR}\rmL(-t)\MelXcheck(t)\rmd\lebesgue(t)\quad\text{and}\quad\lincheck_k \coloneqq \int_{\IR}\rmL(-t)\MelXcheck^k(t)\rmd\lebesgue(t).
	\end{equation*}
    Then, for all $N,M\in\IN$ and $k\in\pRz$, we have 
    \begin{align}
        &\IE_{f|\varepsilon}^{N,M}\brE{\abs{\linhat_k-\lin(f)}^2} = \IE_{f|\varepsilon}^{N,M}\brE{\abs{\linhat_k-\lincheck_k}^2}+\abs{\lincheck_k -\lin(f)}^2,\label{lifumu:unknown:biasvarianceI} \\
        &\IE_{f|\varepsilon}^{N,M}\brE{\abs{\linhat_k-\lincheck}^2} =  \IE_{f|\varepsilon}^{N,M}\brE{\abs{\linhat_k-\lincheck_k}^2}+\abs{\int_{[-k,k]^C}\rmL(t)\MelXcheck(t)\rmd\lebesgue(t)}^2, \label{lifumu:unknown:biasvarianceII} \\
        \text{and}\quad
        &\IE_{f|\varepsilon}^{N,M}\brE{\abs{\linhat_k-\lincheck_k}^2}\leq \frac{1}{N}\Constfeps\ConstfY\norm{\rmL\indk{k}\MelUhatrez\1_{\frakMl}}_{\IL^2}^2.\label{lifumu:unknown:biasvarianceIII}
    \end{align}
\end{lemma}

\begin{proof}[Proof of \cref{lem:lifumu:basicIdentity}]
	Fix $N,M\in\IN$ and $k\in\pRz$. We observe
	\begin{align*}
		 \IE_{f|\varepsilon}^{N,M}\brE{\linhat_k} &= \int_{[-k,k]}\rmL(-t)\MelY(t)\MelUhatrez(t)\1_{\frakMl}(t)\rmd\lebesgue(t) 
          = \lincheck_k.
	\end{align*}
	Thus, \cref{lifumu:unknown:biasvarianceI} follows directly from a squared bias-variance decomposition, i.e.
	\begin{align*}
		\IE_{f|\varepsilon}^{N,M}\brE{\abs{\linhat_k-\lin(f)}^2} & = \IE_{f|\varepsilon}^{N,M}\brE{\abs{\linhat_k- \IE_{f|\varepsilon}^{N,M}\brE{\linhat_k}+ \IE_{f|\varepsilon}^{N,M}\brE{\linhat_k}-\lin(f)}^2}\\
		& = \IE_{f|\varepsilon}^{N,M}\brE{\abs{\linhat_k-\lincheck_k}^2}+\abs{\lincheck_k -\lin(f)}^2.
	\end{align*}
	By replacing $\lin$ with $\lincheck$ in the last display, \cref{lifumu:unknown:biasvarianceII} follows analogously.
    For \cref{lifumu:unknown:biasvarianceIII} observe
	\begin{align*}
		&N\times  \IE_{f|\varepsilon}^{N,M}\brE{\abs{\linhat_k-\lincheck_k}^2} \\
        &\hspace{1cm}= N\times\IE_{f|\varepsilon}^{N,M}\brE{\abs{\int_{[-k,k]}\rmL(-t)\MelUhatrez(t)\1_{\frakMl}(t)\br{\MelYhat(t)-\MelY(t)}\rmd\lebesgue(t)}^2}\\
		& \hspace{1cm}\leq \Constfeps\ConstfY\norm{\rmL\MelUhatrez\1_{\frakMl}\indk{k}}_{\IL^2}^2,
	\end{align*}
	where we have used \cref{lem:MellinTrick}.
    This shows the claim.
\end{proof}

We need the following helping lemma.

\begin{lemma}\label{lem:lifumu:errorestimates:unknown:basicstuff}
    \mbox{}\\
    Grant \cref{ass:lifumu:A1} and \cref{lifumu:ass:A3}.
    Let $M\in\IN$.
    We then have for all $t\in\IR$ and $k\in\pRz$
    \begin{align}
        &\IE_{\varepsilon}^M\left[|\MelU(t)\MelUhatrez(t)|^2\1_{\frakMl}(t)\right]
          \leq 2\sigmabar_{\varepsilon,2}^2,\label{lem:2:p1}\\
        & \IE_\varepsilon^M\brE{\1_{\frakMl^C}(t)}\leq4\times\sigmabar_{\varepsilon,2}^2\br{1\lor\br{\abs{\MelU(t)}^2M}}^{-1},\label{lem:2:p2}\\
        \and\quad&\IE_{\varepsilon}^M\left[|\MelU(t)-\MelUhat(t)|^2|\MelUhat^\dagger(t)|^2\1_{\frakMl}(t)\right]
          \leq 2\const\times \frac{ \sigmabar_{\varepsilon,4}^2+ \sigmabar_{\varepsilon,2}^2}{1\lor(|\MelU(t)|^2M)}, \label{lem:2:p3}
      \end{align}
      where $\const\in\pRz$ is some numerical constant.
\end{lemma}

\begin{proof}[Proof of \cref{lem:lifumu:errorestimates:unknown:basicstuff}]
    First, observe that also $\MelUhat$ is an unbiased estimator of $\MelU$, that is 
  \begin{equation*}
    \forall t\in\IR:\; \IE_\varepsilon^M\brE{\MelUhat(t)} = \MelU(t).
  \end{equation*}
  Further, we have 
  \begin{align*}
    \forall t\in\IR:\;M\times\IE_\varepsilon^M\brE{\abs{\MelUhat(t)-\MelU(t)}^2} & \leq \IE\brE{\varepsilon^{2(c-1)}}\leq \sigmabar_{\varepsilon,2}^2,
  \end{align*}
  and after applying \cref{lem:Petro} with $\rmp = 4$ 
  \begin{equation*}
     \forall t\in\IR:\;M^2\times\IE_\varepsilon^M\brE{\abs{\MelUhat(t)-\MelU(t)}^4}  \leq \const\times \IE\brE{\varepsilon^{4(c-1)}} \leq \const\times\sigmabar_{\varepsilon,4}^2,
  \end{equation*}
  where $\const\in\IR_{\geq 1}$ is a numerical constant.
  We use those bounds without further reference.
  Let $t\in\IR$ be fixed.
  We start with \cref{lem:2:p1}.
  Using that
  $\abs{\MelUhat(t)\MelUhatrez(t)}^2\1_{\frakMl}(t)=\1_{\frakMl}(t)$
  and $\abs{\MelU(t)}^2\leq 2\abs{\MelUhat(t)}^2+2\abs{\MelU(t)-\MelUhat(t)}^2$, we
  have
  \begin{align*}
    \IE\left[\abs{\MelU(t)\MelUhatrez(t)}^2\1_{\frakMl}(t)\right]&\leq 2
    \IE\left[\abs{\MelU(t)-\MelUhat(t)}^2\abs{\MelUhatrez(t)}^2\1_{\frakMl}(t)
      + \1_{\frakMl}(t)\right] \\&\leq 2\left(M\times\IE\brE{\abs{\MelU(t)-\MelUhat(t)}^2}+1\right)\\&\leq
    2\left(\IE\brE{\varepsilon_1^{2(c-1)}}+1\right)
    = 2\sigmabar_{\varepsilon,2}^2,
  \end{align*}
  which shows \eqref{lem:2:p1}.
    Note \cref{lem:2:p2} is trivially satisfied if 
    \begin{equation*}
      1 \leq 4\sigmabar_{\varepsilon,2}^2\times\br{1\lor \br{\abs{\MelU(t)}^2M}}^{-1}.
    \end{equation*}
    Thus, assume $ \br{1\lor \br{\abs{\MelU(t)}^2M}} > 4\sigmabar_{\varepsilon,2}^2\geq 4$, 
    which implies $\frac{1}{M} < \frac{\abs{\MelU(t)}^2}{4}$.
    Hence, 
    \begin{equation*}
      \brK{\abs{\MelUhat(t)}^2 < M^{-1}} \subseteq \brK{\abs{\MelUhat(t)}^2 < \frac{\abs{\MelU(t)}^2}{4}}.
    \end{equation*}
    Observing that
    \begin{equation*}
      \abs{\MelUhat(t)}^2 < \frac{\abs{\MelU(t)}^2}{4}\quad\text{implies}\quad\abs{\MelUhat(t)}<2\abs{\MelUhat(t)} < \abs{\MelU(t)},
    \end{equation*}
    together with the reversed triangle inequality yields
    \begin{equation*}
      \abs{\MelUhat(t)-\MelU(t)} \geq \abs{\MelU(t)} - \abs{\MelUhat(t)} \geq \frac{1}{2}\abs{\MelU(t)}.
    \end{equation*}
    As a consequence, applying Markov's inequality, we obtain
    \begin{align*}
      \IE_\varepsilon^M\brE{\1_{\frakMl^C}(t)} & = \IP_\varepsilon^M\brE{\abs{\MelUhat(t)}^2 < M^{-1}}
        \leq \IP_\varepsilon^M\brE{ \abs{\MelUhat(t)-\MelU(t)} \geq \frac{1}{2}\abs{\MelU(t)}}\\
       & \leq 4\times\IE\brE{\abs{\MelUhat(t)-\MelU(t)}^2}\abs{\MelU(t)}^{-2}
        \leq 4\times\sigmabar_{\varepsilon,2}^2\br{\abs{\MelU(t)}^2M}^{-1}.
    \end{align*}
    As the left-hand side is a probability, we thus obtain
    \begin{equation*}
      \IE_\varepsilon^M\brE{\1_{\frakMl^C}(t)}\leq4\times\sigmabar_{\varepsilon,2}^2\br{1\lor\br{\abs{\MelU(t)}^2M}}^{-1}.
    \end{equation*}
  Lastly, consider \cref{lem:2:p3}. Evidently, we have on the
  one hand 
  \begin{equation*}
    \IE_{\varepsilon}^M\left[|\MelU(t)-\MelUhat(t)|^2|\MelUhat^\dagger(t)|^2\1_{\frakMl}(t)\right]\leq
    MM^{-1} \sigmabar_{\varepsilon,2}^2\leq  \sigmabar_{\varepsilon,2}^2,
  \end{equation*}
  while on the other hand  
  \begin{align*}
    &\IE_{\varepsilon}^M\left[|\MelU(t)-\MelUhat(t)|^2|\MelUhat^\dagger(t)|^2\1_{\frakMl}(t)\right]\\
    &\hspace{0.5cm}\leq
    2\IE_{\varepsilon}^M\left[\abs{\MelU(t)-\MelUhat(t)}^2\abs{\MelUhat^\dagger(t)}^2\1_{\frakMl}(t)\br{\frac{|\MelU(t)-\MelUhat(t)|^2}{|\MelU(t)|^2}+\frac{|\MelUhat(t)|^2}{|\MelU(t)|^2}}\right]\\
    &\hspace{0.5cm}\leq 2 \frac{M\IE\big[|\MelU(t)-\MelUhat(t)|^4\big]}{|\MelU(t)|^2}+2 \frac{\IE_{\varepsilon}^M\big[|\MelU(t)-\MelUhat(t)|^2\big]}{|\MelU(t)|^2}\\
    &\hspace{0.5cm}\leq 2\const\times ( \sigmabar_{\varepsilon,4}^2+ \sigmabar_{\varepsilon,2}^2) (|\MelU(t)|^2M)^{-1},
  \end{align*}
  where $\const\in\IR_{\geq 1}$ is a numerical constant after applying \cref{lem:Petro} with $\rmp=4$.
  Combining both bounds we obtain \eqref{lem:2:p3}, which completes the
  proof.
\end{proof}
  
We finally obtain the following Lemma.
\begin{lemma}\label{lem:lifumu:unknown:errorestimates}
    \mbox{}\\
	In the statistical regime \assuRef{SRII}, grant \cref{ass:lifumu:A1} and \cref{lifumu:ass:A3}. 
    Then for all $M\in\IN$ and $k\in\pRz$, we have
    \begin{align}
         &\IE_{\varepsilon}^M\brE{\norm{\rmL\indk{k}\MelUhatrez\1_{\frakMl}}_{\IL^2}^2} \leq 2\times\sigmabar_{\varepsilon,2}^2\norm{\rmL\MelUrez\indk{k}}_{\IL^2}^2,\label{lifumu:unknown:errorestimatesI}\\
         &\IE_{\varepsilon}^M\brE{\abs{\lincheck_k-\lin}^2} \leq 3\abs{\int_{[-k,k]^C}\rmL(t)\MelX(t)\rmd\lebesgue(t)}^2\notag\\
            &\hspace{3.05cm} +  \const_1\times \br{\sigmabar_{\varepsilon,2}^2+\sigmabar_{\varepsilon,4}^2}\norm{\rmL\indk{k}\MelX\br{1\lor (M\abs{\MelU}^2)}^{-\frac{1}{2}}}_{\IL^1}^2, \label{lifumu:unknown:errorestimatesII}\\
         &\IE_{\varepsilon}^M\brE{\abs{\int_\IR\rmL(t)\MelXcheck(t)\indkco{k}(t)\rmd\lebesgue(t)}^2} \leq 2\times\sigmabar_{\varepsilon,2}^2\norm{\rmL\MelX\indkco{k}}_{\IL^1}^2,\label{lifumu:unknown:errorestimatesIII}\\
        \and\quad &\IE_{\varepsilon}^M\brE{\abs{\lincheck-\lin}^2} \leq \const_2\times \br{\sigmabar_{\varepsilon,2}^2+\sigmabar_{\varepsilon,4}^2}\norm{\rmL\MelX\br{1\lor (M\abs{\MelU}^2)}^{-\frac{1}{2}}}_{\IL^1}^2 \label{lifumu:unknown:errorestimatesIV}
    \end{align}
    for some numerical constants $\const_1,\const_2\in\IR_{\geq 1}$.
\end{lemma}

\begin{proof}[Proof of \cref{lem:lifumu:unknown:errorestimates}]
    In the proof we will make use of \cref{lem:2:p1} - \cref{lem:2:p3}.
	We start with \cref{lifumu:unknown:errorestimatesI}. 
    To that end, fix in the following $k\in\pRz$.
    Fubini's theorem and an application of \cref{lem:2:p1} yields
	\begin{align*}
		&\IE_{\varepsilon}^M\brE{\norm{\rmL\indk{k}\MelUhatrez\1_{\frakMl}}_{\IL^2}^2}\\
        &\hspace{1cm}= \IE_\varepsilon^M\brE{\int_{\IR}\abs{\rmL(t)}^2\abs{\MelUrez(t)}^2\abs{\MelUhatrez(t)\MelU(t)}^2\1_{\frakMl}(t)\indk{k}(t)\rmd\lebesgue(t)}\\
		&\hspace{1cm} \leq 2\sigmabar_{\varepsilon,2}^2\times\int_{\IR}\abs{\rmL(t)}^2\abs{\MelUrez(t)}^2\indk{k}(t)\rmd\lebesgue(t)\\
		&\hspace{1cm} = 2\times\sigmabar_{\varepsilon,2}^2\norm{\rmL\MelUrez\indk{k}}_{\IL^2}^2.
	\end{align*}
	Similarly, applying the generalized Minkowski inequality and \cref{lem:2:p1} yields
	\begin{align*}
		& \IE_{\varepsilon}^M\brE{\abs{\int_\IR\rmL(t)\MelXcheck(t)\indkco{k}(t)\rmd\lebesgue(t)}^2}\\
		& \hspace{1cm}\leq \IE_{\varepsilon}^M\brE{\norm{\rmL\MelXcheck\indkco{k}}_{\IL^1}^2}\\
		& \hspace{1cm}= \IE_{\varepsilon}^M\brE{\abs{\int_\IR \abs{\rmL(t)}\abs{\MelX(t)}\indkco{k}(t)\abs{\MelU(t)\MelUhatrez(t)\1_{\frakMl}(t)}\rmd\lebesgue(t)}^2}\\
		&\hspace{1cm}\leq 2\times\sigmabar_{\varepsilon,2}^2\norm{\rmL\MelX\indkco{k}}_{\IL^1}^2.
	\end{align*}
	This shows \cref{lifumu:unknown:errorestimatesIII}. For \cref{lifumu:unknown:errorestimatesIV} observe that
	$\MelUhatrez(t)\MelUhat(t)\1_{\frakMl}(t) = \1_{\frakMl}(t)$ for all $t\in\IR$ and hence applying the generalized Minkowski inequality yields 
	\begin{align*}
		\frac{1}{2}\IE_{\varepsilon}^M\brE{\abs{\lincheck-\lin(f)}^2}
		& \leq \IE_{\varepsilon}^M\brE{\abs{\int_\IR\rmL(t)\MelX(t)\MelUhatrez(t)\1_{\frakMl}(t)\br{\MelU(t)-\MelUhat(t)}\rmd\lebesgue(t)}^2}\\
			&\hspace{1cm}+\IE_{\varepsilon}^M\brE{\abs{\int_\IR\rmL(t)\MelX(t)\1_{\frakMl^C}(t)\rmd\lebesgue(t)}^2}\\
		&\leq \abs{\int_\IR \abs{\rmL(t)}\abs{\MelX(t)}\IE_{\varepsilon}^M\brE{\abs{\MelUhatrez(t)}^2\1_{\frakMl}(t)\abs{\MelU(t)-\MelUhat(t)}^2}^\frac{1}{2}\rmd\lebesgue(t)}^2\\
			&\hspace{1cm}+\abs{\int_\IR \abs{\rmL(t)}\abs{\MelX(t)}\IE_{\varepsilon}^M\brE{\1_{\frakMl^C}(t)}^\frac{1}{2}\rmd\lebesgue(t)}^2\\
		&\leq 2\const\times\br{\sigmabar_{\varepsilon,2}^2+\sigmabar_{\varepsilon,4}^2}\times\abs{\int_\IR \abs{\rmL(t)}\abs{\MelX(t)}\br{1\lor (M\abs{\MelU(t)}^2)}^{-\frac{1}{2}}\rmd\lebesgue(t)}^2\\
		&\hspace{1cm}+4\sigmabar_{\varepsilon,2}^2\abs{\int_\IR \abs{\rmL(t)}\abs{\MelX(t)}\br{1\lor (M\abs{\MelU(t)}^2)}^{-\frac{1}{2}}\rmd\lebesgue(t)}^2\\
		& \leq \const\times \br{\sigmabar_{\varepsilon,2}^2+\sigmabar_{\varepsilon,4}^2}\norm{\rmL\MelX\br{1\lor (M\abs{\MelU}^2)}^{-\frac{1}{2}}}_{\IL^1}^2,
	\end{align*}
	where we have additionally used \cref{lem:2:p2} and \cref{lem:2:p3}, and where
    the numerical value of $\const$ changed from line to line.
	This shows \cref{lifumu:unknown:errorestimatesIV}.
	Similarly, we show \cref{lifumu:unknown:errorestimatesII} by observing that 
	\begin{align*}
		&\frac{1}{3}\IE_{\varepsilon}^M\brE{\abs{\lincheck_k-\lin}^2}- \abs{\int_{[-k,k]^C}\rmL(t)\MelX(t)\rmd\lebesgue(t)}^2\\
            &\hspace{1cm}  = \frac{1}{3}\IE_\varepsilon^M\brE{\abs{\int_\IR\rmL(t)\br{\MelXcheck^k(t)-\MelX(t)}\rmd\lebesgue(t)}^2}- \abs{\int_{[-k,k]^C}\rmL(t)\MelX(t)\rmd\lebesgue(t)}^2\\
		    &\hspace{1cm}  \leq \IE_\varepsilon^M\brE{\abs{\int_\IR\rmL(t)\br{\MelXcheck^k(t)-\MelX^k(t)\1_{\frakMl}(t)}\rmd\lebesgue(t)}^2}\\
                &\hspace{2cm}  +\IE_\varepsilon^M\brE{\abs{\int_\IR\rmL(t)\br{\MelX^k(t)\1_{\frakMl}(t)-\MelX^k(t)}\rmd\lebesgue(t)}^2}\\
                &\hspace{2cm} +\abs{\int_\IR\rmL(t)\br{\MelX^k(t)-\MelX(t)}\rmd\lebesgue(t)}^2- \abs{\int_{[-k,k]^C}\rmL(t)\MelX(t)\rmd\lebesgue(t)}^2\\
            &\hspace{1cm}  = \IE_\varepsilon^M\brE{\abs{\int_{[-k,k]}\rmL(t)\br{\MelXcheck(t)-\MelX(t)\1_{\frakMl}(t)}\rmd\lebesgue(t)}^2}\\
                &\hspace{2cm}  +\IE_\varepsilon^M\brE{\abs{\int_{[-k,k]}\rmL(t)\MelX(t)\1_{\frakMl^C}(t)\rmd\lebesgue(t)}^2}\\
		&\hspace{1cm}\leq \const\times \br{\sigmabar_{\varepsilon,2}^2+\sigmabar_{\varepsilon,4}^2}\norm{\rmL\indk{k}\MelX\br{1\lor (M\abs{\MelU}^2)}^{-\frac{1}{2}}}_{\IL^1}^2,
	\end{align*}
	using the same arguments as above. This shows the claim.
\end{proof}

\section{Probabilistic Inequalities}\label{sec:probabilisticInequalities}


\begin{lemma}[Petrov's Inequality]\label{lem:Petro}
  \mbox{}\\
  For $N\in\IN$, let $(Z_i)_{i\in\nset{N}}$ be a family of independent and zero-mean random variables.
  Let $\rmp\in\IR_{\geq 2}$. 
  Then there exists a constant $\const = \const(\rmp)\in\IR_{\geq 1}$, such that
  \begin{equation*}
    \IE\brE{\abs{\sum_{i\in\nset{N}}Z_i}^\rmp} \leq \const(\rmp)\times N^{\frac{\rmp}{2}-1}\sum_{i\in\nset{N}}\IE\brE{\abs{Z_i}^\rmp}.
  \end{equation*}
\end{lemma}
\begin{proof}[Proof of \cref{lem:Petro}]
    The claim is proven in \cite[Theorem 2.10]{Petrov_1995}.
\end{proof}

\begin{lemma}[Bernstein Inequality]\label{lem:BernsteinInequality_general}
    \mbox{}\\
    For $N\in\IN$, let $(Z_i)_{i\in\nset{N}}$ be a family of independent and real-valued random variables
    that satisfy
    \begin{equation}\label{eq:condition:general:Bernstein}
        \forall m\in\IN_{\geq 2}:\; \frac{1}{N}\sum_{i\in\nset{N}}\IE\brE{\abs{Z_i}^m} \leq \frac{m!}{2}\rmv^2\rmb^{m-2}
    \end{equation}
    for some finite constants $\rmb,\rmv\in\pRzz$. 
    We then have for all $x\in\pRzz$
    \begin{align}
    \IP\br{\abs{\frac{1}{N}\sum_{i\in\nset{N}}\br{Z_i-\IE[Z_i]}} \geq x} &\leq 2\exp\br{-Ny}\label{eq:Bern:1}\\
       &\leq2\exp\br{-\frac{Nx^2}{2\br{\rmv^2+\rmb x}}} \label{eq:Bern:2}\\
       &\leq2\max\brK{\exp\br{-\frac{Nx^2}{4\rmv^2}},\exp\br{-\frac{Nx}{4\rmb}}}, \label{eq:Bern:3}
    \end{align}
    where $y$ is defined by the equation $x = \rmv\sqrt{2y}+\rmb y$.
    Further, for any $K\in\IR_{\geq 1}$ we have 
    \begin{equation}\label{eq:Bern:4}
       \IE\brE{\br{\abs{\frac{1}{\sqrt{N}}\sum_{i\in\nset{N}}\br{Z_i-\IE[Z_i]}}^2-\br{4\rmv^2+32\rmb^2\frac{\ln(K)}{N}}\ln(K)}_+} \leq\frac{8}{K}\br{\rmv^2+\frac{16\rmb^2}{N}}.
    \end{equation}
\end{lemma}

\begin{proof}[Proof of \cref{lem:BernsteinInequality_general}]
    In \cite[Lemma 8]{BirgeMassart_1998} \cref{eq:Bern:1} and \cref{eq:Bern:2} are proven without the absolute value on the left-hand side and the factor 2
    on the right-hand side. The claim follows immediately by applying the results to $Z_i'\coloneqq -Z_i$ as well.
    \Cref{eq:Bern:3} then follows from \cref{eq:Bern:2} by distinguishing the cases $x\leq \frac{\rmv^2}{\rmb}$ and $x>\frac{\rmv^2}{\rmb}$. 
     To prove \cref{eq:Bern:4}, we observe that for any real valued random variable $X$ and $a\in\pRZz$
     \begin{align*}
         \IE\brE{\br{X^2-a}_+} 
            & \leq  \int_{\pRzz} \IP\br{\abs{X}\geq \sqrt{t+a}}\rmd \lebesgue(t).
     \end{align*}
     Thus, 
     \begin{align*}
         &\IE\brE{\br{\abs{\frac{1}{\sqrt{N}}\sum_{i\in\nset{N}}\br{Z_i-\IE[Z_i]}}^2-\br{4\rmv^2+32\rmb^2\frac{\ln(K)}{N}}\ln(K)}_+}\\
            &\hspace{1cm} \leq \int_{\pRzz} \IP\br{\abs{\frac{1}{N}\sum_{i\in\nset{N}}\br{Z_i-\IE[Z_i]}}\geq \sqrt{\br{4\rmv^2+32\rmb^2\frac{\ln(K)}{N}}\frac{\ln(K)}{N}+\frac{t}{N}}}\rmd \lebesgue(t)\\
            %
            & \hspace{1cm} \leq \frac{2}{K}\int_0^\infty \max\brK{\exp\br{-\frac{t}{4\rmv^2}},\exp\br{-\frac{\sqrt{Nt}}{4\sqrt{2}\rmb}}}  \rmd t,
     \end{align*}
    where the last step is due to \cref{eq:Bern:3}.
    Exploiting
     \begin{equation*}
        \sqrt{2(a+b)}\geq \sqrt{a} + \sqrt{b},\quad \int_0^\infty a\exp(-ax)\rmd x = 1 \quad \and\quad \int_0^\infty a^2\exp(-a\sqrt{x})\rmd x = 2,
     \end{equation*}
     which hold for all $a,b\in\pRz$,
     the last display is upper bounded by 
     \begin{equation*}
        \frac{2}{K}\br{4\rmv^2+\frac{64\rmb^2}{N}} = \frac{8}{K}\br{\rmv^2+\frac{16\rmb^2}{N}},
     \end{equation*}
     which shows the claim.
\end{proof}

\end{document}